\documentclass[11pt]{article}
\usepackage{amsthm,amsfonts,amssymb,amsmath}
\usepackage{graphicx}
\usepackage[english]{babel}
\usepackage{enumerate}
\usepackage{tikz}
\usepackage{xcolor}
\usepackage[normalem]{ulem}
\usepackage{cancel}
\usepackage[hidelinks]{hyperref}
\usepackage{authblk}
\usepackage[T1]{fontenc}

\newtheorem{theorem}{Theorem}
\newtheorem{lemma}[theorem]{Lemma}
\newtheorem{proposition}[theorem]{Proposition}
\newtheorem{corollary}[theorem]{Corollary}
\newtheorem{observation}[theorem]{Observation}
\newtheorem{definition}[theorem]{Definition} 
\newtheorem{conjecture}[theorem]{Conjecture}
\newtheorem{problem}[theorem]{Problem}

\DeclareMathOperator{\rad}{rad}
\DeclareMathOperator{\diam}{diam}
\DeclareMathOperator{\ecc}{ecc}
\DeclareMathOperator{\girth}{girth}

\begin{document}

\title{Global coalitions in graphs}

\author[1]{Marcin Anholcer\thanks{The author was supported by the National Science Centre, Poland (grant number: 2023/51/B/HS4/00829).}}
\author[2]{Nazli Besharati}
\author[3]{Doost Ali Mojdeh}
\author[3]{Mohammad Reza Samadzadeh}
\affil[1]{\scriptsize{}Pozna\'n University of Economics and Business, Institute of Informatics and Quantitative Economy}
\affil[ ]{Al.Niepodleg{\l}o\'sci 10, 61-875 Pozna\'n, Poland, \textit{marcin.anholcer@ue.poznan.pl}}
\affil[ ]{}
\affil[2]{\scriptsize{}Department of Mathematics, Faculty of Sciences}
\affil[ ]{Payame Noor University, Tehran, Iran, \textit{nbesharati@pnu.ac.ir}}
\affil[ ]{}
\affil[3]{\scriptsize{}Department of Mathematics, Faculty of Mathematical Sciences}
\affil[ ]{University of Mazandaran, Babolsar, Iran, \textit{damojdeh@umz.ac.ir, m.samadzadeh02@umail.umz.ac.ir}}
\affil[ ]{}

\date{}
\maketitle

\begin{abstract}
Let $G=(V, E)$ be a graph. A subset $S \subseteq V$ is called a global dominating set of $G$ if it is a dominating set of both $G$ and its complement $\overline{G}$. Two disjoint subsets  $V_1, V_2 \subseteq V$  form a global coalition if neither $V_1$ nor $V_2$ is a global dominating set,
but their union $V_1 \cup V_2$ is. A global coalition partition (abbreviated $gc$-partition) of $G$ is a vertex partition $\pi$ of $V(G)$ such that every class of $\pi$ forms a global coalition with some other class of $\pi$. The global coalition number $GC(G)$ is the maximum cardinality of a $gc$-partition of $G$.
In this paper, we initiate the study of these concepts. We prove that every graph of order at least two admits a $gc$-partition, give general lower bounds on $GC(G)$, and bound the number of global coalitions in which a single class of a $gc$-partition can be involved. We show that $GC(G)$ coincides with the coalition number $C(G)$ whenever $G$ is disconnected or has radius at least three, and consequently that $GC(G)=C(\overline{G})$ for every join of two graphs. We determine the global coalition number of complete multipartite graphs with at most three parts, wheels, fans, paths, cycles and trees of radius two, characterise the graphs $G$ with $GC(G)=|V(G)|$ in terms of antitwins, relate global coalitions to perfect coalitions, and study unicyclic graphs of radius at most two. Finally, we interpret $gc$-partitions in the language of fair division of indivisible goods and pose several open problems.\\
\textbf{2020 Mathematics Subject Classification:} 05C69, 91B32.\\
\textbf{JEL Classification:} C72, C78, D63.\\
\textbf{Keywords}: Global domination, coalition, global coalition partition, antitwins, trees, unicyclic graphs, fair division.
\end{abstract}


\section{Introduction and preliminary results}\label{sec-intro}
\subsection{Introduction}
Throughout this paper, all graphs are assumed to be simple, finite, and undirected.
We first recall the terminology used in the paper. For notions not defined here, we refer the reader to \cite{Bondy}.
Let $G$ be a graph with vertex set $V(G)$ and edge set $E(G)$, and let $n=|V(G)|$ be the order of $G$. We denote the degree of a vertex $v$ in $G$ by $\deg(v)$.
A vertex $v$ of degree $1$ is called a  \emph{leaf}, and its neighbour is called a \emph{support vertex}.
A vertex of degree $n-1$ is called a \emph{full vertex}, while a vertex of degree $0$ is called an
\emph{isolated vertex}.
The minimum and maximum degrees among the vertices of $G$ are denoted by $\delta(G)$ and $\Delta(G)$, respectively.
Two vertices are \emph{neighbours} if they are adjacent. The \emph{open neighborhood} of a vertex $v$ in $G$, denoted $N_G(v)$, is the set of neighbors of $v$ in $G$, while the \emph{closed neighborhood} of $v$ is the set $N_G[v] = N_G(v) \cup \{v\}$.
Given a subset $S \subseteq V(G)$, the open neighborhood of $S$ is defined as
$N_G(S) = \bigcup_{v \in S} N_G(v)$, and the closed neighborhood of  $S$ is defined as  $N_G[S] = \bigcup_{v \in S} N_G[v]$.
When the graph is clear from the context, we write $N(v)$, $N[v]$, $N(S)$ and $N[S]$.
For a set $S$ of vertices of $G$, the subgraph induced by $S$ is denoted by $G[S]$.
Let $X$ and $Y$  be two subsets of vertices in a graph $G$. The set of edges with one endpoint in $X$ and the other in $Y$
is denoted by $ [X, Y]$.
If every vertex in $X$ is adjacent to every vertex in  $Y$, then $[X,Y]$ is said to be full.
If there are no edges between the vertices of $X$ and $Y$, then $[X,Y]$ is said to be empty.
A vertex set $U$ is referred to as a \emph{singleton} if $\lvert U \rvert =1$, and as a \emph{non-singleton} if $\lvert U \rvert \geq 2$.
The double star  $S_{p,q}$ ($p,q\ge 1$) is a tree consisting of two adjacent support vertices: one adjacent to $p$ leaves, and the other to $q$ leaves.
We denote the path, cycle and complete graph of order $n$ by $P_n$, $C_n$ and $K_n$, respectively. The complete bipartite graph with partite sets of order $n$ and $m$ is denoted by $K_{n,m}$, and the complete $m$-partite graph with partite sets of orders $n_1,\dots,n_m$ by $K_{n_1,\dots,n_m}$.
The complete graph $K_3$ is called a \emph{triangle}, and a graph is \emph{triangle-free} if it has no $K_3$ as an induced subgraph.
Let $G$ be a graph of order $n$. The complement of $G$, denoted by $\overline{G}$, is defined as the graph with the same vertex set as $G$, i.e., $V(\overline{G}) = V(G)$, and edge set $E(\overline{G})=E(K_n)-E(G)$.
The \emph{join} $G_1+G_2$ of two vertex-disjoint graphs $G_1$ and $G_2$ is obtained from their disjoint union $G_1\cup G_2$ by adding all edges between $V(G_1)$ and $V(G_2)$. Note that $\overline{G_1+G_2}=\overline{G_1}\cup\overline{G_2}$. The \emph{wheel} $W_n$ ($n\ge 3$) and the \emph{fan} $F_n$ ($n\ge 2$) are the joins $K_1+C_n$ and $K_1+P_n$, respectively. The vertex of the $K_1$ is called the \emph{hub}.
In a connected graph $G$, the \emph{eccentricity}   of a vertex $v$, denoted  $\ecc(v)$, is defined as the greatest distance between $v$ and any other vertex in $G$.
The \emph{diameter} of the graph, denoted $\diam(G)$, is the maximum eccentricity among all vertices in $G$, while the \emph{radius},
denoted $\rad(G)$, is the minimum eccentricity. A vertex $v \in V(G)$ with $\ecc(v)=\rad(G)$ is called a \emph{central vertex} of $G$. The \emph{girth} of a graph containing a cycle is the length of its shortest cycle, and is denoted by $\girth(G)$.

A set $S \subseteq V(G)$ is a \emph{dominating set} of the graph $G$ if $N[S] = V(G)$. The minimum cardinality of a dominating set is the \emph{domination number} $\gamma(G)$ of $G$.
A partition of the vertices of $G$ into dominating sets is called a \emph{domatic partition}.
The \emph{domatic number} of $G$, denoted by $d(G)$, is the maximum number of sets in any domatic partition of $G$.
A dominating set $S$ is \emph{perfect} if every vertex of $V(G)\setminus S$ has exactly one neighbour in $S$. For further reading about domination, see e.g. \cite{haynes-hs, HHH-core}.

Domatic partitions describe situations in which the vertex set can be divided into groups, each of which is "strong enough" to dominate the graph. Haynes et al. \cite{haynes3} proposed a natural relaxation of this idea, motivated by coalitions of political parties: a group need not be strong on its own if it can become strong together with another group. Formally, let $G=(V, E)$ be a graph. Two disjoint sets $V_1, V_2 \subseteq V$ form a \emph{coalition} if neither of them is a dominating set of $G$, but their union $V_1 \cup V_2$ is. A \emph{coalition partition} of $G$ is a vertex partition $\pi$ such that each class $V_i \in \pi$ either is a singleton dominating set, or is not a dominating set on its own but forms a coalition with some class $V_j \in \pi$.
The \emph{coalition number} of $G$, denoted by $C(G)$, is the maximum cardinality of a coalition partition of $G$.

The coalition number was investigated by the same authors in \cite{haynes-ub,haynes1,haynes2,haynes4}, and by Bakhshesh, Henning and Pradhan \cite{BHP-trees, BHP-chains}, who, among other things, determined the coalition number of trees. Since then, several variants of the problem have been considered, like connected \cite{alikhani}, total \cite{ABG-total,barat}, independent \cite{ABGK,samadzade-m}, paired \cite{samadzade-mn}, perfect \cite{mojdeh-s}, double \cite{HeMo-24, HeMo-25}, and $k$-coalitions \cite{jafari}.

Another research direction considers sets that simultaneously dominate several graphs on the same vertex set. The most popular instance is \emph{global domination}, introduced independently by Sampathkumar \cite{sampathkumar} and by Brigham and Dutton \cite{brigham} (the latter in a more general setting of factor domination): a set $S \subseteq V(G)$ is a \emph{global dominating set} of $G$ if it is a dominating set of both $G$ and $\overline{G}$, and the \emph{global domination number} $\gamma_g(G)$ is the minimum cardinality of a global dominating set. See \cite{BCD-survey,desormeaux} for further reading. The partition counterpart of global domination, the \emph{global domatic number}, was introduced by Zelinka \cite{zelinka} (as \emph{complementarily domatic number}) and, independently, by Sampathkumar \cite{sampathkumar}.

In this paper, we combine the two concepts described above: we require that the union of two allied classes dominate not only $G$ but also $\overline{G}$. The resulting concept of a \emph{global coalition} is not only a formal variant. As we show in Section~\ref{sec-gl-co}, it differs from the classical coalition only for graphs $G$ such that both $G$ and $\overline{G}$ are connected and have radius at most two, a class that includes many dense and many "star-like" graphs. For such graphs, the requirement of dominating the complement interacts with domination in $G$ in a non-trivial way, which explains the behaviour of the parameter on the graph classes studied below.

\begin{definition}\label{def-gc}
	Let $G=(V, E)$ be a graph. Two disjoint subsets $V_i, V_j \subseteq V$  are said to form a \emph{global coalition} if neither $V_i$ nor $V_j$ is a global dominating set of $G$, but their union $V_i \cup V_j$ is a global dominating set of $G$. In this case, we also say that $V_i$ and $V_j$ are \emph{$gc$-partners}.
\end{definition}

\begin{definition}
	\label{def-gc-p}
	A \emph{global coalition partition}, abbreviated \emph{$gc$-partition}, of a graph $G$ is a vertex partition $\pi$ such that for every class $V_i \in \pi$, $V_i$ is not a global dominating set of $G$, but there exists another class $V_j \in \pi$ such that $V_i$ and $V_j$ form a global coalition.
	The \emph{global coalition number} of $G$, denoted by $GC(G)$, is the maximum cardinality of a $gc$-partition of $G$. A $gc$-partition of $G$ of cardinality $GC(G)$ is called a $GC(G)$-partition.
\end{definition}

Note that, in contrast with coalition partitions, we do not consider singleton global dominating classes. This is a natural assumption, since a singleton $\{v\}$ is a global dominating set only if $v$ is a full vertex of both $G$ and $\overline{G}$, which is impossible when $n\ge 2$. On the other hand, the graph $K_1$ has no $gc$-partition, so in the remainder of the paper we assume that $n\ge 2$.

The paper is organised as follows. In the remainder of this section, we determine $GC(G)$ for complete graphs, complete bipartite and tripartite graphs, wheels and fans, and we bound the number of global coalitions in which a single class of a $gc$-partition can participate. In Section~\ref{sec-gl-exist} we show that every graph of order at least two admits a $gc$-partition and derive two general lower bounds on $GC(G)$. Section~\ref{sec-gl-co} is devoted to the relation between $GC(G)$ and $C(G)$: we prove that the two parameters coincide whenever $G$ is disconnected or $\rad(G)\ge 3$, which yields $GC(G)=C(\overline{G})$ for all joins, and we determine $GC$ for paths, cycles and trees of radius two. In Section~\ref{sec-gl-perf} we characterise the graphs with $GC(G)=n$ in terms of antitwins and relate global coalitions to perfect coalitions. Section~\ref{sec-gl-uni} deals with unicyclic graphs of radius at most two. In Section~\ref{sec-fair}, we discuss a connection between $gc$-partitions and fair division of indivisible goods. In Section~\ref{sec-concl} we present some open problems.
\subsection{Preliminary  results}\label{sub-gl-preli}
For every graph $G$, $\overline{\overline{G}}=G$, so a set is a global dominating set of $G$ if and only if it is a global dominating set of $\overline{G}$. As a direct implication of  
Definition~\ref{def-gc-p}, we have the following result.
\begin{observation}\label{obs-compl}
	For any graph $G$, we have $GC(G)= GC(\overline{G})$.
\end{observation}
The next two observations determine the global coalition number of complete graphs and complete bipartite graphs.
\begin{observation}\label{obs-kn}
	For any complete graph $K_n$ ($n\geq 2$), we have $GC(K_n)=2$.
\end{observation}
\begin{proof}
Since $\overline{K_n}$ has no edges, the only global dominating set of $K_n$ is $V(K_n)$. Hence every $gc$-partition consists of exactly two classes, and $\{\{v\},V(K_n)\setminus\{v\}\}$ is one of them.
\end{proof}
\begin{observation} \label{obs-kbp}
	For any complete bipartite graph $K_{n,m}$ ($n,m \geq 1$), we have $GC(K_{n,m})=n+m$.
\end{observation}
\begin{proof}
Let $X$ and $Y$ be the partite sets. A set $S$ dominates $\overline{K_{n,m}}=K_n\cup K_m$ if and only if it meets both $X$ and $Y$, and such a set also dominates $K_{n,m}$. Thus $\{x\}$ and $\{y\}$ are $gc$-partners for every $x\in X$ and $y\in Y$, and the partition of $V(K_{n,m})$ into one-vertex sets is a $gc$-partition of $K_{n,m}$.
\end{proof}

The argument above extends to complete multipartite graphs: a set $S$ is a global dominating set of $G=K_{n_1,\dots,n_m}$ ($m\ge 2$) if and only if $S$ meets every partite set, since $\overline{G}$ is the disjoint union of the cliques induced by the partite sets. We call such a set a \emph{transversal}. The following proposition gives a general lower bound. Note that $n_1\ge\dots\ge n_m$, so the condition $n_m\ge t$ means in particular that every partite set has at least $t$ vertices.

\begin{proposition}\label{prop-multi}
	Let $G= K_{n_1,n_2,\dots ,n_m}$ be a complete $m$-partite graph with partite sets $V_1,V_2, \dots ,V_m$, where $n_1\ge n_2\ge \cdots \ge n_m$ and $m\ge 2$. Then
	\[
	GC(G)\ \ge\ \max\Bigl\{\, t+\sum_{i=1}^{t}(n_i-t+1)\ :\ 1\le t\le m,\  n_m\ge t \Bigr\}.
	\]
	 In particular, $GC(G)\ge n_1+1$, and $GC(G)\ge n_1+n_2$ whenever $n_m\ge 2$.
\end{proposition}
\begin{proof}

Select $t$ satisfying the conditions of the statement, so that $n_j\ge n_m\ge t$ for every $j$. For every $i\le t$ choose a set $H_i\subseteq V\setminus V_i$ meeting every partite set $V_j$, $j\ne i$, in such a way that the sets $H_1,\dots,H_t$ are pairwise disjoint and $\bigcup_{j>t}V_j\subseteq \bigcup_{i\le t}H_i$. This is possible, since every $V_j$ with $j>t$ has at least $t$ vertices, so it can be split among $H_1,\dots, H_t$ with each $H_i$ receiving at least one vertex, and every $V_j$ with $j\le t$ has at least $t-1$ vertices (indeed, at least $t$), one for each $H_i$ with $i\ne j$. The remaining vertices, namely $n_j-(t-1)$ vertices of each $V_j$ with $j\le t$, form one-vertex classes. The $t$ sets $H_i$ together with these singletons partition $V(G)$ into $t+\sum_{i\le t}(n_i-t+1)$ classes. None of the classes is a transversal. If $t\ge 2$, then $H_i\cup H_j$ is a transversal for $i\ne j$. If $t=1$, then $n_1\ge 1$ singletons are present. Finally, each singleton $\{v\}$ with $v\in V_i$ forms a global coalition with $H_i$. Thus, the partition is a $gc$-partition. The values $t=1$ and $t=2$ give the two bounds stated at the end of the Proposition.
\end{proof}

For $t=m$, the construction in the proof works if $n_m\ge m-1$, since there are no partite sets to split. We do not include this case in the statement, because it never increases the maximum. If $n_m=m-1$, then $t=m-1$ is admissible and gives the value $n-(m-1)(m-2)\geq n-m(m-2)$, which is the value obtained for $t=m$. For three partite sets, we can determine the exact value, which shows that the lower bound from Proposition~\ref{prop-multi} is sharp.

\begin{theorem}\label{thm-tripartite}
Let $G=K_{n_1,n_2,n_3}$ with $n_1\ge n_2\ge n_3\ge 1$. Then
\[
GC(G)=\begin{cases}
n_1+1, & \text{if } n_3=1,\\
n_1+n_2, & \text{if } 2\le n_3\le 3,\\
n_1+n_2+n_3-3, & \text{if } n_3\ge 3.
\end{cases}
\]
\end{theorem}
\begin{proof}
The lower bounds follow from Proposition~\ref{prop-multi} with $t=1$, $t=2$ and $t=3$, respectively. For the upper bounds, let $\pi$ be a $gc$-partition of $G$ and, for a class $A\in\pi$, let $T(A)\subseteq\{1,2,3\}$ be the set of indices $i$ with $A\cap V_i\ne\emptyset$. No class is a transversal, so $|T(A)|\in\{1,2\}$. We call $A$ a \emph{single} or a \emph{pair} respectively. Two classes are $gc$-partners if and only if the union of their index sets is $\{1,2,3\}$. Let $x_i$ be the number of singles $A$ with $T(A)=\{i\}$, let $p$ be the number of pairs, and let $c_i$ be the number of classes meeting $V_i$, so that $c_i\le n_i$. Counting incidences between classes and partite sets gives $c_1+c_2+c_3=\sum_i x_i+2p$, so
\begin{equation}\label{eq-tri}
|\pi| = \sum_i x_i + p = c_1+c_2+c_3-p \le n_1+n_2+n_3-p.
\end{equation}
Let $I=\{i : x_i>0\}$. A single with index set $\{i\}$ needs a partner, which must be a pair with index set $\{1,2,3\}\setminus\{i\}$. Thus $p\ge |I|$, and every $V_j$ with $j\in I$ meets at least $|I|-1$ pairs (the ones that are complementary to the other elements of $I$).

If $|I|=3$, then \eqref{eq-tri} gives $|\pi|\le n_1+n_2+n_3-3$. Moreover, each $V_i$ meets at least two pairs, so $n_3\ge 2$, and if $n_3=2$ then $x_3=0$, a contradiction. Hence $n_3\ge 3$ in this case.
If $I=\{i,j\}$ with $l$ the third index, then, writing $x_{ij}$ for the number of pairs with index set $\{i,j\}$ and similarly for the other pairs, we have $x_i\le n_i-x_{ij}-x_{il}$ and $x_j\le n_j-x_{ij}-x_{jl}$, whence
$|\pi|=x_i+x_j+x_{ij}+x_{il}+x_{jl}\le n_i+n_j-x_{ij}\le n_1+n_2$. If moreover $n_3=1$, then the pairs complementary to $i$ and to $j$ both meet $V_3$ unless $3\in\{i,j\}$. In the former case $x_{il}+x_{jl}\le 1$ contradicts $x_{il},x_{jl}\ge 1$, and in the latter case, say $j=3$, we get $|\pi|\le n_i+n_3-x_{i3}\le n_1+1$.
If $I=\{i\}$, then $|\pi|=x_i+p\le (n_i-x_{ij}-x_{il})+x_{ij}+x_{il}+x_{jl}=n_i+x_{jl}\le n_i+\min\{n_j,n_l\}$, which is at most $n_1+n_2$, and at most $n_1+1$ when $n_3=1$.
Finally, if $I=\emptyset$, then all classes are pairs and $|\pi|=p=(c_1+c_2+c_3)/2\le n/2$, which is at most $n_1+n_2$, and at most $n_1+1$ when $n_3=1$.
By taking the cases together, $GC(G)\le n_1+1$ if $n_3=1$, and $GC(G)\le\max\{n_1+n_2,\,n_1+n_2+n_3-3\}$ if $n_3\ge 2$, where the second term applies only when $n_3\ge 3$.
\end{proof}
	We conjecture that the lower bound of Proposition~\ref{prop-multi} is always tight (Conjecture~\ref{conj-multi}).

Next, we determine the global coalition number of two other classes of graphs, namely wheel graphs $W_n$ and fan graphs $F_n$. Both contain a full vertex $u$ (the hub) adjacent to all vertices of an induced subgraph: a cycle $C_n$ for $W_n$ (with $n\geq 3$) and a path $P_n$ for $F_n$ (with $n \geq 2$).
The following simple observation describes the global dominating sets of any graph with a full vertex.
\begin{lemma}\label{lem-full}
Let $G=K_1+H$, where $u$ is the vertex of $K_1$ and $H$ has order at least two. A set $S\subseteq V(G)$ is a global dominating set of $G$ if and only if $u\in S$ and $S\setminus\{u\}$ is a dominating set of $\overline{H}$.
\end{lemma}
\begin{proof}
The vertex $u$ is isolated in $\overline{G}=K_1\cup\overline{H}$, so every global dominating set contains $u$, and then $S$ dominates $G$ automatically. Since $u$ is isolated in $\overline{G}$, the set $S$ dominates $\overline{G}$ if and only if $S\setminus\{u\}$ dominates $\overline{H}$.
\end{proof}

For $W_n$, it is well known that $W_3 \simeq
K_4$ and  $\gamma_g(W_3)=4$.  If $n\ge 4$ and $V(C_n)= \{v_1,v_2,\cdots, v_n\}$, then $\gamma_g(W_n)=3$ and any set  $\{u,v_i,v_{j}\}$ with $j\notin\{i, i+2, i-2\}$ (indices modulo $n$), is a $\gamma_g$-set of $W_n$. Indeed, by Lemma~\ref{lem-full}, $\{u\}\cup T$ with $T\subseteq V(C_n)$ is a global dominating set of $W_n$ if and only if $T$ dominates $\overline{C_n}$, i.e., if and only if for every $i$ either $v_i\in T$ or $T\not\subseteq\{v_{i-1},v_{i+1}\}$. For $n\ge 4$ this holds exactly when $|T|\ge 3$, or $|T|=2$ and $T\ne\{v_{i-1},v_{i+1}\}$ for all $i$.

Based on these observations, we now establish the following result.

\begin{proposition}\label{prop-wheel}
	For the wheel $W_n$, $n\ge 3$, we have $GC(W_n)=n-1$.
\end{proposition}
\begin{proof}
For $n=3$ the claim follows from Observation~\ref{obs-kn}, so let $n\ge 4$. Let $V(W_n)=\{u,v_1,v_2,\cdots, v_n\}$, and let $\pi$ be a $gc$-partition of $W_n$.  Since the hub vertex $u$ belongs to every global dominating set of $W_n$,  every class $A \in \pi$ not containing $u$ is a $gc$-partner of the class $U$ containing $u$. Since $U$ is not a global dominating set, the description of the global dominating sets given above shows that  $U=\{u\}$ or $U=\{u,v_i\}$ or $U=\{u,v_i,v_{i+2}\}$ for some $i$ (indices modulo $n$). If $U=\{u\}$, then each $gc$-partner of $U$ must contain at least two vertices from the cycle $C_n$. In this case,  $U$ has at most $\lfloor\frac{n}{2}\rfloor$ $gc$-partners, and $|\pi|\le \lfloor n/2\rfloor+1\le n-1$.

If $U=\{u,v_i\}$, then a class $A\subseteq V(C_n)\setminus\{v_i\}$ is a $gc$-partner of $U$ unless $A=\{v_{i-2}\}$ or $A=\{v_{i+2}\}$. Hence the vertices $v_{i-2}$ and $v_{i+2}$ lie in non-singleton classes, and the $n-1$ vertices of $V(C_n)\setminus\{v_i\}$ are covered by at most $n-2$ classes, so $|\pi|\le n-1$. 

Finally, if $U=\{u,v_i,v_{i+2}\}$ , then the remaining $n-2$ vertices are covered by at most $n-2$ classes and again $|\pi|\le n-1$.
On the other hand, $\{u,v_1,v_3\}$ together with the singletons of the remaining $n-2$ vertices is a $gc$-partition of $W_n$, because every one-vertex class is a $gc$-partner of $\{u,v_1,v_3\}$. Hence $GC(W_n)=n-1$.
\end{proof}

For the fan graph $F_n$, it is well known that $F_2 \simeq K_3$ and  $\gamma_g(F_2)=3$.
Assume $n\ge 3$ and let $V(P_n)= \{v_1,v_2,\cdots, v_n\}$. Then $\gamma_g(F_n)=3$ and any set  $\{u,v_i,v_{j}\}$ with $|i-j|\ne 2$
is a $\gamma_g$-set of $F_n$. More precisely, by Lemma~\ref{lem-full}, a set $\{u\}\cup T$ with $T\subseteq V(P_n)$ is a global dominating set of $F_n$ if and only if $T$ dominates $\overline{P_n}$, which for $n\ge 2$ happens exactly when $|T|\ge 3$, or $|T|=2$ and $T\ne\{v_{i-1},v_{i+1}\}$ for all $2\le i\le n-1$.

The first version of this paper claimed that $GC(F_n)=n-1$ for all $n\ge 2$. This is false for $n\in\{2,3\}$: $GC(F_2)=GC(K_3)=2$ by Observation~\ref{obs-kn}, and $\{\{u,v_2\},\{v_1\},\{v_3\}\}$ is a $gc$-partition of $F_3$, so $GC(F_3)=3$. The reason is that for $n\le 3$ the vertex $v_2$ has no vertex at distance two on the path. The correct statement is the following.

\begin{proposition}\label{prop-fan}
	For the fan $F_n$ we have $GC(F_2)=2$, $GC(F_3)=3$ and $GC(F_n)=n-1$ for $n\ge 4$.
\end{proposition}
\begin{proof}
The values for $n\le 3$ were proved above. Let $n\ge 4$, let $\pi$ be a $gc$-partition of $F_n$, and let $U$ be the class containing $u$. As in the proof of Proposition~\ref{prop-wheel}, every other class is a $gc$-partner of $U$, and $U\in\{\{u\},\{u,v_i\},\{u,v_i,v_{i+2}\}\}$. If $U=\{u\}$, every partner has at least two vertices and $|\pi|\le \lfloor n/2\rfloor+1\le n-1$. If $U=\{u,v_i\}$, then a class $A$ is a partner of $U$ unless $A=\{v_{i-2}\}$ or $A=\{v_{i+2}\}$. Since $n\ge 4$, at least one of $v_{i-2},v_{i+2}$ exists and lies in a non-singleton class, so $|\pi|\le 1+(n-2)=n-1$. If $U=\{u,v_i,v_{i+2}\}$, then $|\pi|\le 1+(n-2)=n-1$. Conversely, $\{u,v_1,v_3\}$ together with the one-vertex classes consisting of the remaining vertices is a $gc$-partition of size $n-1$.
\end{proof}
Wheels, fans and complete multipartite graphs are joins. In Section~\ref{sec-gl-co} (Corollary~\ref{cor-join}), we show that for every join $G$ the parameter $GC(G)$ equals the coalition number of the disconnected graph $\overline{G}$. The propositions above may thus also be read as results about $C(K_1\cup\overline{C_n})$, $C(K_1\cup\overline{P_n})$ and $C(K_{n_1}\cup\dots\cup K_{n_m})$.

Our final result in this subsection provides an upper bound on the number of global coalitions in which each member of a $gc$-partition can participate. The bound below is slightly stronger than the one in the first version of this paper, which was $\max\{\Delta(G)+1,\min\{n-|A|,n-\delta(G)\}\}$, and it is symmetric with respect to $G$ and $\overline{G}$, in accordance with Observation~\ref{obs-compl}. Note that $\Delta(\overline{G})+1=n-\delta(G)$.
\begin{theorem}
	\label{gc-bound}
	Let $G$ be a graph of order $n$ with minimum degree $\delta (G)$ and maximum degree $\Delta(G)$, and let $\pi$ be a $gc$-partition of $G$.
	If $A \in \pi$, then the set $A$  belongs to at most $$\min\bigl\{\, n-|A|,\ \max\{\Delta(G)+1,\ n-\delta(G)\}\bigr\}$$
	global coalitions. Moreover, this bound is sharp.
\end{theorem}
\begin{proof}
Let $t$ denote the number of sets in $\pi$  that form a global coalition with $A$. Since the $gc$-partners of $A$ are pairwise disjoint non-empty subsets of $V(G)\setminus A$, we have $t\le n-|A|$.
Since $A \in \pi$, it follows that $A$ is not a global dominating set of $G$.
Assume first that $A$ does not dominate $G$. Let $x \in V(G)$ be a vertex that is not dominated by $A$.
Every $gc$-partner of $A$ must dominate $x$, that is, it must contain at least one vertex of $N[x]$.
Thus,  $t \leq |N[x]|\le \Delta(G)+1$. Assume now that $A$ dominates $G$.
It follows that $A$ does not dominate $\overline{G}$. Let $y \in V(\overline{G})$ be a vertex that is not dominated by $A$ in $\overline{G}$. Now every $gc$-partner of $A$ must  contain  at least one vertex of $N_{\overline{G}} [y]$, so $t\le |N_{\overline{G}}[y]|\le \Delta(\overline{G})+1=n-\delta(G)$. In both cases $t\le\max\{\Delta(G)+1,n-\delta(G)\}$, and the bound follows.

To see the sharpness, consider the family $G_k$ ($k \geq 2$) of graphs with  vertex set $V(G_k)=\lbrace u, v_i, w_i : 1 \leq i \leq k \rbrace$ and  edge set $E(G_k)= \lbrace uv_i : 1 \leq i \leq k\rbrace \cup  \lbrace v_i w_i : 1 \leq i \leq k \rbrace \cup \lbrace w_i w_{i+1}: 1\leq i\leq k\ (\rm{mod}\ k)\rbrace$ (for $k=2$ the last set consists of the single edge $w_1w_2$).  Figure \ref{fig-G4} illustrates such a graph for $k=4$. Observe that $n=2k+1$, $\Delta (G_k)=k$ and $\delta(G_k)=2$. Now consider the $gc$-partition $\lbrace \lbrace u\rbrace, \lbrace v_i : 1 \leq i \leq k\rbrace , \lbrace w_1\rbrace ,\lbrace w_2 \rbrace ,\dots ,\lbrace w_k\rbrace \rbrace$, where the  set  $A=\lbrace v_i:  1 \leq i \leq k\rbrace$ forms a global coalition with all $k+1$ other members of the partition, while the bound equals $\min\{n-|A|,\max\{k+1,2k-1\}\}=k+1$. This completes the proof.
\end{proof}
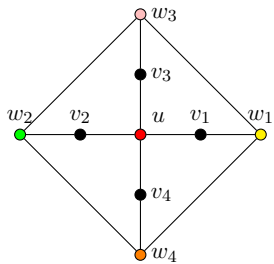
\begin{figure}[!htbp]
	\centering
	\begin{tikzpicture}[scale=.4, transform shape]
		
		\node [draw, shape=circle,fill=red] (v1) at  (0,0) {};
		\node [draw, shape=circle,fill=black] (v2) at  (2,0) {};
		\node [draw, shape=circle,fill=black] (v3) at  (-2,0) {};
		\node [draw, shape=circle,fill=black] (v4) at  (0,2) {};	
		\node [draw, shape=circle,fill=black] (v5) at  (0,-2) {};
		\node [draw, shape=circle,fill=yellow] (v6) at  (4,0) {};
		\node [draw, shape=circle,fill=green] (v7) at  (-4,0) {};		
		
		\node [draw, shape=circle,fill=pink] (v8) at  (0,4) {};
		\node [draw, shape=circle,fill=orange] (v9) at  (0,-4) {};
		
		\node [scale=2] at (0.6,0.6) {$u$};
		\node [scale=2] at (2,0.6) {$v_1$};
		\node [scale=2] at (-2,0.6) {$v_3$};
		\node [scale=2] at (0.7,2) {$v_2$};
		\node [scale=2] at (0.7,-2) {$v_4$};
		\node [scale=2] at (4,0.6) {$w_1$};
		\node [scale=2] at (-4,0.6) {$w_3$};
		\node [scale=2] at (0.8,4) {$w_2$};
		\node [scale=2] at (0.8,-4) {$w_4$};
		
		\draw (v1)--(v2);
		\draw (v1)--(v3);
		\draw(v1)--(v4);
		\draw (v1)--(v5);
		\draw (v2)--(v6);
		\draw(v3)--(v7);
		\draw (v4)--(v8);		
		\draw (v5)--(v9);
		
		\draw (v6)--(v8);
		\draw (v7)--(v9);
		\draw(v6)--(v9);
		\draw (v7)--(v8);
	\end{tikzpicture}

	\caption{The graph $G_4$.}
	\label{fig-G4}
\end{figure}


\section{Global coalition existence and general lower bounds}\label{sec-gl-exist}
In this section, we show that every graph $G$ of order at least two has a global coalition partition, and we present two general lower bounds on $GC(G)$. For this purpose, we recall the following concept, introduced by Zelinka \cite{zelinka} under the name \emph{complementarily domatic partition} and, independently, by Sampathkumar \cite{sampathkumar}.

\begin{definition}\label{def-gdom}
	A global domatic partition of a graph $G$ is a partition $\pi = \lbrace V_1, V_2, \dots, V_k \rbrace$ of $V (G)$ such that each set $V_i$ is a global dominating set of $G$.
	The maximum cardinality of a global domatic partition of $G$ is called the global domatic number of $G$, denoted by $d_g(G)$.
\end{definition}
Zelinka \cite{zelinka} proved, among other things, that $d_g(G)=d_g(\overline{G})$, that $d_g(G)\le\min\{\delta(G)+1,\,n-\Delta(G)\}$, that $d_g(G)=d(G)$ when $G$ is disconnected, and that $d_g(G)=1$ when $G$ has an isolated or a full vertex. We will use only the trivial fact that $d_g(G)\ge 1$, since $V(G)$ is a global dominating set.
\begin{theorem}
	\label{gb-exist}
	For any graph $G$ of order at least two, there is a $gc$-partition of $G$.
\end{theorem}
\begin{proof}
Let $\pi = \lbrace V_1,V_2, \dots ,V_k\rbrace$ be a global domatic partition of $G$, where  $d_g (G)=k$.
For each $1 \leq i \leq k-1$, if $V_i$ is not a minimal global dominating set, then we move $V_i \setminus V^\prime _i$ to $V_k$, where $V^\prime _i$ is a minimal global dominating set contained in $V_i$. After this step $V_i=V'_i$ for $i\le k-1$.
Note that $\lvert V_i\rvert \geq 2$ for each $1 \leq i \leq k$, because a singleton is never a global dominating set of a graph of order at least two. Now for each $1 \leq i \leq k-1$, we partition $V_i$ into  two non-empty sets $V_{i,1}$ and $V_{i,2}$, which yields a collection $\pi^\prime$ of $2(k-1)$ sets. Observe that the sets $V_{i,1}$ and $V_{i,2}$ are global coalition partners, since $V_i$ is a minimal global dominating set.
Now, consider $V_k$. If it  is a minimal global dominating set, then we partition it into two non-empty  subsets
and add them to $\pi^\prime$, which makes $\pi^\prime$ a $gc$-partition of $G$.
If $V_k$ is not a minimal global dominating set, then we partition $V^\prime _k$ into
two non-empty sets $V^\prime _{k,1}$ and $V^\prime _{k,2}$, where $V^\prime _k$ is a minimal global dominating set contained in $V_k$, and add them to $\pi^\prime$.
Let $A=V_k \setminus V^\prime _k$. The set $A$ is not a global dominating set, for otherwise we would have $d_g (G)\geq k+1$,    contradicting the maximality of $k$. If $A$ forms a global coalition with some member of $\pi^\prime $, then adding it to $\pi^\prime$ yields a $gc$-partition of $G$, as desired. Otherwise, $A$ is not a $gc$-partner of $V'_{k,2}$, and since neither $A$ nor $V'_{k,2}$ is a global dominating set, their union $V'_{k,2}\cup A$ is not a global dominating set either, whereas $(V'_{k,2}\cup A)\cup V'_{k,1}=V_k$ is. Hence we may replace $V'_{k,2}$ by $V'_{k,2}\cup A$ in $\pi'$, and the resulting partition is a $gc$-partition of $G$.
Thus, in all cases, we construct a $gc$-partition of $G$, completing the proof.
\end{proof}

Considering Theorem \ref{gb-exist}, its proof, and the fact that any global dominating set in a non-trivial graph has cardinality at least two, we readily derive the following result.
\begin{corollary}\label{cor-2dg}
	For any non-trivial graph $G$ with global domatic number $d_g (G)$, we have $GC(G) \geq 2 d_g (G)$.
\end{corollary}
The bound of Corollary~\ref{cor-2dg} is attained, for example, by complete graphs ($d_g(K_n)=1$, $GC(K_n)=2$) and by the cycle $C_6$ ($d_g(C_6)=3$ by \cite{zelinka} and $GC(C_6)=6$ by Theorem~\ref{gccn} below). Our second lower bound is of a different nature and explains all the lower-bound constructions used in Sections~\ref{sec-gl-co} and~\ref{sec-gl-uni}. For a graph $G$ of order $n\ge 2$, let $\mu(G)$ denote the minimum cardinality of an inclusion-maximal subset of $V(G)$ that is not a global dominating set of $G$. Since the empty set is not a global dominating set, a subset of this type exists.

\begin{proposition}\label{prop-mu}
For every graph $G$ of order $n\ge 2$ we have $GC(G)\ge n+1-\mu(G)$.
\end{proposition}
\begin{proof}
Let $S$ be an inclusion-maximal set which is not a global dominating set of $G$, with $|S|=\mu(G)$. Then $S\ne V(G)$, and by maximality $S\cup\{x\}$ is a global dominating set for every $x\in V(G)\setminus S$. Since no one-vertex class is a global dominating set, $\{S\}\cup\{\{x\}: x\in V(G)\setminus S\}$ is a $gc$-partition of $G$ with $n-\mu(G)+1$ classes, in which $S$ is a $gc$-partner of every other class.
\end{proof}
Proposition~\ref{prop-mu} provides another proof of Theorem~\ref{gb-exist}. A good way to obtain a maximal non-global-dominating set is to start from the open neighborhood $N(v)$ of a vertex $v$ of eccentricity two: the set $N(v)$ dominates $G$ but not $\overline{G}$ (the vertex $v$ has no neighbor in $N(v)$ in $\overline{G}$), and whenever every set $N(v)\cup\{w\}$ with $w\notin N(v)$ dominates $\overline{G}$, we obtain $GC(G)\ge n-\deg(v)+1$. This is exactly the situation encountered for trees of radius two (Theorem~\ref{gc-rad2}) and for the unicyclic graphs of Section~\ref{sec-gl-uni}. The bound of Proposition~\ref{prop-mu} is not always reached. For instance $\mu(P_4)=2$ while $GC(P_4)=4$.


\section{Global coalition vs. coalition}\label{sec-gl-co}
This section aims to investigate the relationship between the global coalition number and the coalition number of a graph.
In this regard, we establish several results for paths, cycles, trees, and graphs $G$ with $\girth(G) \geq6$.
First, we state and slightly extend the following result from \cite{mojdeh-ac}.

\begin{lemma} 	\label{lemm-rad3}
	
	Let $G$ be a graph of order at least two.
	\begin{enumerate}[(i)]
	\item If $G$ is disconnected or $\rad(G)\ge 3$, then every dominating set of $G$ is a dominating set of $\overline{G}$ \emph{(\cite{mojdeh-ac})}.
	\item If $G$ is connected and $\rad(G)\le 2$, then $G$ has a dominating set which is not a dominating set of $\overline{G}$.
	\end{enumerate}
\end{lemma}
\begin{proof}
(i) Let $D$ be a dominating set of $G$ and let $v\in V(G)\setminus D$. If $G$ is disconnected, then $D$ contains a vertex of a component not containing $v$, and this vertex is adjacent to $v$ in $\overline{G}$. If $\rad(G)\ge 3$, then there is a vertex $w$ at distance at least three from $v$. A vertex $d\in D$ dominating $w$ is at distance at least two from $v$, thus adjacent to $v$ in $\overline{G}$.
(ii) Let $v$ be a vertex with $\ecc(v)\le 2$. Since $G$ is connected of order at least two, $N(v)\ne\emptyset$, and every vertex other than $v$ is either a neighbour of $v$ or has a common neighbour with $v$. Thus $N(v)$ dominates $G$. However, $v$ has no neighbor in $N(v)$ in $\overline{G}$, so $N(v)$ does not dominate $\overline{G}$.
\end{proof}
Using Lemma \ref{lemm-rad3}, we can determine the global coalition number of graphs $G$ with $\rad(G)\geq 3$, and of disconnected graphs, in terms of their coalition number as follows.

\begin{theorem}\label{rad-c-gc}
	If $G$ is a graph with $\rad(G) \geq 3$, or a disconnected graph of order at least two, then $GC(G)=C(G)$.
\end{theorem}
\begin{proof}
Since $\rad(G) \geq 3$ or $G$ is disconnected, it follows that $G$ contains no full vertices.
Let $\pi$ be a $c$-partition of $G$, and  let $A \in \pi$. Since $G$ has no full vertices,  the set $A$ cannot be a singleton dominating set of $G$, and so it must have a $c$-partner in $\pi$.
Let $B$ be a $c$-partner of $A$. By definition, neither $A$ nor $B$ is a dominating set of $G$,
implying that neither $A$ nor $B$ is a global dominating set of $G$.
However, since $A \cup B$ forms a dominating set of $G$,   it follows from Lemma \ref{lemm-rad3}(i)
that $A \cup B$ is a global dominating set of $G$.
Therefore, $\pi$ is a $gc$-partition of $G$, implying that $GC(G) \geq C(G)$.
Now, assume that $\pi^\prime$ is  a $gc$-partition of $G$. Let $C \in \pi^\prime$, and let $D$ be a $gc$-partner of $C$.
Since $C$ is not a global dominating set of $G$,  it follows from Lemma \ref{lemm-rad3}(i) that $C$ is not a dominating set of $G$. A similar argument shows that $D$ is not a dominating set of $G$ either.
Further, since $C \cup D$ is a global dominating set of $G$, it follows that $C \cup D$ is also a dominating set of $G$.
Thus, $C$ forms a coalition with $D$,  implying that the partition $\pi^\prime$ is a $c$-partition of $G$.
Consequently, $C(G) \geq GC(G)$, and so $C(G)=GC(G)$. This completes the proof.
\end{proof}
Combining Theorem~\ref{rad-c-gc} with Observation~\ref{obs-compl}, we obtain a reduction for all graphs with a disconnected complement, that is, for all joins.
\begin{corollary}\label{cor-join}
Let $G$ be a graph of order at least two such that $\overline{G}$ is disconnected or $\rad(\overline{G})\ge 3$. Then $GC(G)=C(\overline{G})$. In particular, $GC(G_1+G_2)=C(\overline{G_1}\cup\overline{G_2})$ for every join $G_1+G_2$ of two graphs.
\end{corollary}
\begin{proof}
By Observation~\ref{obs-compl} and Theorem~\ref{rad-c-gc} applied to $\overline{G}$, we get $GC(G)=GC(\overline{G})=C(\overline{G})$. The second statement follows since $\overline{G_1+G_2}=\overline{G_1}\cup\overline{G_2}$ is disconnected.
\end{proof}
Since a graph and its complement cannot both be disconnected, Theorem~\ref{rad-c-gc} and Corollary~\ref{cor-join} show that $GC(G)$ can differ from both $C(G)$ and $C(\overline{G})$ only if $G$ and $\overline{G}$ are connected and both have radius at most two. Lemma~\ref{lemm-rad3}(ii) shows that this is exactly the class of graphs for which global domination is a significantly stronger requirement than domination in both $G$ and $\overline{G}$ separately. Note also that Corollary~\ref{cor-join} gives a second proof of Observations~\ref{obs-kn} and~\ref{obs-kbp} and of Propositions~\ref{prop-wheel} and~\ref{prop-fan}, provided the corresponding coalition numbers of the disconnected complements are known. We preferred the direct arguments, which are short.

\subsection{Trees}\label{sub-gl-tree}
As a direct implication of Theorem~\ref{rad-c-gc}, we have the following.

\begin{corollary}\label{cor-tree-rad3}
	If $T$ is a tree with $\rad(T) \geq 3$, then $GC(T)=C(T)$.
\end{corollary}

The coalition number of trees was studied by Bakhshesh, Henning and Pradhan \cite{BHP-trees}, who showed in particular that $C(T)\le 6$ for every tree $T$ and characterised the trees attaining this bound. By Corollary~\ref{cor-tree-rad3} the same applies to $GC(T)$ when $\rad(T)\ge 3$. Trees of radius one are stars, for which $GC(K_{1,n-1})=n$ by Observation~\ref{obs-kbp}. It remains to consider trees of radius two.

Before we determine the global coalition number of trees $T$ with $\rad(T)=2$, let us first recall the following well-known theorem.

\begin{figure}[!htbp]
	\centering
	\begin{tikzpicture}[scale=.4, transform shape]
		\node [draw, shape=circle,fill=blue] (v1) at  (0,0) {};
		\node [draw, shape=circle,fill=black] (v2) at  (0,3) {};
		\node [draw, shape=circle,fill=red] (v3) at  (-2.5,1.5) {};
		\node [draw, shape=circle,fill=green] (v4) at  (2.5,1.5) {};	
		\node [draw, shape=circle,fill=yellow] (v5) at  (-1.5,-2.5) {};
		\node [draw, shape=circle,fill=pink] (v6) at  (1.5,-2.5) {};
		
		\node [draw, shape=circle,fill=blue] (v7) at  (0,6) {};
		\node [draw, shape=circle,fill=blue] (v8) at  (-2.5,4.5) {};
		\node [draw, shape=circle,fill=blue] (v9) at  (2.5,4.5) {};

		\node [scale=2] at (-0.4,3.7) {$1$};
		\node [scale=2] at (-0.4,0.7) {$2$};
		\node [scale=2] at (0,6.7) {$2$};
		\node [scale=2] at (-2.9,5.2) {$2$};
		\node [scale=2] at (2.9,5.2) {$2$};
		\node [scale=2] at (-2.9,2.2) {$3$};
		\node [scale=2] at (2.9,2.2) {$4$};
		\node [scale=2] at (-1.5,-3.2) {$6$};
		\node [scale=2] at (1.5,-3.2) {$5$};
		
		\draw (v1)--(v2);
		\draw (v1)--(v3);
		\draw(v1)--(v4);
		\draw (v1)--(v5);
		\draw (v1)--(v6);
		\draw(v2)--(v7);
		\draw(v2)--(v8);		
		
		\draw (v2)--(v9);
		
		\node [draw, shape=circle,fill=black] (v1) at  (-10,1.5) {};
		\node [draw, shape=circle,fill=blue] (v2) at  (-10,4.5) {};
		\node [draw, shape=circle,fill=blue] (v3) at  (-12.5,3) {};
		\node [draw, shape=circle,fill=blue] (v4) at  (-7.5,3) {};	
		\node [draw, shape=circle,fill=blue] (v5) at  (-11.5,-1) {};
		\node [draw, shape=circle,fill=blue] (v6) at  (-8.5,-1) {};
		\node [draw, shape=circle,fill=blue] (u6) at  (-12.5,1) {};		
		\node [draw, shape=circle,fill=green] (v8) at  (-12.5,6) {};
		\node [draw, shape=circle,fill=yellow] (v9) at  (-7.5,6) {};	
		\node [draw, shape=circle,fill=pink] (v10) at  (-15,4) {};	
		\node [draw, shape=circle,fill=orange] (v14) at  (-14,-1.5) {};				
		\node [draw, shape=circle,fill=white] (v15) at  (-13.5,-2.5) {};
		\node [draw, shape=circle,fill=teal] (v16) at  (-13,-3.5) {};		
		\node [draw, shape=circle,fill=red] (v17) at  (-7,-3) {};
		
		\node [scale=2] at (-10,5.2) {$2$};
		\node [scale=2] at (-10.4,2.2) {$1$};
		\node [scale=2] at (-12.5,3.7) {$2$};
		\node [scale=2] at (-7.5,3.7) {$2$};
		\node [scale=2] at (-12.5,1.7) {$2$};
		\node [scale=2] at (-11.5,-.3) {$2$};
		\node [scale=2] at (-8.5,-.3) {$2$};
		\node [scale=2] at (-12.5,6.7) {$3$};
		\node [scale=2] at (-7.5,6.7) {$4$};
		\node [scale=2] at (-7,-3.7) {$5$};
		\node [scale=2] at (-13,-4.2) {$6$};
		\node [scale=2] at (-13.9,-3.2) {$7$};
		\node [scale=2] at (-14.6,-1.9) {$8$};
		\node [scale=2] at (-15.6,4) {$9$};

		\draw (v1)--(v2);
		\draw (v1)--(v3);
		\draw(v1)--(v4);
		\draw (v1)--(v5);
		\draw (v1)--(v6);
		\draw(v2)--(v8);
		\draw (v1)--(u6);		
		
		\draw (v2)--(v9);
		\draw (v3)--(v10);
		\draw(v5)--(v14);
		\draw(v5)--(v15);	
		\draw(v5)--(v16);	
		\draw(v6)--(v17);	
	\end{tikzpicture}

	\caption{Examples of $GC(T)$-partitions for trees with radius two.}\label{fig-T-d-star}
\end{figure}
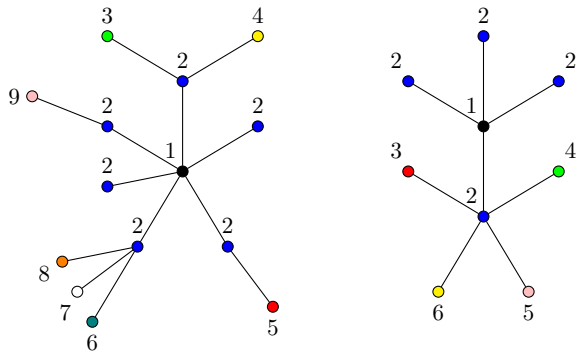

\begin{theorem}
	\label{rad-diam}
	(\cite{chartrand}) For any connected graph $G$,
	$ 	\rad(G) \leq \diam(G) \leq 2 \rad(G).	$
\end{theorem}
Note that trees $T$ with $\rad(T)=2$ and $\diam(T)=4$ possess a unique vertex $v$ with  $\ecc(v)=2$.
For such trees, we define $\ell$ to be the number of vertices in $T$ that are at distance $2$ from $v$. Thus $\ell=n-\deg(v)-1$. However, the trees $T$ with $\rad(T)=2$ and $\diam(T)=3$ are
double stars, which possess two
vertices with eccentricity $2$ (see Figure \ref{fig-T-d-star}).

With this setup, we are now prepared to determine the global coalition number of trees $T$ with $\rad(T)=2$. Note that $P_4=S_{1,1}$ is excluded below: $GC(P_4)=4$ by Theorem~\ref{gc-paths}, while the formula would give $3$.

\begin{theorem}
	\label{gc-rad2}
	For any tree $T \neq P_4$  with $\rad(T)=2$,
	\[
	GC(T) = \left\{
	\begin{array}{ll}
		\ell + 2 \ =\ n-\deg(v)+1 & \mbox{if\ $\diam(T)=4$, where $v$ is the central vertex of $T$}, \\
		\max\{p,q\}+2 & \mbox{if $T=S_{p,q}$}.
	\end{array}
	\right.
	\]
\end{theorem}
\begin{proof}
By Theorem \ref{rad-diam}, either $\diam(T)=3$ or $\diam(T)=4$. Accordingly, we consider two cases. In both cases, the lower bound is the special case of Proposition~\ref{prop-mu} in which the maximal non-global-dominating set is the open neighbourhood of a vertex of eccentricity two. We give the partition explicitly.

\textbf{Case 1.} $\diam(T)=3$.

In this case, $T$ is a double star. Let $T=S_{p,q}$ with $p \geq q$. Since $T\ne P_4$, we have $p\ge 2$. Let $a,b \in V(T)$ be the support vertices of $T$
such that $N(a)=\lbrace b,v_1,v_2,\dots ,v_p\rbrace$ and $N(b)=\lbrace a,u_1,u_2,\dots ,u_q\rbrace$, and let $L_a=\{v_1,\dots,v_p\}$, $L_b=\{u_1,\dots,u_q\}$ and $L=L_a\cup L_b$.
Then  the partition $\lbrace \lbrace b\rbrace ,\lbrace a,u_1,u_2,\dots ,u_q\rbrace ,\lbrace v_1\rbrace ,\dots ,\lbrace v_p\rbrace \rbrace$
is a $gc$-partition of $T$ of size $p+2$, where the set  $\lbrace a,u_1,u_2,\dots,u_q\rbrace$ $=N(b)$ forms
a global coalition with all other members of the partition. Thus, $GC(T) \geq p+2$. Now we show that $GC(T)\le p+2$.

We first describe the global dominating sets of $T$. Every set containing both $a$ and $b$ is a global dominating set. A set $S$ with $a\in S$, $b\notin S$ dominates $T$ if and only if $L_b\subseteq S$, and it dominates $\overline{T}$ if and only if $b$ has a neighbor in $S$ in $\overline{T}$, i.e., if and only if $S\cap L_a\ne\emptyset$ (the other vertices are then dominated by $a$ or by $L_b$). Symmetrically, a set $S$ with $b\in S$, $a\notin S$ is a global dominating set if and only if $L_a\subseteq S$ and $S\cap L_b\ne\emptyset$. Finally, a set avoiding $a$ and $b$ dominates $T$ only if it contains $L$, and $L$ is a global dominating set.

Let $\pi$ be a $gc$-partition of $T$, and let $V_a$ and $V_b$ be the classes containing $a$ and $b$. They are distinct, as a class containing $\{a,b\}$ would be a global dominating set. Every other class is contained in $L$. If two classes $C,D\notin\{V_a,V_b\}$ are $gc$-partners, then $C\cup D=L$ and $\pi=\{V_a,V_b,C,D\}$, so $|\pi|=4\le p+2$. Thus, we may assume that every class other than $V_a$ and $V_b$ is a $gc$-partner of $V_a$ or of $V_b$. We consider three subcases.

\textbf{Subcase 1.1.} $L_b\subseteq V_a$. Then $V_a\cap L_a=\emptyset$, for otherwise $V_a$ would be a global dominating set. Thus $V_a=\{a\}\cup L_b$ and $V_b\subseteq\{b\}\cup L_a$. A class $C\subseteq L_a$ cannot be a $gc$-partner of $V_b$, because $V_b\cup C$ contains $b$ but not $a$, and it does not meet $L_b$ (recall that $L_b\subseteq V_a$). Hence all classes other than $V_a$ and $V_b$ are $gc$-partners of $V_a$; they are pairwise disjoint subsets of $L_a\setminus V_b$, and so $|\pi|\le p+2$.

\textbf{Subcase 1.2.} $L_a\subseteq V_b$. By symmetry, $|\pi|\le q+2\le p+2$.

\textbf{Subcase 1.3.} $L_b\not\subseteq V_a$ and $L_a\not\subseteq V_b$. Every $gc$-partner $C\ne V_b$ of $V_a$ satisfies $L_b\subseteq V_a\cup C$, so $C\supseteq L_b\setminus V_a\ne\emptyset$. Since the classes are disjoint, $V_a$ has at most one $gc$-partner other than $V_b$. The same holds for $V_b$, and therefore $|\pi|\le 4\le p+2$.

\textbf{Case 2.} $\diam(T)=4$.

Let $a\in V(T)$ be the vertex with $\ecc(a)=2$. Let $U_1=N(a)$, and $U_2=V(T) \setminus N[a]$,
where $U_2=\lbrace v_1,v_2,\dots ,v_\ell \rbrace$. Observe that each vertex in $U_2$ has distance $2$ from $a$, and that each vertex $y\in U_2$ is a leaf whose unique neighbour, denoted $p(y)$, lies in $U_1$. Let $A=\{x\in U_1: N(x)\cap U_2\ne\emptyset\}$ and $B=U_1\setminus A$. The vertices of $B$ are leaves adjacent to $a$, and for $x\in A$ we write $\mathrm{ch}(x)=N(x)\cap U_2$. Since $\diam(T)=4$, we have $|A|\ge 2$.
Now the partition $\lbrace \lbrace a\rbrace ,U_1, \lbrace v_1 \rbrace ,\dots ,\lbrace v_\ell \rbrace \rbrace$ is a $gc$-partition of $T$ of
cardinality $\ell +2$, where the set $U_1$ forms a global coalition with all other   partition classes. Indeed, $U_1\cup\{a\}=N[a]$ dominates $T$ and contains a non-neighbor of every vertex of $U_1$ and of $U_2$ in $\overline{T}$, and $U_1\cup\{y\}$ ($y\in U_2$) dominates $T$, while in $\overline{T}$ the vertex $a$ is dominated by $y$ and every other vertex of $U_2$ is adjacent to $y$.
Thus, we have $GC(T) \geq \ell + 2$.
We now show that  $GC(T)\le \ell+2$. Let $\pi$ be a $gc$-partition of $T$, and let $V_a$ be the class containing $a$. We use two claims.

\textbf{Claim 1.} \emph{$V_a$ has at most two $gc$-partners.} We first show that $V_a$ does not dominate $T$. Suppose it does, i.e., for every $y\in U_2$ we have $y\in V_a$ or $p(y)\in V_a$. As $V_a$ is not a global dominating set, some vertex $u$ is not dominated by $V_a$ in $\overline{T}$. Since $a\in V_a$ is adjacent in $\overline{T}$ to every vertex of $U_2$, we have $u\in U_1\setminus V_a$, and the neighbors of $u$ in $\overline{T}$, namely all vertices of $(U_1\setminus\{u\})\cup(U_2\setminus \mathrm{ch}(u))$, avoid $V_a$. Thus $V_a\subseteq\{a\}\cup\mathrm{ch}(u)$. Take $x\in A\setminus\{u\}$ and $y\in\mathrm{ch}(x)$. Then $y\notin\mathrm{ch}(u)$, so $y\notin V_a$, and $p(y)=x\notin V_a$, a contradiction. Hence there is $y_0\in U_2$ with $y_0\notin V_a$ and $p(y_0)\notin V_a$. Every $gc$-partner of $V_a$ must dominate $y_0$, hence contains $y_0$ or $p(y_0)$, and so $V_a$ has at most two $gc$-partners.

\textbf{Claim 2.} \emph{If $C,D\in\pi$ are $gc$-partners with $a\notin C\cup D$, then $C\cup D\supseteq B$ and, for every $x\in A$, $C\cup D$ contains $x$ or all of $\mathrm{ch}(x)$. In particular $|C\cup D|\ge|U_1|$, with equality only if $C\cup D=B\cup R$, where $R$ contains exactly one vertex from $\{x\}\cup\mathrm{ch}(x)$ for each $x\in A$ (and this vertex is $x$ unless $|\mathrm{ch}(x)|=1$).} Indeed, $C\cup D$ dominates $T$ and does not contain $a$. A vertex of $B$ is dominated only by itself and by $a$, and a vertex $y\in U_2$ is dominated only by itself and by $p(y)$.

If no two classes of $\pi\setminus\{V_a\}$ are $gc$-partners, then every class other than $V_a$ is a $gc$-partner of $V_a$, and $|\pi|\le 3\le \ell+2$ by Claim~1. Otherwise, fix $gc$-partners $C,D$ with $a\notin C\cup D$. The classes of $\pi\setminus\{C,D\}$ partition $V(T)\setminus(C\cup D)$, which by Claim~2 has at most $n-|U_1|=\ell+1$ vertices. Hence $|\pi|\le \ell+3$, and equality forces $|C\cup D|=|U_1|$ and all classes other than $C$ and $D$ to be singletons. Suppose that $|\pi|=\ell+3$. Then $V_a=\{a\}$. The $gc$-partner of $\{a\}$ is not a singleton $\{z\}$, because $\{a,z\}$ dominates $T$ only if every vertex of $U_2$ is $z$ or a child of $z$, which contradicts $|A|\ge 2$. So the partner of $\{a\}$ is $C$ or $D$, say $C$. Then $\{a\}\cup C$ dominates $T$, so $C$ contains $x$ or $\mathrm{ch}(x)$ for every $x\in A$. By Claim~2, $C\cup D=B\cup R$, and we conclude that $R\subseteq C$ and $D\subseteq B$. In particular, $B\ne\emptyset$. Now choose $x_1\in A$ and a vertex $w\in(\{x_1\}\cup\mathrm{ch}(x_1))\setminus(C\cup D)$, which exists since $R$ contains only one vertex of $\{x_1\}\cup\mathrm{ch}(x_1)$. The singleton $\{w\}$ is a class of $\pi$. Let $F$ be its $gc$-partner. As shown above, $F\ne\{a\}$, so $a\notin\{w\}\cup F$ and, as in Claim~2, $B\subseteq\{w\}\cup F$. Since $w\notin B$, we get $B\subseteq F$. If $F=C$, this contradicts $\emptyset\ne D\subseteq B$ and $C\cap D=\emptyset$. If $F$ is a singleton, then $B=F=D$. So in all cases $F=D=B$, and then $\{w\}\cup F=\{w\}\cup B$ must dominate $T$. But for $x_2\in A\setminus\{x_1\}$ neither $x_2$ nor its children belong to $\{w\}\cup B$, a contradiction. Therefore $|\pi|\le\ell+2$, which completes the proof.
\end{proof}

Figure \ref{fig-T-d-star} illustrates a $GC(T)$-partition for two trees of radius two. One with $\diam(T)=3$ and another with $\diam(T)=4$. To illustrate more clearly, a distinct number and colour are assigned to each partition class. As seen in Figure \ref{fig-T-d-star}, the set of vertices with blue colour forms a global
coalition with all other partition classes.


\subsection{Graphs $G$ with $\girth(G) \geq 6$}\label{sub-gl-girth-6}
The following lemma allows us to express the global coalition number of graphs with $\girth(G) \geq 6$ in terms of their coalition number. Note that the assumption that $G$ contains a cycle is necessary, as the center of the star $K_{1,3}$ is a dominating set of $K_{1,3}$ which does not dominate $\overline{K_{1,3}}$, and indeed $GC(K_{1,m})=m+1\ne 3=C(K_{1,m})$ for $m\ge 3$. 

\begin{lemma}
	\label{gl-girth6}
	Let $G$ be a graph containing a cycle, with $\girth(G) \geq 6$. Then every dominating set of $G$ is a dominating set of $\overline{G}$.
\end{lemma}
\begin{proof}

 If $G$ is disconnected, the claim follows from Lemma~\ref{lemm-rad3}(i), so let $G$ be connected. We show that $\rad (G) \geq 3$. Suppose, to the contrary, that $\rad (G) < 3$. If $\rad (G)=1$, then $G$ has a full vertex $v$. A cycle of $G$ contains an edge $xy$ with $v\notin\{x,y\}$, and then $vxy$ is a triangle, contradicting $\girth(G)\ge 6$. Hence, we may assume that $\rad (G)=2$. Let $v \in V(G)$ be a vertex with $\ecc (v)=2$. Let $C$ be a cycle in $G$ of length $\girth(G)$. First, we show that $v \notin V(C)$. Suppose to the contrary that $v\in V(C)$. Let $u \in V(C)$ be such that $d_C (u,v)=3$. It follows that $d_G (u,v)=3$, for otherwise a $u$--$v$ path of length at most two together with the $u$--$v$ path of length three along $C$ would contain a cycle of length at most five, and $G$ would contain a cycle of length less than $6$, which is a contradiction. This contradicts the assumption that $\ecc(v)=2$. Hence, $v \notin V(C)$. Note also that $v$ has at most one neighbour in $V(C)$, for otherwise a cycle of length less than $|V(C)|$ would be created, which is a contradiction. Therefore, we can choose an edge $ab \in E(C)$ such that $v$ has no neighbour in $\lbrace a,b\rbrace$. Let $P_1=(a,x,v)$ be a path in $G$ connecting $a$ and $v$, and let $P_2=(b,y,v)$ be a path in $G$ connecting $b$ and $v$. If $x=y$, then we have a triangle, which is a contradiction. Otherwise, $v, x, a, b, y$ are five distinct vertices, and we have a $5$-cycle, which is again a contradiction. Hence, $\rad (G) \geq 3$. Now, by Lemma \ref{lemm-rad3}(i), the desired result follows.
\end{proof}
\begin{corollary}\label{trm-girth6}
	If $G$ is a graph containing a cycle, with $\girth(G) \geq 6$, then $GC(G)=C(G)$.
\end{corollary}
\begin{proof}
Considering the proof of Lemma \ref{gl-girth6}, $G$ is disconnected or $\rad(G) \geq 3$. The desired result now follows from Theorem \ref{rad-c-gc}.
\end{proof}
\subsection{Paths and cycles}\label{sub-gl-pa-cy}
In \cite{haynes3}, Haynes et al. determined the coalition number of paths as follows.
\begin{theorem}
	\label{c-paths}
	(\cite{haynes3})	For the path $P_n$,
	$$
	C(P_n)= \begin{cases}
		n & \text{if }  n \leq 4, \\
		4 & \text{if }  n=5, \\
		5 & \text{if }  6\leq n\leq 9, \\
		6 & \text{if }  n\geq 10.
	\end{cases}
	$$
\end{theorem}
Using Theorems \ref{rad-c-gc}, \ref{rad-diam} and \ref{c-paths}, we can determine the global coalition number of paths as follows.
\begin{theorem}\label{gc-paths}
	For the path $P_n$,
	$$
	GC(P_n)= \begin{cases}
		n & \text{if }  2\leq n \leq 4, \\
		4 & \text{if }  n=5, \\
		5 & \text{if }  6\leq n\leq 9, \\
		6 & \text{if }  n\geq 10.
	\end{cases}
	$$
\end{theorem}
\begin{proof}
The result is easily verified for values of $n$ such that  $2\leq n\leq 4$: $P_2=K_2$ and $P_3=K_{1,2}$ are covered by Observations~\ref{obs-kn} and~\ref{obs-kbp}, and in $P_4=(v_1,v_2,v_3,v_4)$ the pairs $\{v_1,v_4\}$ and $\{v_2,v_3\}$ are global dominating sets, so the partition into singletons is a $gc$-partition.  Now, consider the case $n=5$. Let
the vertex set of $P_5$ be $\lbrace v_1,v_2,v_3,v_4,v_5\rbrace$.  It is straightforward to observe that the partition of $V(P_5)$ into singleton sets does not constitute a valid global coalition partition
for $P_5$ (for instance, no set $\{v_3,v_i\}$ is a global dominating set). However,  the partition $\lbrace \lbrace v_1\rbrace ,\lbrace v_2,v_3\rbrace ,\lbrace v_4\rbrace ,\lbrace v_5\rbrace \rbrace$ is a $gc$-partition of $P_5$. Hence, $GC(P_5)=4$. It remains to prove the result when $n\geq 6$. In this case, we observe that $\diam(P_n)\geq 5$. Therefore, by Theorem \ref{rad-diam}, it follows that $\rad(P_n) \geq 3$. Consequently, by Theorem \ref{rad-c-gc}, we obtain  $GC(P_n)=C(P_n)$.
This completes the proof.
\end{proof}
In \cite{haynes3}, Haynes et al. determined the coalition number of cycles as follows.
\begin{theorem}
	\label{c-cycles}
	(\cite{haynes3})	For the cycle $C_n$,
	$$
	C(C_n)= \begin{cases}
		n & \text{if }  n \leq 6, \\
		5 & \text{if }  n=7,\\
		6 & \text{if }  n\geq 8.
	\end{cases}
	$$
\end{theorem}
The following theorem determines the global coalition number of cycles.
\begin{theorem}
	\label{gccn}
	For the cycle $C_n$,
	$$
	GC(C_n)= \begin{cases}
		2 & \text{if }  n=3, \\
		4 & \text{if }  n =4,5, \\
		5 & \text{if }  n=7, \\
		6 & \text{otherwise}.
	\end{cases}
	$$
\end{theorem}
\begin{proof}
Consider the cycle $C_n$ with $V(C_n)=\lbrace v_1,v_2,\dots ,v_n\rbrace$.
For $n=3$, we have $GC(C_3)=GC(K_3)=2$ by Observation~\ref{obs-kn}.
Now, let $n=4$.  The singleton partition of $V(C_4)$ is a valid $gc$-partition of $C_4$, since $\{v_1,v_2\}$ and $\{v_3,v_4\}$ are global dominating sets. Hence, $GC(C_4)=4$. Next, assume $n=5$. It is easy to check that the singleton partition of $V(C_5)$  does not form a valid $gc$-partition of $C_5$, because no two vertices of $C_5$ form a global dominating set. However, the partition
$\lbrace \lbrace v_1,v_3\rbrace ,\lbrace v_2\rbrace,\lbrace v_4\rbrace, \lbrace v_5\rbrace  \rbrace$
is a  valid $gc$-partition of $C_5$,  implying that $GC(C_5)=4$.
For $n\geq 6$, the desired result follows from Corollary \ref{trm-girth6} and Theorem \ref{c-cycles}.
\end{proof}


\section{Graphs with GC(G)=n and perfect coalitions}\label{sec-gl-perf}
We begin with a characterisation of the graphs whose global coalition number equals their order. Following Olariu \cite{olariu}, two vertices $x$ and $y$ of a graph $G$ are called \emph{antitwins} if every vertex $z\in V(G)\setminus\{x,y\}$ is adjacent to exactly one of $x$ and $y$.
\begin{theorem}\label{thm-antitwin}
Let $G$ be a graph of order $n\ge 2$. Then $GC(G)=n$ if and only if every vertex of $G$ has an antitwin.
\end{theorem}
\begin{proof}
We have $GC(G)=n$ if and only if the partition of $V(G)$ into singletons is a $gc$-partition. Since no singleton is a global dominating set, this happens if and only if for every vertex $x$ there is a vertex $y\ne x$ such that $\{x,y\}$ is a global dominating set. Now $\{x,y\}$ dominates $G$ if and only if every $z\notin\{x,y\}$ is adjacent to $x$ or to $y$, and it dominates $\overline{G}$ if and only if every such $z$ is non-adjacent to $x$ or to $y$. Both conditions hold simultaneously if and only if $x$ and $y$ are antitwins.
\end{proof}
For instance, in a star every leaf and the centre are antitwins, in $P_4=(v_1,v_2,v_3,v_4)$ the pairs $\{v_1,v_4\}$ and $\{v_2,v_3\}$ are antitwins, and $C_5$ has no antitwins. Note also that, since $x,y$ are antitwins in $G$ if and only if they are antitwins in $\overline{G}$, Theorem~\ref{thm-antitwin} is consistent with Observation~\ref{obs-compl}.

In the remainder of this section, we examine the relation between the concept of global coalition and that of perfect coalition, which was introduced and investigated in \cite{mojdeh-s}. In \cite{mojdeh-s}, the authors defined the concepts of perfect coalition and perfect coalition partition as follows.
\begin{definition}\label{def-prc} (\cite{mojdeh-s})
	A perfect coalition in a graph $G$ consists of two disjoint sets of vertices $V_i$ and $V_j$ such that
	\begin{enumerate}
		\item[(i)]
		Neither $V_i$ nor $V_j$ is a dominating set of $G$.
		\item[(ii)]
		Each vertex in $V(G) \setminus V_i$ has at most one neighbour in $V_i$, and each vertex in $V(G) \setminus V_j$ has at most one neighbour in $V_j$.
		\item[(iii)]
		$V_i \cup V_j$ is a perfect dominating set of $G$.
	\end{enumerate}
\end{definition}
\begin{definition}\label{def-prcp} (\cite{mojdeh-s})
	A perfect coalition partition (abbreviated $prc$-partition) in a graph $G$ is a vertex partition $\pi= \lbrace V_1, V_2, \dots, V_k \rbrace$ such that for each set $V_i$ of $\pi$ either $V_i$ is a singleton dominating set, or there exists a set  $V_j \in \pi$  that forms a perfect coalition with $V_i$. The maximum cardinality of a $prc$-partition in $G$ is called the perfect coalition number of $G$, denoted by $PRC(G)$. A $prc$-partition of $G$  with cardinality  $PRC(G)$ is called a $PRC(G)$-partition.
\end{definition}
\begin{theorem}\label{thm-gc-prc}
	If $G$ is a graph with no full vertices, then $GC(G) \geq PRC(G)$.
\end{theorem}
\begin{proof}
Let $\pi$ be a $prc$-partition of $G$. Since $G$ has no full vertices, each member of $\pi$ must have a $prc$-partner in $\pi$. Consider an arbitrary pair $(A,B)$ of $prc$-partners in $\pi$.
By definition, neither $A$ nor $B$ is a global dominating set of $G$. Furthermore, since each vertex in $V(G) \setminus (A \cup B)$ has exactly one neighbour in $A \cup B$ and $|A\cup B|\ge 2$, it follows that each such vertex has a non-neighbour in $A\cup B$, i.e., $A \cup B$ forms a dominating set of $\overline{G}$.
Thus, $A \cup B$ is a global dominating set of $G$, which implies that the sets $A$ and $B$ are $gc$-partners in $G$.
Since the pair  $(A, B)$ was chosen arbitrarily,  we conclude that $\pi$ is a $gc$-partition of $G$. Therefore, $GC(G) \geq PRC(G)$.
\end{proof}
\begin{theorem}
	\label{glb-n}
	Let $G$ be a graph of order $n$ with no full vertices. Then $GC(G)=n$ if and only if $PRC(G)=n$.
\end{theorem}
\begin{proof}

By Theorem~\ref{thm-antitwin}, $GC(G)=n$ if and only if every vertex of $G$ has an antitwin. On the other hand, $PRC(G)=n$ if and only if the partition of $V(G)$ into singletons is a $prc$-partition. Since $G$ has no full vertices, no singleton is a dominating set, so this holds if and only if every singleton $\{x\}$ has a $prc$-partner $\{y\}$. Condition (ii) of Definition~\ref{def-prc} is trivially satisfied by singletons, and $\{x,y\}$ is a perfect dominating set if and only if every vertex outside $\{x,y\}$ has exactly one neighbor in $\{x,y\}$, i.e., if and only if $x$ and $y$ are antitwins. Thus, both conditions are equivalent to every vertex having an antitwin.
\end{proof}

Next we characterize trees $T$ with $GC(T)=\lvert V(T)\rvert$. For this purpose, we shall make use of the following result, which has been obtained in \cite{mojdeh-s}.
\begin{corollary} \label{tree-n-prc}
	(\cite{mojdeh-s}) Let $T$ be a tree of order $n$. Then $PRC(T)=n$ if and only if $T \in \lbrace P_1,P_2,P_4\rbrace$.
\end{corollary}
\begin{observation}\label{obs-tree-n}
	Let $T$ be a tree of order $n>1$. Then $GC(T)=n$ if and only if either $T$ is a star or $T=P_4$.
\end{observation}
\begin{proof}
If $T$ is a star, then by Observation \ref{obs-kbp}, $GC(T)=n$. Also, it is clear that $GC(P_4)=4$. Conversely, assume that $GC(T)=n$. If $T$ has a full vertex, then $T$ is a star. Otherwise, it follows from Theorem \ref{glb-n} that $PRC(T)=n$, and so, by Corollary \ref{tree-n-prc}, we have $T=P_4$.
\end{proof}

Note that the difference between the global coalition number and the perfect coalition number of a graph can be arbitrarily large. For example, consider the double star $S_{p,q}$ with $2 \leq p \leq q$. By Theorem \ref{gc-rad2}, we have $GC(S_{p,q})=q+2$, while it is straightforward to verify that $PRC(S_{p,q})=2$.

In \cite{mojdeh-s}, the authors also provided a characterisation of triangle-free graphs $G$ of order $n \geq 4$ for which $PRC(G)=n$, which we now recall.
\begin{definition}(\cite{mojdeh-s})
	\label{family2}
	Let $\mathcal{T}_1$ denote the family of bipartite graphs $H- M$, where $H=K_{r,r}$ is a complete bipartite graph with $r \geq 2$ and $M$ is a perfect matching in $H$.
	Let $\mathcal{T}_2$ denote the family of bipartite graphs $H- M$, where $H=K_{r,s}$ is a
	complete bipartite graph with $r,s \geq 2$, and $M$ (possibly empty) is a  matching in $H$ such that $\lvert M\rvert < \min \lbrace r,s\rbrace$.
\end{definition}
\begin{theorem}
	\label{t-free-prc}
	(\cite{mojdeh-s}) Let $G$ be a triangle-free graph of order $n\geq 4$.
	Then $PRC(G)=n$ if  and only if $G \in \mathcal{T}_1 \cup \mathcal{T}_2$.
\end{theorem}
Using Theorems~\ref{glb-n} and~\ref{t-free-prc}, we obtain the following result.
\begin{corollary}\label{cor-tfree}
	Let $G$ be a triangle-free graph of order $n\geq 4$. Then $GC(G)=n$ if  and only if either $G$ is a star or $G \in \mathcal{T}_1 \cup \mathcal{T}_2$.
\end{corollary}
\begin{proof}
If $G$ has a full vertex, then, being triangle-free, $G$ is a star, and $GC(G)=n$ by Observation~\ref{obs-kbp}. Otherwise, the claim follows from Theorems~\ref{glb-n} and~\ref{t-free-prc}. Conversely, stars satisfy $GC(G)=n$, and graphs in $\mathcal{T}_1\cup\mathcal{T}_2$ have no full vertices.
\end{proof}
We remark that Corollary~\ref{cor-tfree} can also be derived directly from Theorem~\ref{thm-antitwin}, by analysing pairs of antitwins in triangle-free graphs. We skip the details.

\section{Global coalition  in unicyclic graphs}\label{sec-gl-uni}
In this section, we investigate the concept of global coalition in unicyclic graphs. A unicyclic graph is
a connected graph that contains exactly one cycle. According to Theorem \ref{rad-c-gc},
the global coalition number of unicyclic graphs $G$ with $\rad(G) \geq 3$  is equal to their coalition number.  Motivated by this result, we study unicyclic graphs $G$ with $\rad (G) \leq 2$. As seen in the  proof of Lemma \ref{gl-girth6}, if $\rad (G) \leq 2$, then  $\girth (G) \leq 5$. Therefore, we investigate three categories of unicyclic graphs: those whose cycle is $C_5$, those with $C_4$, and those with $C_3$, where each graph in these categories has radius at most two.

Throughout this section, we adopt the following notation. If $u$ is a support vertex in $G$, then the number of leaves adjacent to $u$ is denoted by $n_u$, and the leaves adjacent to $u$ are denoted by
$u_1,u_2, \cdots, u_{n_u}$. Further,  we define $L_u=\{u_1,u_2, \cdots, u_{n_u}\}$. A vertex $s$ outside the cycle which is adjacent to a vertex $v$ and to at least one leaf is called a \emph{strong child} of $v$. In the figures of this section, dotted lines indicate that the number of the corresponding leaves or strong children is arbitrary, but at least one.

All lower bounds of this section are instances of Proposition~\ref{prop-mu}, in the form given by the following lemma.

\begin{lemma}\label{lem-uni-lb}
Let $G$ be a connected graph of order $n\ge 3$ and let $v$ be a vertex of $G$ with $\ecc(v)=2$ such that $N(v)\cup\{w\}$ dominates $\overline{G}$ for every $w\in V(G)\setminus N[v]$. Then $\{N(v)\}\cup\{\{w\}: w\notin N(v)\}$ is a $gc$-partition of $G$, and consequently $GC(G)\ge n-\deg(v)+1$.
\end{lemma}

\begin{proof}
The set $N(v)$ dominates $G$ because $\ecc(v)=2$, but it does not dominate $\overline{G}$, as $v$ has no neighbor in $N(v)$ in $\overline{G}$. The set $N[v]$ dominates $G$ and $\overline{G}$ (in $\overline{G}$, every vertex outside $N[v]$ is adjacent to $v$). For $w\notin N[v]$, the set $N(v)\cup\{w\}$ dominates $G$ and, by assumption, $\overline{G}$. Hence $N(v)$ is a $gc$-partner of every singleton, and no singleton is a global dominating set.
\end{proof}

In each of the families below, the hypothesis of Lemma~\ref{lem-uni-lb} is satisfied by the central vertex $a$ (the vertex denoted by $a$ in Figures~\ref{u5-rad2}--\ref{u3-rad2}). Indeed, $\deg(a)\ge 2$, and a vertex $w'$ at distance two from $a$ has at most two neighbours in $N(a)$, with equality only if $w'$ is the vertex opposite to $a$ on a $4$-cycle, in which case $w'$ has no other neighbours. Thus, for every $w\notin N[a]$, each vertex $w'\ne w$ at distance two from $a$ is non-adjacent to some vertex of $N(a)\cup\{w\}$, which implies that $N(a)\cup\{w\}$ dominates $\overline{G}$.

\subsection{Unicyclic graphs with cycle $C_5$}\label{sub-sec-u-5}
\begin{figure}[h]
	\centering
	\begin{tikzpicture}[scale=.4, transform shape]

        \node [draw, shape=circle,fill=black] (v1) at  (-10,1.5) {};
		\node [draw, shape=circle,fill=black] (v2) at  (-8.5,.5) {};
		\node [draw, shape=circle,fill=black] (v3) at  (-11.5,.5) {};	
		\node [draw, shape=circle,fill=black] (v4) at  (-11.5,-1.5) {};
		\node [draw, shape=circle,fill=black] (v5) at  (-8.5,-1.5) {};

        \node [draw, shape=circle,fill=black] (v6) at  (-11,2.5) {};
		\node [draw, shape=circle,fill=black] (v7) at  (-9,2.5) {};

		\draw (v1)--(v2);
		\draw (v1)--(v3);
		\draw(v4)--(v5);
		\draw (v3)--(v4);
		\draw (v2)--(v5);
        \draw (v1)--(v6);
		\draw (v1)--(v7);
       \node [scale=2] at (-10,2.5) {$\cdots$};
       
       \node [scale=2] at (-10,.8) {$a$};
       \node [scale=2] at (-8,.5) {$b$};
       \node [scale=2] at (-8,-1.5) {$c$};
       \node [scale=2] at (-12,-1.5) {$d$};
		\node [scale=2] at (-12,.5) {$e$};
        \node [scale=2] at (-9.7,-.5) {$U_{5_1}$};


        \node [draw, shape=circle,fill=black] (v1) at  (-3,1.5) {};
		\node [draw, shape=circle,fill=black] (v2) at  (-1.5,.5) {};
		\node [draw, shape=circle,fill=black] (v3) at  (-4.5,.5) {};	
		\node [draw, shape=circle,fill=black] (v4) at  (-4.5,-1.5) {};
		\node [draw, shape=circle,fill=black] (v5) at  (-1.5,-1.5) {};

\node [draw, shape=circle,fill=black] (v6) at  (-4,2.5) {};
		\node [draw, shape=circle,fill=black] (v7) at  (-2,2.5) {};

		\node [draw, shape=circle,fill=black] (v8) at  (-5.5,1.5) {};
		\node [draw, shape=circle,fill=black] (v9) at  (-5.5,-.5) {};

		\draw (v1)--(v2);
		\draw (v1)--(v3);
		\draw(v4)--(v5);
		\draw (v3)--(v4);
		\draw (v2)--(v5);
        \draw (v1)--(v6);
		\draw (v1)--(v7);
       \draw (v3)--(v8);
       \draw (v3)--(v9);
       \node [scale=2] at (-3,2.5) {$\cdots$};
       \node [scale=2] at (-5.5,.7) {$\vdots$};
		
        \node [scale=2] at (-2.7,-.5) {$U_{5_2}$};


        \node [draw, shape=circle,fill=black] (v1) at  (3,1.5) {};
		\node [draw, shape=circle,fill=black] (v2) at  (4.5,.5) {};
		\node [draw, shape=circle,fill=black] (v3) at  (1.5,.5) {};	
		\node [draw, shape=circle,fill=black] (v4) at  (1.5,-1.5) {};
		\node [draw, shape=circle,fill=black] (v5) at  (4.5,-1.5) {};

\node [draw, shape=circle,fill=black] (v6) at  (2,2.5) {};
		\node [draw, shape=circle,fill=black] (v7) at (4,2.5) {};

		\node [draw, shape=circle,fill=black] (v8) at  (.5,1.5) {};
		\node [draw, shape=circle,fill=black] (v9) at  (.5,-.5) {};
		
\node [draw, shape=circle,fill=black] (v10) at  (5.5,1.5) {};
		\node [draw, shape=circle,fill=black] (v11) at  (5.5,-.5) {};

		\draw (v1)--(v2);
		\draw (v1)--(v3);
		\draw(v4)--(v5);
		\draw (v3)--(v4);
		\draw (v2)--(v5);

\draw (v1)--(v6);
		\draw (v1)--(v7);
       \draw (v3)--(v8);
       \draw (v3)--(v9);
       
        \draw (v2)--(v10);
       \draw (v2)--(v11);
       
       \node [scale=2] at (3,2.5) {$\cdots$};
       \node [scale=2] at (.5,.7) {$\vdots$};
       
       \node [scale=2] at (5.5,.7) {$\vdots$};

        \node [scale=2] at (3.3,-.5) {$U_{5_3}$};



		\node [draw, shape=circle,fill=black] (v1) at  (10,1.5) {};
		\node [draw, shape=circle,fill=black] (v2) at  (11.5,.5) {};
		\node [draw, shape=circle,fill=black] (v3) at  (8.5,.5) {};	
		\node [draw, shape=circle,fill=black] (v4) at  (8.5,-1.5) {};
		\node [draw, shape=circle,fill=black] (v5) at  (11.5,-1.5) {};

        \node [draw, shape=circle,fill=black] (v8) at  (7.5,1.5) {};
		\node [draw, shape=circle,fill=black] (v9) at  (7.5,-.5) {};
		
        \node [draw, shape=circle,fill=black] (v10) at  (12.5,1.5) {};
		\node [draw, shape=circle,fill=black] (v11) at  (12.5,-.5) {};

		\draw (v1)--(v2);
		\draw (v1)--(v3);
		\draw(v4)--(v5);
		\draw (v3)--(v4);
		\draw (v2)--(v5);
		
\draw (v3)--(v8);
       \draw (v3)--(v9);
       
        \draw (v2)--(v10);
       \draw (v2)--(v11);
       
       \node [scale=2] at (7.5,.7) {$\vdots$};
       
       \node [scale=2] at (12.5,.7) {$\vdots$};

        \node [scale=2] at (10.3,-.5) {$U_{5_4}$};



        \node [draw, shape=circle,fill=black] (v1) at  (17,1.5) {};
		\node [draw, shape=circle,fill=black] (v2) at  (18.5,.5) {};
		\node [draw, shape=circle,fill=black] (v3) at  (15.5,.5) {};	
		\node [draw, shape=circle,fill=black] (v4) at  (15.5,-1.5) {};
		\node [draw, shape=circle,fill=black] (v5) at  (18.5,-1.5) {};

        \node [draw, shape=circle,fill=black] (v6) at  (16,2.5) {};
		\node [draw, shape=circle,fill=black] (v7) at  (18,2.5) {};

        \node [draw, shape=circle,fill=black] (v12) at  (15,3.5) {};
		\node [draw, shape=circle,fill=black] (v13) at  (15,1.5) {}; 

        \node [draw, shape=circle,fill=black] (v14) at  (19,3.5) {};
		\node [draw, shape=circle,fill=black] (v15) at  (19,1.5) {};

		\draw (v1)--(v2);
		\draw (v1)--(v3);
		\draw(v4)--(v5);
		\draw (v3)--(v4);
		\draw (v2)--(v5);
        \draw (v1)--(v6);
		\draw (v1)--(v7);

       \draw (v6)--(v12);
       \draw (v6)--(v13);
        \draw (v7)--(v14);
       \draw (v7)--(v15);
       \node [scale=2] at (16.5,2.8) {$f$};
       \node [scale=2] at (17.5,2.8) {$h$};

       \node [scale=2] at (17,3.5) {$\cdots$};
       
       \node [scale=2] at (15,2.8) {$\vdots$};
       
       \node [scale=2] at (19,2.8) {$\vdots$};
       
       \node [scale=2] at (17,.8) {$a$};
		
        \node [scale=2] at (17.3,-.5) {$U_{5_5}$};
        


        \node [draw, shape=circle,fill=black] (v1) at  (-10,-5) {};
		\node [draw, shape=circle,fill=black] (v2) at  (-8.5,-6) {};
		\node [draw, shape=circle,fill=black] (v3) at  (-11.5,-6) {};	
		\node [draw, shape=circle,fill=black] (v4) at  (-11.5,-8) {};
		\node [draw, shape=circle,fill=black] (v5) at  (-8.5,-8) {};

        \node [draw, shape=circle,fill=black] (v6) at  (-11,-4) {};
		\node [draw, shape=circle,fill=black] (v7) at  (-9,-4) {};

        \node [draw, shape=circle,fill=black] (v12) at  (-12,-3) {};
		\node [draw, shape=circle,fill=black] (v13) at  (-12,-5) {}; 

        \node [draw, shape=circle,fill=black] (v14) at  (-8,-3) {};
		\node [draw, shape=circle,fill=black] (v15) at  (-8,-5) {}; 

        \node [draw, shape=circle,fill=black] (v16) at  (-9.2,-2.7) {};
		\node [draw, shape=circle,fill=black] (v17) at  (-10.8,-2.7) {}; 
		
		\draw (v1)--(v2);
		\draw (v1)--(v3);
		\draw(v4)--(v5);
		\draw (v3)--(v4);
		\draw (v2)--(v5);
        \draw (v1)--(v6);
		\draw (v1)--(v7);

       \draw (v6)--(v12);
       \draw (v6)--(v13);
        \draw (v7)--(v14);
       \draw (v7)--(v15);

        \draw (v1)--(v16);
		\draw (v1)--(v17);

       \node [scale=3] at (-10,-4) {$\ldots$};
       
       \node [scale=2] at (-12,-3.7) {$\vdots$};
       
       \node [scale=2] at (-8,-3.7) {$\vdots$};
		
        \node [scale=2] at (-9.7,-7) {$U_{5_6}$};
        
        \node [scale=2] at (-10,-2.7) {$\cdots$};
        
         \node [scale=2] at (-10,-5.7) {$a$};



        \node [draw, shape=circle,fill=black] (v1) at  (-3,-5) {};
		\node [draw, shape=circle,fill=black] (v2) at  (-1.5,-6) {};
		\node [draw, shape=circle,fill=black] (v3) at  (-4.5,-6) {};	
		\node [draw, shape=circle,fill=black] (v4) at  (-4.5,-8) {};
		\node [draw, shape=circle,fill=black] (v5) at  (-1.5,-8) {};

        \node [draw, shape=circle,fill=black] (v6) at  (-4,-4) {};
		\node [draw, shape=circle,fill=black] (v7) at  (-2,-4) {};

        \node [draw, shape=circle,fill=black] (v12) at  (-5,-3) {};
		\node [draw, shape=circle,fill=black] (v13) at  (-5,-5) {}; 

        \node [draw, shape=circle,fill=black] (v14) at  (-1,-3) {};
		\node [draw, shape=circle,fill=black] (v15) at  (-1,-5) {}; 

         \node [draw, shape=circle,fill=black] (v18) at  (-6,-5) {};
		 \node [draw, shape=circle,fill=black] (v19) at  (-6,-7) {};

		\draw (v1)--(v2);
		\draw (v1)--(v3);
		\draw(v4)--(v5);
		\draw (v3)--(v4);
		\draw (v2)--(v5);
        \draw (v1)--(v6);
		\draw (v1)--(v7);

       \draw (v6)--(v12);
       \draw (v6)--(v13);
        \draw (v7)--(v14);
       \draw (v7)--(v15);
       
       \draw (v3)--(v18);
       \draw (v3)--(v19);

       \node [scale=2] at (-3,-3.3) {$\cdots$};
       
       \node [scale=2] at (-5,-3.7) {$\vdots$};
       
       \node [scale=2] at (-1,-3.7) {$\vdots$};
		
        \node [scale=2] at (-2.7,-7) {$U_{5_7}$};
        
        \node [scale=2] at (-6,-5.8) {$\vdots$};
        
         \node [scale=2] at (-3,-5.7) {$a$};


        \node [draw, shape=circle,fill=black] (v1) at  (4,-5) {};
		\node [draw, shape=circle,fill=black] (v2) at  (5.5,-6) {};
		\node [draw, shape=circle,fill=black] (v3) at  (2.5,-6) {};	
		\node [draw, shape=circle,fill=black] (v4) at  (2.5,-8) {};
		\node [draw, shape=circle,fill=black] (v5) at  (5.5,-8) {};

        \node [draw, shape=circle,fill=black] (v6) at  (3,-4) {};
		\node [draw, shape=circle,fill=black] (v7) at  (5,-4) {};

        \node [draw, shape=circle,fill=black] (v12) at  (2,-3) {};
		\node [draw, shape=circle,fill=black] (v13) at  (2,-5) {}; 

        \node [draw, shape=circle,fill=black] (v14) at  (6,-3) {};
		\node [draw, shape=circle,fill=black] (v15) at  (6,-5) {}; 

        \node [draw, shape=circle,fill=black] (v16) at  (4.8,-2.7) {};
		\node [draw, shape=circle,fill=black] (v17) at  (3.2,-2.7) {}; 

         \node [draw, shape=circle,fill=black] (v20) at  (1,-5) {};
		 \node [draw, shape=circle,fill=black] (v21) at  (1,-7) {}; 
		
		\draw (v1)--(v2);
		\draw (v1)--(v3);
		\draw(v4)--(v5);
		\draw (v3)--(v4);
		\draw (v2)--(v5);
        \draw (v1)--(v6);
		\draw (v1)--(v7);

       \draw (v6)--(v12);
       \draw (v6)--(v13);
        \draw (v7)--(v14);
       \draw (v7)--(v15);

        \draw (v1)--(v16);
		\draw (v1)--(v17);

       \draw (v3)--(v21);
		\draw (v3)--(v20);
        
       \node [scale=3] at (4,-4) {$\ldots$};
       
       \node [scale=2] at (2,-3.7) {$\vdots$};
       
       \node [scale=2] at (6,-3.7) {$\vdots$};
		
        \node [scale=2] at (4.3,-7) {$U_{5_8}$};
        
        \node [scale=2] at (4,-2.7) {$\cdots$};
        
        \node [scale=2] at (1,-5.8) {$\vdots$};
      
         \node [scale=2] at (4,-5.7) {$a$};



        \node [draw, shape=circle,fill=black] (v1) at  (11,-5) {};
		\node [draw, shape=circle,fill=black] (v2) at  (12.5,-6) {};
		\node [draw, shape=circle,fill=black] (v3) at  (9.5,-6) {};	
		\node [draw, shape=circle,fill=black] (v4) at  (9.5,-8) {};
		\node [draw, shape=circle,fill=black] (v5) at  (12.5,-8) {};

        \node [draw, shape=circle,fill=black] (v6) at  (10,-4) {};
		\node [draw, shape=circle,fill=black] (v7) at  (12,-4) {};

        \node [draw, shape=circle,fill=black] (v12) at  (9,-3) {};
		\node [draw, shape=circle,fill=black] (v13) at  (9,-5) {}; 

        \node [draw, shape=circle,fill=black] (v14) at  (13,-3) {};
		\node [draw, shape=circle,fill=black] (v15) at  (13,-5) {}; 

         \node [draw, shape=circle,fill=black] (v18) at  (8,-5) {};
		 \node [draw, shape=circle,fill=black] (v19) at  (8,-7) {};

        \node [draw, shape=circle,fill=black] (v22) at  (14,-5) {};
		 \node [draw, shape=circle,fill=black] (v23) at  (14,-7) {}; 
		
		\draw (v1)--(v2);
		\draw (v1)--(v3);
		\draw(v4)--(v5);
		\draw (v3)--(v4);
		\draw (v2)--(v5);
        \draw (v1)--(v6);
		\draw (v1)--(v7);

       \draw (v6)--(v12);
       \draw (v6)--(v13);
        \draw (v7)--(v14);
       \draw (v7)--(v15);
       
       \draw (v3)--(v18);
       \draw (v3)--(v19);
       
       \draw (v2)--(v22);
       \draw (v2)--(v23);

       \node [scale=2] at (11,-3.3) {$\cdots$};
       
       \node [scale=2] at (9,-3.7) {$\vdots$};
       
       \node [scale=2] at (13,-3.7) {$\vdots$};
		
        \node [scale=2] at (11.3,-7) {$U_{5_9}$};
        
        \node [scale=2] at (8,-5.8) {$\vdots$};
        
        \node [scale=2] at (14,-5.8) {$\vdots$};
        
         \node [scale=2] at (11,-5.7) {$a$};
        

        \node [draw, shape=circle,fill=black] (v1) at  (19,-5) {};
		\node [draw, shape=circle,fill=black] (v2) at  (20.5,-6) {};
		\node [draw, shape=circle,fill=black] (v3) at  (17.5,-6) {};	
		\node [draw, shape=circle,fill=black] (v4) at  (17.5,-8) {};
		\node [draw, shape=circle,fill=black] (v5) at  (20.5,-8) {};

        \node [draw, shape=circle,fill=black] (v6) at  (18,-4) {};
		\node [draw, shape=circle,fill=black] (v7) at  (20,-4) {};

        \node [draw, shape=circle,fill=black] (v12) at  (17,-3) {};
		\node [draw, shape=circle,fill=black] (v13) at  (17,-5) {}; 

        \node [draw, shape=circle,fill=black] (v14) at  (21,-3) {};
		\node [draw, shape=circle,fill=black] (v15) at  (21,-5) {}; 

        \node [draw, shape=circle,fill=black] (v16) at  (19.8,-2.7) {};
		\node [draw, shape=circle,fill=black] (v17) at  (18.2,-2.7) {}; 

         \node [draw, shape=circle,fill=black] (v20) at  (16,-5) {};
		 \node [draw, shape=circle,fill=black] (v21) at  (16,-7) {}; 

         \node [draw, shape=circle,fill=black] (v24) at  (22,-5) {};
		 \node [draw, shape=circle,fill=black] (v25) at  (22,-7) {}; 
		
		\draw (v1)--(v2);
		\draw (v1)--(v3);
		\draw(v4)--(v5);
		\draw (v3)--(v4);
		\draw (v2)--(v5);
        \draw (v1)--(v6);
		\draw (v1)--(v7);

       \draw (v6)--(v12);
       \draw (v6)--(v13);
        \draw (v7)--(v14);
       \draw (v7)--(v15);

        \draw (v1)--(v16);
		\draw (v1)--(v17);

       \draw (v3)--(v21);
		\draw (v3)--(v20);

       \draw (v2)--(v24);
       \draw (v2)--(v25);
        
       \node [scale=3] at (19,-4) {$\ldots$};
       
       \node [scale=2] at (17,-3.7) {$\vdots$};
       
       \node [scale=2] at (21,-3.7) {$\vdots$};
		
        \node [scale=2] at (19.3,-7) {$U_{5_{10}}$};
        
        \node [scale=2] at (19,-2.7) {$\cdots$};
        
        \node [scale=2] at (16,-5.8) {$\vdots$};

         \node [scale=2] at (22,-5.8) {$\vdots$};
         
          \node [scale=2] at (19,-5.7) {$a$};

	\end{tikzpicture}

	\caption{Unicyclic graphs with cycle $C_5$ and radius two.}
	\label{u5-rad2}
\end{figure}
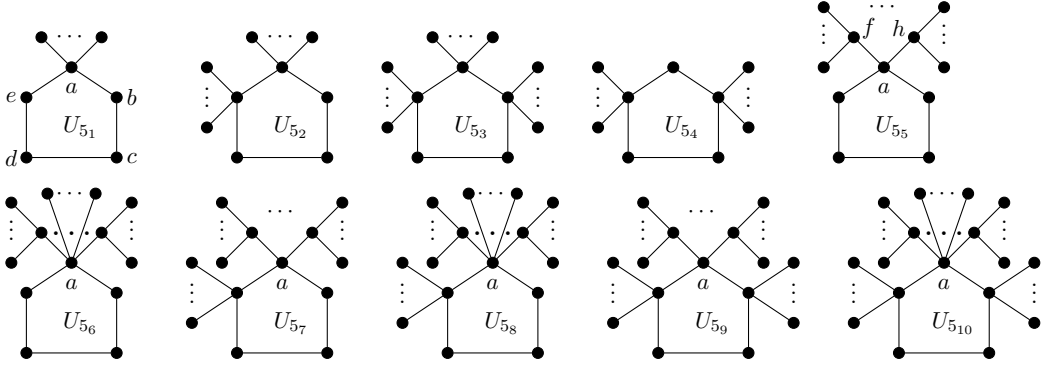

In this subsection, we study the class of unicyclic graphs with order at least six and radius two, whose cycle is $C_5$. The unicyclic graphs that satisfy these conditions are depicted in Figure \ref{u5-rad2}. Let $C=(a,b,c,d,e)$ be the cycle and let $a$ be a central vertex. In the families $U_{5_1}$--$U_{5_4}$ the vertex $a$ has no strong children, and the sets of support vertices are $\{a\}$, $\{a,e\}$, $\{a,b,e\}$ and $\{b,e\}$, respectively. In the families $U_{5_5}$--$U_{5_{10}}$ the vertex $a$ has at least one strong child, and the cycle vertices adjacent to leaves are: none ($U_{5_5}$), $a$ ($U_{5_6}$), $e$ ($U_{5_7}$), $a$ and $e$ ($U_{5_8}$), $b$ and $e$ ($U_{5_9}$), and $a$, $b$ and $e$ ($U_{5_{10}}$).

Given the family of graphs presented in Figure \ref{u5-rad2} and the explanations given at the beginning of this section, we have the following observation.

\begin{observation}\label{obs-u5}
	A graph  $G$ of order at least six is a unicyclic graph containing the cycle $C_5$ and satisfying $\rad(G)\le 2$ if and only if    $G \in U_{5_i}$ for some $1\le i \le 10$. (The notation $U_{5_i}$ denotes the family of graphs shown in Figure \ref{u5-rad2}.)
\end{observation}
\begin{proof}
Every graph in the families $U_{5_i}$ is unicyclic with cycle $C_5$, and $a$ has eccentricity two. Conversely, let $G$ be unicyclic with cycle $C=(a,b,c,d,e)$, of order at least six and radius at most two, and let $v$ be a central vertex. A vertex $x\notin V(C)$ is joined to $C$ by a unique path ending at some vertex of $C$, and the two cycle vertices at distance two from that vertex are at distance at least three from $x$. Hence $v\in V(C)$, say $v=a$. Every vertex is at distance two from $a$, so the pendant trees attached at $c$ and $d$ are trivial, the pendant trees attached at $b$ and $e$ consist of leaves only, and the pendant tree at $a$ has depth at most two, i.e., $a$ is adjacent to leaves and to strong children only. Enumerating the possibilities (whether $a$ has strong children, and which of $a$, $b$, $e$ are support vertices) up to the symmetry exchanging $b$ with $e$ gives exactly the ten families of Figure~\ref{u5-rad2}. Note that a graph with leaves at $b$ only is isomorphic to a graph in $U_{5_1}$.
\end{proof}

In what follows, we determine the exact values of the global coalition number for some classes of graphs introduced above, and establish bounds for others. To aid our discussion, we assume that $G$ is a unicyclic graph of order $n$ with cycle $C_5 =(a,b,c,d,e)$ and radius two.

\begin{proposition}\label{prop-u5-1}
	If $G \in U_{5_1}\cup U_{5_4}$, then $GC(G)=n-1$.
\end{proposition}

\begin{proof}
Let $G \in U_{5_1}$. Let $\pi$ be a $GC(G)$-partition.  Since $G$ has no $\gamma_g$-set consisting of two vertices (equivalently, by Theorem~\ref{thm-antitwin}, some vertex has no antitwin; for instance, $c$ has none),  it follows that $\pi$ contains a partition class with at least two vertices. Thus, $\lvert \pi \rvert \leq n-1$.   On the other hand, consider the vertex partition $\phi$ of $G$, where $\{a,d\} \in \phi$ and the other partition classes are singleton sets.  Then $\phi$ is a $gc$-partition of $G$, where any singleton set is a $gc$-partner of $\lbrace a,d\rbrace$. Thus, $\lvert \pi\rvert \geq n-1$, and so $GC(G)=n-1$.
\vspace{1mm}

Let $G \in U_{5_4}$. Let $\pi$ be a $GC(G)$-partition.  Since $G$ has no $\gamma_g$-set with two vertices,   $\pi$ contains a partition class with at least two vertices. Thus, $\lvert \pi \rvert \leq n-1$.   On the other hand, consider the vertex partition $\phi$ of $G$, where $\lbrace b,e\rbrace \in \phi$ is the set of support vertices of $G$, and the other partition classes are singleton sets.  Then $\phi$ is a $gc$-partition of $G$, where any singleton set is a $gc$-partner of $\lbrace b,e\rbrace$. Thus, $\lvert \pi\rvert \geq n-1$, and so $GC(G)=n-1$.
\end{proof}

\begin{theorem}\label{theo-u5-1}
	If $G \in U_{5_2}$, then $GC(G)=\max \{4+n_a, 4+n_e\}$, where the vertices $a$ and $e$ are the support vertices of $G$.
\end{theorem}
\begin{proof}
Assume, by symmetry, that $n_a \geq n_e$. Then the partition   $\lbrace \lbrace a,d\rbrace \cup L_e,\lbrace e\rbrace, \lbrace b\rbrace ,\lbrace c\rbrace ,\lbrace a_i\rbrace : 1 \leq i \leq n_a \rbrace$ is a $gc$-partition of $G$ of size $n_a +4$, where the set $\lbrace a,d\rbrace \cup L_e$ forms a global coalition with all other members. Thus, $GC(G) \geq \max \lbrace n_a +4, n_e +4 \rbrace$. Now we show that $GC(G) \leq \max \lbrace n_a +4, n_e +4\rbrace$. Let $\pi$ be a $gc$-partition of $G$. Note that a set $S\supseteq\{a,e\}$ is a global dominating set if and only if $S\cap\{b,c,d\}\ne\emptyset$. We consider two cases regarding $\pi$:

\textbf{Case 1}. $\lbrace a,e\rbrace \subseteq V_i$ for some $V_i \in \pi$.

If there exists a pair $(A,B)$ of $gc$-partners such that $V_i \notin \lbrace A,B\rbrace$, then we must have $(L_a \cup L_e) \subset (A \cup B)$, and $\lbrace b,c,d\rbrace \cap (A \cup B) \neq \emptyset$, which means that $\lvert \pi\rvert \leq 5$, as desired. Thus, we may assume  that $V_i \in \lbrace A,B\rbrace$ for any pair $(A,B)$ of $gc$-partners in $\pi$. On the other hand, $c \notin V_i$, for otherwise $V_i$ would be a global dominating set. Therefore, each $gc$-partner of $V_i$ must dominate $c$, implying that $V_i$ has at most three $gc$-partners, and so $\lvert \pi \rvert \leq 4$, as desired.

\textbf{Case 2}.  No member of $\pi$ contains $\lbrace a,e\rbrace$.

Let $V_a$ and $V_e$ be the members of $\pi$ containing $a$ and $e$, respectively. Assume first that there exists a pair $(A,B)$ of $gc$-partners such that $\lbrace A,B\rbrace \cap \lbrace V_a,V_e\rbrace =\emptyset$. It follows that $(L_a \cup L_e) \subset (A \cup B)$, and $ (A \cup B) \cap \lbrace b,c,d\rbrace \neq \emptyset$, implying that $\lvert \pi \rvert \leq 6$. Now we show that $\lvert \pi \rvert \neq 6$. Suppose to the contrary that $|\pi|=6$. Then we must have $A \cup B= L_a \cup L_e \cup \lbrace k\rbrace$, where $k \in \lbrace b,c,d\rbrace$, and all other partition classes must be singleton sets. In this case, it is easily seen that no singleton subset of $ \lbrace b,c,d\rbrace \setminus \lbrace k\rbrace$ has a $gc$-partner, a contradiction. Thus, $\lvert \pi \rvert \leq 5$, as desired.  Hence, we may assume that for any pair $(A,B)$ of $gc$-partners, $\lbrace A,B\rbrace \cap \lbrace V_a,V_e\rbrace  \neq \emptyset$.
If the only $gc$-partner of $V_a$ outside $\{V_a,V_e\}$ is also the only such partner of $V_e$, then $|\pi|\le 3$.
Assume next that the sets $V_a$ and $V_e$ have distinct $gc$-partners in $\pi \setminus \lbrace V_a , V_e \rbrace$. Let $C$ and $D$ be $gc$-partners of $V_a$ and $V_e$, respectively. It follows that $L_e \subset (V_a \cup C)$, $L_a \subset (V_e \cup D)$, $ (V_a \cup C) \cap \lbrace b,c,d\rbrace \neq \emptyset$ and $(V_e \cup D) \cap \lbrace b,c,d\rbrace \neq \emptyset$. Consequently, $\lvert \pi \rvert \leq 5$, as desired. Hence, we may assume that one of the sets in $\lbrace V_a, V_e\rbrace$ forms a global coalition with all members of $\pi \setminus \lbrace V_a, V_e \rbrace$, and that there are at least two such members. Assume that $V_a$ has this property. It follows that $L_e \subset V_a$. If $V_a \cap \lbrace b,c,d\rbrace =\emptyset$, then $V_a$ has at most three $gc$-partners in $\pi\setminus\{V_a,V_e\}$, since each $gc$-partner of $V_a$ must dominate $c$. In this case, we have $\lvert \pi \rvert \leq 5$, as desired. Otherwise,  $V_a \cap \lbrace b,c,d\rbrace \neq \emptyset$, which implies that $\lvert \pi \rvert \leq n_a +4$, as desired. Similarly, we can show that if $V_e$ forms a global coalition with all members of $\pi \setminus \lbrace V_a, V_e \rbrace$, then $\lvert \pi \rvert \leq n_e +4$. Hence, $GC(G) \leq \max \lbrace n_a +4, n_e +4\rbrace$, and so $GC(G) = \max \lbrace n_a +4, n_e +4\rbrace$.
\end{proof}

\begin{theorem}\label{theo-u5-2}
	If $G \in U_{5_3}$, then $GC(G)=n_b +n_e +4$, where the vertices $a$, $b$ and $e$ are the support vertices of $G$.
\end{theorem}

\begin{proof}
The partition $\lbrace L_a \cup \lbrace b,e\rbrace , \lbrace a\rbrace ,\lbrace c\rbrace ,\lbrace d\rbrace , \lbrace b_i \rbrace : 1 \leq i \leq n_b , \lbrace e_i \rbrace : 1 \leq i \leq n_e \rbrace$ is a $gc$-partition of $G$ of size $n_b +n_e +4$, where the set $L_a \cup \lbrace b,e\rbrace$ $=N(a)$ forms a global coalition with all other members, so $GC(G) \geq n_b +  n_e +4$. Now we show that $GC(G) \leq n_b +  n_e +4$. Let $\pi$ be a $gc$-partition of $G$. We consider two cases:

\textbf{Case 1}. $\lbrace b,e \rbrace \subseteq V_i$ for some $V_i \in \pi$.

Note that $a \notin V_i$, for otherwise $V_i$ would be a global dominating set. Let $V_a$ be the member of $\pi$ containing $a$. If there exists a pair $(A,B)$ of $gc$-partners such that $\lbrace A,B\rbrace \cap \lbrace V_a,V_i \rbrace =\emptyset$, then we must have $(L_a \cup L_b \cup L_e ) \subset (A \cup B)$, and $\lbrace c,d\rbrace \cap (A \cup B) \neq \emptyset$, implying that $\lvert \pi \rvert \leq 5$, as desired. Hence, we may assume that each set in $\pi \setminus \lbrace V_a ,V_i \rbrace$ is a $gc$-partner of $V_a$ or $V_i$. Now we show that $V_a$ has at most two $gc$-partners in $\pi \setminus \lbrace V_i\rbrace$. Suppose, to the contrary, that it has at least three. Every $gc$-partner $C\ne V_i$ of $V_a$ satisfies $L_b\cup L_e\subseteq V_a\cup C$. If $L_b\cup L_e\not\subseteq V_a$, all such partners contain the non-empty set $(L_b\cup L_e)\setminus V_a$, so there is at most one of them. It follows that $(L_b \cup L_e) \subset V_a$. Thus, $\lbrace c,d\rbrace \cap V_a = \emptyset$, for otherwise $V_a$ would be a global dominating set. Now each $gc$-partner of $V_a$ in $\pi \setminus \lbrace V_i\rbrace$ must contain $c$ or $d$, which contradicts the assumption. Hence, $V_a$ has at most two $gc$-partners in $\pi \setminus \lbrace V_i\rbrace$. Assume first that $V_a$ has two $gc$-partners in $\pi \setminus \lbrace V_i\rbrace$. Let $A$ and $B$ be the $gc$-partners of $V_a$. It follows as above that $(L_b \cup L_e) \subset V_a$, $\{c,d\}\cap V_a=\emptyset$, and either $c \in  A$ and $d \in B$, or $c \in B$ and $d \in A$. If $L_a \subset V_i$, then $\lvert \pi \rvert =4$, as desired. Otherwise, $V_i$ has at most one $gc$-partner in $\pi \setminus \lbrace V_a\rbrace$ (such a partner must contain $L_a\setminus V_i$), which implies that $\lvert \pi \rvert \leq 5$, as desired. Now assume that $V_a$ has one $gc$-partner in $\pi \setminus \lbrace V_i\rbrace$. Let $A$ be the $gc$-partner of $V_a$. It follows that $(L_b \cup L_e) \subset (A \cup V_a)$. If $L_a \subset V_i$, then clearly $\lvert \pi \rvert \leq 5$. Otherwise $V_i$ has at most one $gc$-partner in $\pi \setminus \lbrace V_a\rbrace$, since every $gc$-partner of $V_i$ must dominate $L_a$. Therefore, since each set in $\pi \setminus \lbrace V_a ,V_i \rbrace$ is a $gc$-partner of $V_a$ or $V_i$, we have $\lvert \pi \rvert \leq 5$, as desired. Hence, we may assume that $V_a$ has no $gc$-partner in $\pi \setminus \lbrace V_i\rbrace$, implying that each set in $\pi \setminus \lbrace V_a ,V_i \rbrace$ is a $gc$-partner of $V_i$. If $L_a \nsubseteq V_i$, then $V_i$ has at most one $gc$-partner in $\pi \setminus \lbrace V_a\rbrace$, implying that $\lvert \pi \rvert \leq 3$. Hence, we may assume $L_a \subset V_i$, which implies that $\lvert \pi \rvert \leq n_b + n_e +4$, as desired.

\textbf{Case 2}. No member of $\pi$ contains $\lbrace b,e\rbrace$.

Suppose first that $\{a,e\}\subseteq V_j$ for some $V_j\in\pi$. A set $S\supseteq\{a,e\}$ is a global dominating set if and only if ($b\in S$ or $L_b\subseteq S$) and $S\cap\{b,c,d\}\ne\emptyset$. If some pair $(A,B)$ of $gc$-partners avoids $V_j$, then $A\cup B$ contains $L_a\cup L_e$, meets $\{c,d\}$ and contains $b$ or $L_b$, so $|A\cup B|\ge n_a+n_e+2$ and $|\pi|\le 2+(n-n_a-n_e-2)=n_b+5\le n_b+n_e+4$, as desired. Otherwise, every class other than $V_j$ is a $gc$-partner of $V_j$. If $V_j\cap\{b,c,d\}=\emptyset$, each partner meets $\{b,c,d\}$ and $|\pi|\le 4$. Otherwise $b\notin V_j$ and $L_b\not\subseteq V_j$ (as $V_j$ is not a global dominating set), so each partner contains $b$ or the non-empty set $L_b\setminus V_j$, whence $|\pi|\le 3$. The case $\{a,b\}\subseteq V_j$ is symmetric. Hence, we may assume that the vertices $a$, $b$ and $e$ are in distinct partition classes.
Let $V_a$, $V_b$ and $V_e$ be the members of $\pi$ containing $a$, $b$ and $e$, respectively. If a pair $(A,B)$ of $gc$-partners avoids $\{V_a,V_b,V_e\}$, then $A\cup B\supseteq L_a\cup L_b\cup L_e$ and $A\cup B$ meets $\{c,d\}$, so $|\pi|\le 2+4=6$, as desired. Hence every class outside $\{V_a,V_b,V_e\}$ is a $gc$-partner of one of $V_a,V_b,V_e$. If  there exists a member $A \in \pi \setminus \lbrace V_a,V_b,V_e\rbrace$ having at least two $gc$-partners in $\lbrace V_a,V_b,V_e\rbrace$, say $V_a$ and $V_b$, then $L_e\subseteq (V_a\cup A)\cap(V_b\cup A)=A$, $L_b\subseteq V_a\cup A$ and $L_a\subseteq V_b\cup A$, so $V_a\cup V_b\cup A$ contains all vertices except possibly $c$, $d$ and $e$, and $|\pi|\le 6$. The other choices of two partners are analogous. Hence, we may assume that every class outside $\{V_a,V_b,V_e\}$ has exactly one $gc$-partner in $\{V_a,V_b,V_e\}$. Let $\mathcal{A}$, $\mathcal{B}$ and $\mathcal{C}$ be the sets of classes outside $\{V_a,V_b,V_e\}$ which are $gc$-partners of $V_a$, $V_b$ and $V_e$, respectively, so that $|\pi|=3+|\mathcal{A}|+|\mathcal{B}|+|\mathcal{C}|$. For $A\in\mathcal{A}$, the set $V_a\cup A$ contains neither $b$ nor $e$, hence it contains $L_b\cup L_e$ and meets $\{c,d\}$. If $L_b\cup L_e\not\subseteq V_a$, then every $A\in\mathcal{A}$ contains the non-empty set $(L_b\cup L_e)\setminus V_a$, so $|\mathcal{A}|\le 1$. If $L_b\cup L_e\subseteq V_a$, then $c,d\notin V_a$ (otherwise $V_a$ would be a global dominating set), each $A\in\mathcal{A}$ contains $c$ or $d$, so $|\mathcal{A}|\le 2$. Moreover, in this case $\mathcal{B}=\mathcal{C}=\emptyset$, because a class $B\in\mathcal{B}$ would satisfy $L_e\subseteq V_b\cup B$, which is impossible as $L_e\subseteq V_a$, and similarly for $\mathcal{C}$. The symmetric statements hold for $\mathcal{B}$ (with $L_a\cup L_e$) and for $\mathcal{C}$ (with $L_a\cup L_b$). Consequently, either one of the three sets has two elements and the other two are empty, or each of them has at most one element. In both cases $|\pi|\le 6\le n_b+n_e+4$. This completes the proof.
\end{proof}

\begin{proposition}\label{prop-u5-2}
	If $G \in \bigcup_{i=5} ^{10} U_{5_i}$, then $GC(G) \geq n-\lvert N(a) \rvert +1$, where $a$ is the central vertex of $G$.
\end{proposition}
\begin{proof}
Consider the vertex partition $\pi$ of $G$, where $N(a) \in \pi$ and all other partition classes are singleton sets. Then $\pi$ is a $gc$-partition of $G$ of size $n- \lvert N(a) \rvert +1$, where the set $N(a)$ forms a global coalition with all other partition classes (Lemma~\ref{lem-uni-lb}). Hence,  $GC(G)\ge n- \lvert N(a) \rvert +1$.
\end{proof}
\subsection{Unicyclic graphs with cycle $C_4$}\label{sub-sec-u-4}
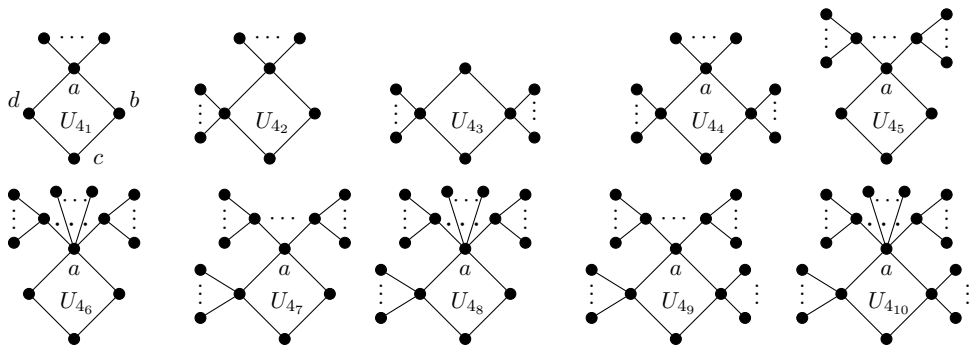
\begin{figure}[h]
	\centering
	\begin{tikzpicture}[scale=.4, transform shape]

		\node [draw, shape=circle,fill=black] (v1) at  (-45,1.5) {};
		\node [draw, shape=circle,fill=black] (v2) at  (-43.5,0) {};
		\node [draw, shape=circle,fill=black] (v3) at  (-46.5,0) {};	
		\node [draw, shape=circle,fill=black] (v4) at  (-45,-1.5) {};
		
		\node [draw, shape=circle,fill=black] (v5) at  (-46,2.5) {};
		\node [draw, shape=circle,fill=black] (v6) at  (-44,2.5) {};

		\draw (v1)--(v2);
		\draw (v1)--(v3);
		\draw(v4)--(v3);
		\draw (v2)--(v4);
		\draw (v1)--(v5);
		\draw (v1)--(v6);
		
		\node [scale=2] at (-45,2.5) {$\cdots$};
		
		\node [scale=2] at (-45,.8) {$a$};
		\node [scale=2] at (-43,.5) {$b$};
		\node [scale=2] at (-44.2,-1.5) {$c$};
		
		\node [scale=2] at (-47,.5) {$d$};
		\node [scale=2] at (-44.9,-.3) {$U_{4_1}$};
		
		
		\node [draw, shape=circle,fill=black] (v1) at  (-38.5,1.5) {};
		\node [draw, shape=circle,fill=black] (v2) at  (-37,0) {};
		\node [draw, shape=circle,fill=black] (v3) at  (-40,0) {};	
		\node [draw, shape=circle,fill=black] (v4) at  (-38.5,-1.5) {};

		\node [draw, shape=circle,fill=black] (v5) at  (-39.5,2.5) {};
		\node [draw, shape=circle,fill=black] (v6) at  (-37.5,2.5) {};

		\node [draw, shape=circle,fill=black] (v7) at  (-40.8,.8) {};
		\node [draw, shape=circle,fill=black] (v8) at  (-40.8,-.8) {};
		
		\draw (v1)--(v2);
		\draw (v1)--(v3);
		\draw(v4)--(v3);
		\draw (v2)--(v4);
		\draw (v1)--(v5);
		\draw (v1)--(v6);
		
		\draw (v3)--(v7);
		\draw (v3)--(v8);
		
		\node [scale=2] at (-38.5,2.5) {$\cdots$};
		
		\node [scale=2] at (-40.8,.2) {$\vdots$};

		\node [scale=2] at (-38.4,-.3) {$U_{4_2}$};
		
		
		\node [draw, shape=circle,fill=black] (v1) at  (-32,1.5) {};
		\node [draw, shape=circle,fill=black] (v2) at  (-30.5,0) {};
		\node [draw, shape=circle,fill=black] (v3) at  (-33.5,0) {};	
		\node [draw, shape=circle,fill=black] (v4) at  (-32,-1.5) {};

		\node [draw, shape=circle,fill=black] (v9) at  (-29.7,.8) {};
		\node [draw, shape=circle,fill=black] (v10) at  (-29.7,-.8) {};

		\node [draw, shape=circle,fill=black] (v7) at  (-34.3,.8) {};
		\node [draw, shape=circle,fill=black] (v8) at  (-34.3,-.8) {};
		
		\draw (v1)--(v2);
		\draw (v1)--(v3);
		\draw(v4)--(v3);
		\draw (v2)--(v4);
		\draw (v2)--(v9);
		\draw (v2)--(v10);
		
		\draw (v3)--(v7);
		\draw (v3)--(v8);
		
		\node [scale=2] at (-29.7,.3) {$\vdots$};
		
		\node [scale=2] at (-34.3,.2) {$\vdots$};

		\node [scale=2] at (-31.9,-.3) {$U_{4_3}$};
		
		
		\node [draw, shape=circle,fill=black] (v1) at  (-24,1.5) {};
		\node [draw, shape=circle,fill=black] (v2) at  (-22.5,0) {};
		\node [draw, shape=circle,fill=black] (v3) at  (-25.5,0) {};	
		\node [draw, shape=circle,fill=black] (v4) at  (-24,-1.5) {};
		
		\node [draw, shape=circle,fill=black] (v5) at  (-25,2.5) {};
		\node [draw, shape=circle,fill=black] (v6) at  (-23,2.5) {};
		
		\node [draw, shape=circle,fill=black] (v9) at  (-21.7,.8) {};
		\node [draw, shape=circle,fill=black] (v10) at  (-21.7,-.8) {};

		\node [draw, shape=circle,fill=black] (v7) at  (-26.3,.8) {};
		\node [draw, shape=circle,fill=black] (v8) at  (-26.3,-.8) {};
		
		\draw (v1)--(v2);
		\draw (v1)--(v3);
		\draw(v4)--(v3);
		\draw (v2)--(v4);
		\draw (v2)--(v9);
		\draw (v2)--(v10);
		
		\draw (v1)--(v5);
		\draw (v1)--(v6);

		\draw (v3)--(v7);
		\draw (v3)--(v8);
		
		\node [scale=2] at (-24,2.5) {$\cdots$};
		
		\node [scale=2] at (-21.7,.3) {$\vdots$};
		
		\node [scale=2] at (-26.3,.2) {$\vdots$};

		\node [scale=2] at (-23.9,-.3) {$U_{4_4}$};
		
		\node [scale=2] at (-24,0.8) {$a$};
		
		\node [draw, shape=circle,fill=black] (v1) at  (-18,1.5) {};
		\node [draw, shape=circle,fill=black] (v2) at  (-16.5,0) {};
		\node [draw, shape=circle,fill=black] (v3) at  (-19.5,0) {};	
		\node [draw, shape=circle,fill=black] (v4) at  (-18,-1.5) {};
		
		\node [draw, shape=circle,fill=black] (v5) at  (-19,2.5) {};
		\node [draw, shape=circle,fill=black] (v6) at  (-17,2.5) {};
		
		\node [draw, shape=circle,fill=black] (v11) at  (-20,3.3) {};
		\node [draw, shape=circle,fill=black] (v12) at  (-20,1.7) {};
		
		\node [draw, shape=circle,fill=black] (v13) at  (-16,3.3) {};
		\node [draw, shape=circle,fill=black] (v14) at  (-16,1.7) {};

		\draw (v1)--(v2);
		\draw (v1)--(v3);
		\draw(v4)--(v3);
		\draw (v2)--(v4);
		\draw (v1)--(v5);
		\draw (v1)--(v6);
		
		\draw (v5)--(v11);
		\draw (v5)--(v12);
		
		\draw (v6)--(v13);
		\draw (v6)--(v14);
		
		\node [scale=2] at (-18,2.5) {$\cdots$};
		
		\node [scale=2] at (-20,2.7) {$\vdots$};
		
		\node [scale=2] at (-16,2.7) {$\vdots$};
		
		
		\node [scale=2] at (-17.9,-.3) {$U_{4_5}$};
		
		\node [scale=2] at (-18,0.8) {$a$};
		
		\node [draw, shape=circle,fill=black] (v1) at  (-45,-4.5) {};
		\node [draw, shape=circle,fill=black] (v2) at  (-43.5,-6) {};
		\node [draw, shape=circle,fill=black] (v3) at  (-46.5,-6) {};	
		\node [draw, shape=circle,fill=black] (v4) at  (-45,-7.5) {};
		
		\node [draw, shape=circle,fill=black] (v5) at  (-46,-3.5) {};
		\node [draw, shape=circle,fill=black] (v6) at  (-44,-3.5) {};
		
		\node [draw, shape=circle,fill=black] (v11) at  (-47,-2.7) {};
		\node [draw, shape=circle,fill=black] (v12) at  (-47,-4.3) {};
		
		\node [draw, shape=circle,fill=black] (v13) at  (-43,-2.7) {};
		\node [draw, shape=circle,fill=black] (v14) at  (-43,-4.3) {};

		\node [draw, shape=circle,fill=black] (v15) at  (-45.6,-2.6) {};
		\node [draw, shape=circle,fill=black] (v16) at  (-44.4,-2.6) {};
		
		\draw (v1)--(v2);
		\draw (v1)--(v3);
		\draw(v4)--(v3);
		\draw (v2)--(v4);
		\draw (v1)--(v5);
		\draw (v1)--(v6);
		
		\draw (v5)--(v11);
		\draw (v5)--(v12);
		
		\draw (v6)--(v13);
		\draw (v6)--(v14);
		
		\draw (v1)--(v15);
		\draw (v1)--(v16);
		
		\node [scale=2.7] at (-45,-3.7) {$\cdots$};
		
		\node [scale=2] at (-47,-3.3) {$\vdots$};
		
		\node [scale=2] at (-43,-3.3) {$\vdots$};
		
		\node [scale=2] at (-44.9,-2.9) {$\cdots$};
		
		\node [scale=2] at (-44.9,-6.3) {$U_{4_6}$};
		
		\node [scale=2] at (-45,-5.2) {$a$};
		
		\node [draw, shape=circle,fill=black] (v1) at  (-38,-4.5) {};
		\node [draw, shape=circle,fill=black] (v2) at  (-36.5,-6) {};
		\node [draw, shape=circle,fill=black] (v3) at  (-39.5,-6) {};	
		\node [draw, shape=circle,fill=black] (v4) at  (-38,-7.5) {};

		\node [draw, shape=circle,fill=black] (v5) at  (-39,-3.5) {};
		\node [draw, shape=circle,fill=black] (v6) at  (-37,-3.5) {};
		
		\node [draw, shape=circle,fill=black] (v11) at  (-40,-2.7) {};
		\node [draw, shape=circle,fill=black] (v12) at  (-40,-4.3) {};
		
		\node [draw, shape=circle,fill=black] (v13) at  (-36,-2.7) {};
		\node [draw, shape=circle,fill=black] (v14) at  (-36,-4.3) {};
		
		\node [draw, shape=circle,fill=black] (v7) at  (-40.8,-5.2) {};
		\node [draw, shape=circle,fill=black] (v8) at  (-40.8,-6.8) {};

		\draw (v1)--(v2);
		\draw (v1)--(v3);
		\draw(v4)--(v3);
		\draw (v2)--(v4);
		\draw (v1)--(v5);
		\draw (v1)--(v6);
		
		\draw (v5)--(v11);
		\draw (v5)--(v12);
		
		\draw (v6)--(v13);
		\draw (v6)--(v14);
		
		\draw (v3)--(v7);
		\draw (v3)--(v8);

		\node [scale=2] at (-40.8,-5.8) {$\vdots$};
		
		\node [scale=2] at (-38,-3.5) {$\cdots$};
		
		\node [scale=2] at (-40,-3.3) {$\vdots$};
		
		\node [scale=2] at (-36,-3.3) {$\vdots$};
		
		\node [scale=2] at (-37.9,-6.3) {$U_{4_7}$};
		
		\node [scale=2] at (-38,-5.2) {$a$};
		
		\node [draw, shape=circle,fill=black] (v1) at  (-32,-4.5) {};
		\node [draw, shape=circle,fill=black] (v2) at  (-30.5,-6) {};
		\node [draw, shape=circle,fill=black] (v3) at  (-33.5,-6) {};	
		\node [draw, shape=circle,fill=black] (v4) at  (-32,-7.5) {};

		\node [draw, shape=circle,fill=black] (v5) at  (-33,-3.5) {};
		\node [draw, shape=circle,fill=black] (v6) at  (-31,-3.5) {};
		
		\node [draw, shape=circle,fill=black] (v11) at  (-34,-2.7) {};
		\node [draw, shape=circle,fill=black] (v12) at  (-34,-4.3) {};
		
		\node [draw, shape=circle,fill=black] (v13) at  (-30,-2.7) {};
		\node [draw, shape=circle,fill=black] (v14) at  (-30,-4.3) {};

		\node [draw, shape=circle,fill=black] (v15) at  (-32.6,-2.6) {};
		\node [draw, shape=circle,fill=black] (v16) at  (-31.4,-2.6) {};
		
		\node [draw, shape=circle,fill=black] (v7) at  (-34.8,-5.2) {};
		\node [draw, shape=circle,fill=black] (v8) at  (-34.8,-6.8) {};
		
		\draw (v1)--(v2);
		\draw (v1)--(v3);
		\draw(v4)--(v3);
		\draw (v2)--(v4);
		\draw (v1)--(v5);
		\draw (v1)--(v6);
		
		\draw (v5)--(v11);
		\draw (v5)--(v12);
		
		\draw (v6)--(v13);
		\draw (v6)--(v14);
		
		\draw (v1)--(v15);
		\draw (v1)--(v16);
		
		\draw (v3)--(v7);
		\draw (v3)--(v8);
		
		\node [scale=2] at (-34.8,-5.8) {$\vdots$};
		
		\node [scale=2.7] at (-32,-3.7) {$\cdots$};
		
		\node [scale=2] at (-34,-3.3) {$\vdots$};
		
		\node [scale=2] at (-30,-3.3) {$\vdots$};
		
		\node [scale=2] at (-31.9,-2.9) {$\cdots$};
		
		\node [scale=2] at (-31.9,-6.3) {$U_{4_8}$};
		
		\node [scale=2] at (-32,-5.2) {$a$};
		
		\node [draw, shape=circle,fill=black] (v1) at  (-25,-4.5) {};
		\node [draw, shape=circle,fill=black] (v2) at  (-23.5,-6) {};
		\node [draw, shape=circle,fill=black] (v3) at  (-26.5,-6) {};	
		\node [draw, shape=circle,fill=black] (v4) at  (-25,-7.5) {};

		\node [draw, shape=circle,fill=black] (v5) at  (-26,-3.5) {};
		\node [draw, shape=circle,fill=black] (v6) at  (-24,-3.5) {};
		
		\node [draw, shape=circle,fill=black] (v11) at  (-27,-2.7) {};
		\node [draw, shape=circle,fill=black] (v12) at  (-27,-4.3) {};
		
		\node [draw, shape=circle,fill=black] (v13) at  (-23,-2.7) {};
		\node [draw, shape=circle,fill=black] (v14) at  (-23,-4.3) {};
		
		\node [draw, shape=circle,fill=black] (v7) at  (-27.8,-5.2) {};
		\node [draw, shape=circle,fill=black] (v8) at  (-27.8,-6.8) {};
		
		\node [draw, shape=circle,fill=black] (v9) at  (-22.7,-5.2) {};
		\node [draw, shape=circle,fill=black] (v10) at  (-22.7,-6.8) {};

		\draw (v1)--(v2);
		\draw (v1)--(v3);
		\draw(v4)--(v3);
		\draw (v2)--(v4);
		\draw (v1)--(v5);
		\draw (v1)--(v6);
		
		\draw (v5)--(v11);
		\draw (v5)--(v12);
		
		\draw (v6)--(v13);
		\draw (v6)--(v14);
		
		\draw (v3)--(v7);
		\draw (v3)--(v8);
		
		\draw (v2)--(v9);
		\draw (v2)--(v10);
		
		\node [scale=2] at (-22.3,-5.8) {$\vdots$};
		
		\node [scale=2] at (-27.8,-5.8) {$\vdots$};
		
		\node [scale=2] at (-25,-3.5) {$\cdots$};
		
		\node [scale=2] at (-27,-3.3) {$\vdots$};
		
		\node [scale=2] at (-23,-3.3) {$\vdots$};
		
		\node [scale=2] at (-24.9,-6.3) {$U_{4_{9}}$};
		
		\node [scale=2] at (-25,-5.2) {$a$};
		\node [draw, shape=circle,fill=black] (v1) at  (-18,-4.5) {};
		\node [draw, shape=circle,fill=black] (v2) at  (-16.5,-6) {};
		\node [draw, shape=circle,fill=black] (v3) at  (-19.5,-6) {};	
		\node [draw, shape=circle,fill=black] (v4) at  (-18,-7.5) {};

		\node [draw, shape=circle,fill=black] (v5) at  (-19,-3.5) {};
		\node [draw, shape=circle,fill=black] (v6) at  (-17,-3.5) {};
		
		\node [draw, shape=circle,fill=black] (v11) at  (-20,-2.7) {};
		\node [draw, shape=circle,fill=black] (v12) at  (-20,-4.3) {};
		
		\node [draw, shape=circle,fill=black] (v13) at  (-16,-2.7) {};
		\node [draw, shape=circle,fill=black] (v14) at  (-16,-4.3) {};

		\node [draw, shape=circle,fill=black] (v15) at  (-18.6,-2.6) {};
		\node [draw, shape=circle,fill=black] (v16) at  (-17.4,-2.6) {};
		
		\node [draw, shape=circle,fill=black] (v7) at  (-20.8,-5.2) {};
		\node [draw, shape=circle,fill=black] (v8) at  (-20.8,-6.8) {};
		
		\node [draw, shape=circle,fill=black] (v9) at  (-15.7,-5.2) {};
		\node [draw, shape=circle,fill=black] (v10) at  (-15.7,-6.8) {};
		
		\draw (v1)--(v2);
		\draw (v1)--(v3);
		\draw(v4)--(v3);
		\draw (v2)--(v4);
		\draw (v1)--(v5);
		\draw (v1)--(v6);
		
		\draw (v5)--(v11);
		\draw (v5)--(v12);
		
		\draw (v6)--(v13);
		\draw (v6)--(v14);
		
		\draw (v1)--(v15);
		\draw (v1)--(v16);
		
		\draw (v3)--(v7);
		\draw (v3)--(v8);
		
		\draw (v2)--(v9);
		\draw (v2)--(v10);
		
		\node [scale=2] at (-15.3,-5.8) {$\vdots$};

		\node [scale=2] at (-20.8,-5.8) {$\vdots$};
		
		\node [scale=2.7] at (-18,-3.7) {$\cdots$};
		
		\node [scale=2] at (-20,-3.3) {$\vdots$};
		
		\node [scale=2] at (-16,-3.3) {$\vdots$};
		
		\node [scale=2] at (-17.9,-2.9) {$\cdots$};
		
		\node [scale=2] at (-17.9,-6.3) {$U_{4_{10}}$};
		
		\node [scale=2] at (-18,-5.2) {$a$};
		
	\end{tikzpicture}

	\caption{Unicyclic graphs with cycle $C_4$ and radius two.}
	\label{u4-rad2}
\end{figure}

In this subsection, we discuss the class of unicyclic graphs of order at least five with radius two, whose cycle is $C_4$. The unicyclic graphs satisfying these conditions are illustrated in Figure \ref{u4-rad2}. Here $C=(a,b,c,d)$ is the cycle and $a$ is a central vertex. In $U_{4_1}$--$U_{4_4}$ the vertex $a$ has no strong children and the sets of support vertices are $\{a\}$, $\{a,b\}$, $\{b,d\}$ and $\{a,b,d\}$, respectively. In $U_{4_5}$--$U_{4_{10}}$ the vertex $a$ has at least one strong child and the cycle vertices adjacent to leaves are: none, $a$, $d$, $a$ and $d$, $b$ and $d$, and $a$, $b$ and $d$, respectively.

In the same way as in Subsection \ref{sub-sec-u-5}, we obtain the following observation. The proof is the same as that of Observation~\ref{obs-u5}, except that now the pendant tree at the vertex $c$ opposite to $a$ is trivial.

\begin{observation}\label{obs-u4-1}
	A graph $G$ of order at least five is a unicyclic graph containing the cycle $C_4$  and satisfying $\rad(G)\le 2$ if and only if   $G \in U_{4_i}$ for some $1\le i \le 10$.
\end{observation}
In what follows, we determine exact values of the global coalition number for some classes of graphs introduced above, and establish bounds for others. To aid our discussion, we assume that $G$ is a unicyclic graph of order $n$ containing the cycle $C_4 =(a,b,c,d)$ and having radius two.

\begin{proposition}\label{prop-u4-1}
	Let $G \in U_{4_1}$, and let $a$ be the support vertex of $G$.  Then
	\[
	GC(G) = \left\{
	\begin{array}{ll}
		\ 5 & \mbox{if\ $n_a =1$}, \\
		n-1 & \mbox{otherwise}.
	\end{array}
	\right.
	\]
\end{proposition}
\begin{proof}
If $n_a =1$, then it is clear that the singleton partition of $V(G)$ is a $gc$-partition of $G$ (the pairs $\{a,b\}$, $\{a,d\}$ and $\{a_1,c\}$ are global dominating sets), so $GC(G)=5$. Now assume $n_a \neq 1$. One can observe that the singleton partition of $V(G)$ is not a $gc$-partition of $G$ (no leaf of $a$ has an antitwin), implying that $GC(G) \leq n-1$. On the other hand, consider the vertex partition $\pi$ of $G$, where $\lbrace a,c\rbrace \in \pi$, and all other partition classes are singleton sets. Then $\pi$ is a $gc$-partition of $G$, where each singleton set forms a global coalition with $\lbrace a,c\rbrace$. Thus, $GC(G) \geq n-1$, and so $GC(G)=n-1$.
\end{proof}
\begin{theorem}\label{theo-u4-1}
	Let $G \in U_{4_2}$, and let $a,b$ be the support vertices of $G$.  Then
	\[
	GC(G) = \left\{
	\begin{array}{ll}
		\ 5 & \mbox{if\ $n_a =n_b  =1$}, \\
		\max \lbrace n_a +3 ,n_b +3 \rbrace & \mbox{otherwise}.
	\end{array}
	\right.
	\]
\end{theorem}
\begin{proof}
Assume, by symmetry, that $n_a \geq n_b$. If $n_a = n_b =1$, then it is easy to verify that the singleton partition of $V(G)$ is not a $gc$-partition of $G$, and that the partition $\lbrace L_a \cup L_b,\lbrace a\rbrace,\lbrace b\rbrace,\lbrace c\rbrace,\lbrace d\rbrace \rbrace$ is a  $gc$-partition of $G$. Thus, in this case, $GC(G)=5$. Hence, we may assume that $n_a \geq 2$. The partition  $\lbrace \lbrace a,c\rbrace \cup L_b, \lbrace b\rbrace ,\lbrace d\rbrace ,\lbrace a_i\rbrace : 1 \leq i \leq n_a \rbrace$ is a $gc$-partition of $G$ of size $n_a +3$, where the set $\lbrace a,c\rbrace \cup L_b$ forms a global coalition with all other members. Thus, $GC(G) \geq \max \lbrace n_a +3, n_b +3 \rbrace$. Now we show that $GC(G) \leq \max \lbrace n_a +3, n_b +3\rbrace$. Let $\pi$ be a $GC(G)$-partition. Since each subset $S \subseteq V(G)$ containing $a$ and $b$ is a global dominating set, it follows that the vertices $a$ and $b$ are in different partition classes. Let $V_a$ and $V_b$ be the members of $\pi$ containing $a$ and $b$, respectively. If there exists a pair $(A,B)$ of $gc$-partners such that $\lbrace A,B\rbrace \cap \lbrace V_a ,V_b \rbrace =\emptyset$, then we must have $(L_a \cup L_b) \subset (A \cup B)$, and $\lbrace c,d\rbrace \cap (A \cup B) \neq \emptyset$, which implies that $\lvert \pi \rvert \leq 5$, as desired. Hence, we may assume that each set in $\pi \setminus \lbrace V_a, V_b\rbrace$ has a $gc$-partner in $\lbrace V_a, V_b\rbrace$.

If there exists a partition class $C \in \pi$ forming a global coalition with both $V_a$ and $V_b$, then we must have $(L_a \cup L_b) \subseteq (C \cup V_a \cup V_b)$, implying that $\lvert \pi \rvert \leq 5$, as desired.

Next, we consider the case in which both the sets $V_a$ and $V_b$ have   $gc$-partners in $\pi \setminus \lbrace V_a, V_b\rbrace$. By the previous paragraph, these partners are distinct. Let $A$ and $B$ be  $gc$-partners of $V_a$ and $V_b$, respectively. It follows that $L_a \subset (B \cup V_b)$, $L_b \subset (A \cup V_a)$,  $(A \cup V_a) \cap \lbrace c,d\rbrace \neq \emptyset$ and $(B \cup V_b) \cap \lbrace c,d\rbrace \neq \emptyset$. Since $A\cup V_a$ and $B\cup V_b$ are disjoint, one of them contains $c$, and the other contains $d$, so these four classes cover $V(G)$. Thus, $\lvert \pi \rvert \leq 4$, as desired.  Hence, for the rest of the proof, we may assume that one of the following cases occurs.

\textbf{Case 1}. For any pair $(A,B)$ of $gc$-partners, $V_a \in \lbrace A,B\rbrace$.

Let $C$ be a $gc$-partner of $V_a$ with $C\ne V_b$. It follows that $L_b \subset (C \cup V_a)$, and $(C \cup V_a) \cap \lbrace c,d\rbrace \neq \emptyset$. Therefore, $\lvert \pi \rvert \leq n_a +4$. Now we show that $\lvert \pi \rvert \neq n_a +4$. Suppose to the contrary that $|\pi|=n_a+4$. In this case, either $C \cup V_a =L_b \cup \lbrace a,c\rbrace$ or $C \cup V_a =L_b \cup \lbrace a,d\rbrace$. Further, each member in $\pi \setminus \lbrace V_a, C\rbrace$ must be a singleton set. Assume first that $C \cup V_a =L_b \cup \lbrace a,c\rbrace$.  Since the singleton $\{d\}$ must form a global coalition with $V_a$ and $b\notin V_a\cup\{d\}$, we must have $L_b\subseteq V_a$, and hence $C=\{c\}$. But then $V_a\cup C=L_b\cup\{a,c\}$ is not a global dominating set, because in $\overline{G}$ the vertex $b$ is adjacent only to $d$ and to the vertices of $L_a$, none of which belongs to $V_a\cup C$. This contradicts the choice of $C$. A similar argument shows that if $C \cup V_a =L_b \cup \lbrace a,d\rbrace$, then we derive a contradiction. Hence, $\lvert \pi \rvert \leq n_a +3$.

\textbf{Case 2}. For any pair $(A,B)$ of $gc$-partners, $V_b \in \lbrace A,B\rbrace$.

Similar to the previous case, we can prove that $\lvert \pi \rvert \leq n_b +3$. Thus, $GC(G) \leq \max \lbrace n_a +3, n_b +3\rbrace$, and so $GC(G) = \max \lbrace n_a +3, n_b +3\rbrace$. This completes the proof.
\end{proof}
\begin{proposition}\label{prop-u4-2}
	If $G \in U_{4_3}$, then $GC(G)=n-1$.
\end{proposition}
\begin{proof}
Clearly, the singleton partition of $V(G)$ is not a $gc$-partition of $G$ (the vertex $a$ has no antitwin), which implies that $GC(G) \leq n-1$. On the other hand, consider the vertex partition $\pi$ of $G$, where one of the partition classes consists of the support vertices of $G$,  and the remaining classes are singleton sets. Then $\pi$ is a $gc$-partition of $G$, where each singleton set forms a global coalition with the set containing the support vertices of $G$. Thus, $GC(G) \geq n-1$, and so $GC(G)=n-1$.
\end{proof}
\begin{proposition}\label{prop-u4-3}
	If $G \in \bigcup_{i=4} ^{10} U_{4_i}$, then $GC(G) \geq n-\lvert N(a) \rvert +1$, where $a$ is the central vertex of $G$.
\end{proposition}
\begin{proof}
Consider the vertex partition $\pi$ of $G$, where $N(a) \in \pi$ and all other partition classes are singleton sets. Then $\pi$ is a $gc$-partition of $G$ of size $n- \lvert N(a) \rvert +1$, where the set $N(a)$ forms a global coalition with all other partition classes (Lemma~\ref{lem-uni-lb}). Hence,  $GC(G) \geq n- \lvert N(a) \rvert +1$.
\end{proof}
\subsection{Unicyclic graphs with cycle $C_3$}\label{sub-sec-u-3}
\begin{figure}[h]
	\centering
	\begin{tikzpicture}[scale=.4, transform shape]

		\node [draw, shape=circle,fill=black] (v1) at  (-45,1.5) {};
		\node [draw, shape=circle,fill=black] (v2) at  (-43.5,0) {};
		\node [draw, shape=circle,fill=black] (v3) at  (-46.5,0) {};

		\node [draw, shape=circle,fill=black] (v4) at  (-46,2.5) {};
		\node [draw, shape=circle,fill=black] (v5) at  (-44,2.5) {};

		\draw (v1)--(v2);
		\draw (v1)--(v3);
		\draw(v2)--(v3);
		
		\draw (v1)--(v4);
		\draw (v1)--(v5);
		
		\node [scale=2] at (-45,2.5) {$\cdots$};
		
		\node [scale=2] at (-44.3,1.4) {$a$};
		\node [scale=2] at (-43,.5) {$b$};

		\node [scale=2] at (-47,.5) {$c$};
		\node [scale=2] at (-44.8,.55) {$U_{3_1}$};
		
		
		\node [draw, shape=circle,fill=black] (v1) at  (-39,1.5) {};
		\node [draw, shape=circle,fill=black] (v2) at  (-37.5,0) {};
		\node [draw, shape=circle,fill=black] (v3) at  (-40.5,0) {};

		\node [draw, shape=circle,fill=black] (v6) at  (-36,1) {};
		\node [draw, shape=circle,fill=black] (v7) at  (-36,-1) {};
		
		\node [draw, shape=circle,fill=black] (v8) at  (-41.5,1) {};
		\node [draw, shape=circle,fill=black] (v9) at  (-41.5,-1) {};

		\draw (v1)--(v2);
		\draw (v1)--(v3);
		\draw(v2)--(v3);
		
		\draw (v2)--(v6);
		\draw (v2)--(v7);
		
		\draw (v3)--(v8);
		\draw (v3)--(v9);
		
		\node [scale=2] at (-41.5,.3) {$\vdots$};
		
		\node [scale=2] at (-36,.3) {$\vdots$};

		\node [scale=2] at (-38.8,.55) {$U_{3_2}$};
		

		\node [draw, shape=circle,fill=black] (v1) at  (-32,1.5) {};
		\node [draw, shape=circle,fill=black] (v2) at  (-30.5,0) {};
		\node [draw, shape=circle,fill=black] (v3) at  (-33.5,0) {};	
		
		\node [draw, shape=circle,fill=black] (v4) at  (-33,2.5) {};
		\node [draw, shape=circle,fill=black] (v5) at  (-31,2.5) {};
		
		\node [draw, shape=circle,fill=black] (v6) at  (-29,1) {};
		\node [draw, shape=circle,fill=black] (v7) at  (-29,-1) {};
		
		\node [draw, shape=circle,fill=black] (v8) at  (-34.5,1) {};
		\node [draw, shape=circle,fill=black] (v9) at  (-34.5,-1) {};

		\draw (v1)--(v2);
		\draw (v1)--(v3);
		\draw(v2)--(v3);
		
		\draw (v1)--(v4);
		\draw (v1)--(v5);
		
		\draw (v2)--(v6);
		\draw (v2)--(v7);
		
		\draw (v3)--(v8);
		\draw (v3)--(v9);
		
		\node [scale=2] at (-32,2.5) {$\cdots$};
		
		\node [scale=2] at (-34.5,.3) {$\vdots$};
		
		\node [scale=2] at (-29,.3) {$\vdots$};

		\node [scale=2] at (-31.8,.55) {$U_{3_3}$};

		
		\node [draw, shape=circle,fill=black] (v1) at  (-25,1.5) {};
		\node [draw, shape=circle,fill=black] (v2) at  (-23.5,0) {};
		\node [draw, shape=circle,fill=black] (v3) at  (-26.5,0) {};

		\node [draw, shape=circle,fill=black] (v4) at  (-26.2,2.5) {};
		\node [draw, shape=circle,fill=black] (v5) at  (-23.8,2.5) {};

		\node [draw, shape=circle,fill=black] (v12) at  (-27.4,3.2) {};
		\node [draw, shape=circle,fill=black] (v13) at  (-27.4,1.8) {};

		\node [draw, shape=circle,fill=black] (v14) at  (-22.6,3.2) {};
		\node [draw, shape=circle,fill=black] (v15) at  (-22.6,1.8) {};

		\draw (v1)--(v2);
		\draw (v1)--(v3);
		\draw(v2)--(v3);
		
		\draw (v1)--(v4);
		\draw (v1)--(v5);
		
		;
		
		\draw (v4)--(v12);
		\draw (v4)--(v13);
		
		\draw (v5)--(v14);
		\draw (v5)--(v15);

		\node [scale=2.7] at (-24.9,2.3) {$\cdots$};
		
		
		\node [scale=2] at (-27.3,2.75) {$\vdots$};
		\node [scale=2] at (-22.5,2.75) {$\vdots$};
		
		\node [scale=2] at (-24.8,.55) {$U_{3_4}$};
		
		\node [scale=2] at (-25.7,1.5) {$a$};

		\node [draw, shape=circle,fill=black] (v1) at  (-18,1.5) {};
		\node [draw, shape=circle,fill=black] (v2) at  (-16.5,0) {};
		\node [draw, shape=circle,fill=black] (v3) at  (-19.5,0) {};

		\node [draw, shape=circle,fill=black] (v4) at  (-19.2,2.5) {};
		\node [draw, shape=circle,fill=black] (v5) at  (-16.8,2.5) {};
		
		\node [draw, shape=circle,fill=black] (v10) at  (-18.8,3.5) {};
		\node [draw, shape=circle,fill=black] (v11) at  (-17.2,3.5) {};

		\node [draw, shape=circle,fill=black] (v12) at  (-20.4,3.2) {};
		\node [draw, shape=circle,fill=black] (v13) at  (-20.4,1.8) {};

		\node [draw, shape=circle,fill=black] (v14) at  (-15.6,3.2) {};
		\node [draw, shape=circle,fill=black] (v15) at  (-15.6,1.8) {};

		\draw (v1)--(v2);
		\draw (v1)--(v3);
		\draw(v2)--(v3);
		
		\draw (v1)--(v4);
		\draw (v1)--(v5);
		
		\draw (v1)--(v10);
		\draw (v1)--(v11);
		
		\draw (v4)--(v12);
		\draw (v4)--(v13);
		
		\draw (v5)--(v14);
		\draw (v5)--(v15);

		\node [scale=2.7] at (-17.9,2.3) {$\cdots$};
		
		\node [scale=2.] at (-17.9,3.5) {$\cdots$};
		
		\node [scale=2] at (-20.3,2.75) {$\vdots$};
		\node [scale=2] at (-15.5,2.75) {$\vdots$};
		
		\node [scale=2] at (-17.8,.55) {$U_{3_5}$};
		
		\node [scale=2] at (-18.7,1.5) {$a$};
		
		\node [draw, shape=circle,fill=black] (v1) at  (-45,-4.5) {};
		\node [draw, shape=circle,fill=black] (v2) at  (-43.5,-6) {};
		\node [draw, shape=circle,fill=black] (v3) at  (-46.5,-6) {};

		\node [draw, shape=circle,fill=black] (v4) at  (-46.2,-3.5) {};
		\node [draw, shape=circle,fill=black] (v5) at  (-43.8,-3.5) {};
		
		
		\node [draw, shape=circle,fill=black] (v8) at  (-47.5,-5) {};
		\node [draw, shape=circle,fill=black] (v9) at  (-47.5,-7) {};
		
		\node [draw, shape=circle,fill=black] (v12) at  (-47.4,-2.8) {};
		\node [draw, shape=circle,fill=black] (v13) at  (-47.4,-4.2) {};

		\node [draw, shape=circle,fill=black] (v14) at  (-42.6,-2.8) {};
		\node [draw, shape=circle,fill=black] (v15) at  (-42.6,-4.2) {};

		\draw (v1)--(v2);
		\draw (v1)--(v3);
		\draw(v2)--(v3);
		
		\draw (v1)--(v4);
		\draw (v1)--(v5);
		
		
		\draw (v3)--(v8);
		\draw (v3)--(v9);
		
		\draw (v4)--(v12);
		\draw (v4)--(v13);
		
		\draw (v5)--(v14);
		\draw (v5)--(v15);
		
		\node [scale=2.7] at (-44.9,-3.7) {$\cdots$};
		
		
		\node [scale=2] at (-47.5,-5.7) {$\vdots$};
		
		\node [scale=2] at (-47.3,-3.25) {$\vdots$};
		\node [scale=2] at (-42.5,-3.25) {$\vdots$};
		
		\node [scale=2] at (-44.8,-5.45) {$U_{3_6}$};
		
		\node [scale=2] at (-45.7,-4.5) {$a$};
		
		\node [draw, shape=circle,fill=black] (v1) at  (-38,-4.5) {};
		\node [draw, shape=circle,fill=black] (v2) at  (-36.5,-6) {};
		\node [draw, shape=circle,fill=black] (v3) at  (-39.5,-6) {};

		\node [draw, shape=circle,fill=black] (v4) at  (-39.2,-3.5) {};
		\node [draw, shape=circle,fill=black] (v5) at  (-36.8,-3.5) {};
		
		\node [draw, shape=circle,fill=black] (v10) at  (-38.8,-2.5) {};
		\node [draw, shape=circle,fill=black] (v11) at  (-37.2,-2.5) {};
		
		\node [draw, shape=circle,fill=black] (v8) at  (-40.5,-5) {};
		\node [draw, shape=circle,fill=black] (v9) at  (-40.5,-7) {};
		
		\node [draw, shape=circle,fill=black] (v12) at  (-40.4,-2.8) {};
		\node [draw, shape=circle,fill=black] (v13) at  (-40.4,-4.2) {};

		\node [draw, shape=circle,fill=black] (v14) at  (-35.6,-2.8) {};
		\node [draw, shape=circle,fill=black] (v15) at  (-35.6,-4.2) {};

		\draw (v1)--(v2);
		\draw (v1)--(v3);
		\draw(v2)--(v3);
		
		\draw (v1)--(v4);
		\draw (v1)--(v5);
		
		\draw (v1)--(v10);
		\draw (v1)--(v11);
		
		\draw (v3)--(v8);
		\draw (v3)--(v9);
		
		\draw (v4)--(v12);
		\draw (v4)--(v13);
		
		\draw (v5)--(v14);
		\draw (v5)--(v15);
		
		\node [scale=2.7] at (-37.9,-3.7) {$\cdots$};
		
		\node [scale=2.] at (-37.9,-2.5) {$\cdots$};
		
		\node [scale=2] at (-40.5,-5.7) {$\vdots$};
		
		\node [scale=2] at (-40.3,-3.25) {$\vdots$};
		\node [scale=2] at (-35.5,-3.25) {$\vdots$};
		
		\node [scale=2] at (-37.8,-5.45) {$U_{3_7}$};
		
		\node [scale=2] at (-38.7,-4.5) {$a$};

		\node [draw, shape=circle,fill=black] (v1) at  (-31,-4.5) {};
		\node [draw, shape=circle,fill=black] (v2) at  (-29.5,-6) {};
		\node [draw, shape=circle,fill=black] (v3) at  (-32.5,-6) {};

		\node [draw, shape=circle,fill=black] (v4) at  (-32.2,-3.5) {};
		\node [draw, shape=circle,fill=black] (v5) at  (-29.8,-3.5) {};
		
		\node [draw, shape=circle,fill=black] (v6) at  (-28.5,-5) {};
		\node [draw, shape=circle,fill=black] (v7) at  (-28.5,-7) {};

		\node [draw, shape=circle,fill=black] (v8) at  (-33.5,-5) {};
		\node [draw, shape=circle,fill=black] (v9) at  (-33.5,-7) {};
		
		\node [draw, shape=circle,fill=black] (v12) at  (-33.4,-2.8) {};
		\node [draw, shape=circle,fill=black] (v13) at  (-33.4,-4.2) {};

		\node [draw, shape=circle,fill=black] (v14) at  (-28.6,-2.8) {};
		\node [draw, shape=circle,fill=black] (v15) at  (-28.6,-4.2) {};

		\draw (v1)--(v2);
		\draw (v1)--(v3);
		\draw(v2)--(v3);
		
		\draw (v1)--(v4);
		\draw (v1)--(v5);
		
		\draw (v2)--(v6);
		\draw (v2)--(v7);
		
		\draw (v3)--(v8);
		\draw (v3)--(v9);
		
		\draw (v4)--(v12);
		\draw (v4)--(v13);
		
		\draw (v5)--(v14);
		\draw (v5)--(v15);
		
		\node [scale=2.7] at (-30.9,-3.7) {$\cdots$};
		
		\node [scale=2.] at (-28.5,-5.7) {$\vdots$};
		
		\node [scale=2] at (-33.5,-5.7) {$\vdots$};
		
		\node [scale=2] at (-33.3,-3.25) {$\vdots$};
		\node [scale=2] at (-28.5,-3.25) {$\vdots$};
		
		\node [scale=2] at (-30.8,-5.45) {$U_{3_8}$};
		
		\node [scale=2] at (-31.7,-4.5) {$a$};
		
		\node [draw, shape=circle,fill=black] (v1) at  (-23,-4.5) {};
		\node [draw, shape=circle,fill=black] (v2) at  (-21.5,-6) {};
		\node [draw, shape=circle,fill=black] (v3) at  (-24.5,-6) {};

		\node [draw, shape=circle,fill=black] (v4) at  (-24.2,-3.5) {};
		\node [draw, shape=circle,fill=black] (v5) at  (-21.8,-3.5) {};
		
		\node [draw, shape=circle,fill=black] (v6) at  (-20.5,-5) {};
		\node [draw, shape=circle,fill=black] (v7) at  (-20.5,-7) {};
		
		\node [draw, shape=circle,fill=black] (v10) at  (-23.8,-2.5) {};
		\node [draw, shape=circle,fill=black] (v11) at  (-22.2,-2.5) {};
		
		\node [draw, shape=circle,fill=black] (v8) at  (-25.5,-5) {};
		\node [draw, shape=circle,fill=black] (v9) at  (-25.5,-7) {};
		
		\node [draw, shape=circle,fill=black] (v12) at  (-25.4,-2.8) {};
		\node [draw, shape=circle,fill=black] (v13) at  (-25.4,-4.2) {};

		\node [draw, shape=circle,fill=black] (v14) at  (-20.6,-2.8) {};
		\node [draw, shape=circle,fill=black] (v15) at  (-20.6,-4.2) {};

		\draw (v1)--(v2);
		\draw (v1)--(v3);
		\draw(v2)--(v3);
		
		\draw (v1)--(v4);
		\draw (v1)--(v5);
		
		\draw (v1)--(v10);
		\draw (v1)--(v11);
		
		\draw (v2)--(v6);
		\draw (v2)--(v7);
		
		\draw (v3)--(v8);
		\draw (v3)--(v9);
		
		\draw (v4)--(v12);
		\draw (v4)--(v13);
		
		\draw (v5)--(v14);
		\draw (v5)--(v15);
		
		\node [scale=2.7] at (-22.9,-3.7) {$\cdots$};
		
		\node [scale=2.] at (-29.5,-5.7) {$\vdots$};
		
		\node [scale=2.] at (-22.9,-2.5) {$\cdots$};
		
		\node [scale=2] at (-20.5,-5.7) {$\vdots$};
		
		\node [scale=2] at (-25.3,-3.25) {$\vdots$};
		\node [scale=2] at (-20.5,-3.25) {$\vdots$};
		
		\node [scale=2] at (-22.8,-5.45) {$U_{3_9}$};
		
		\node [scale=2] at (-23.7,-4.5) {$a$};

		\node [draw, shape=circle,fill=black] (v1) at  (-17,-4.5) {};
		\node [draw, shape=circle,fill=black] (v2) at  (-15.5,-6) {};
		\node [draw, shape=circle,fill=black] (v3) at  (-18.5,-6) {};	
		\node [draw, shape=circle,fill=black] (v17) at  (-17,-3.5) {};
		\node [draw, shape=circle,fill=black] (v18) at  (-17,-2.5) {};
		
		\node [draw, shape=circle,fill=black] (v19) at  (-18,-1.5) {};
		\node [draw, shape=circle,fill=black] (v20) at  (-16,-1.5) {};

		\draw (v1)--(v2);
		\draw (v1)--(v3);
		\draw(v2)--(v3);
		
		\draw (v1)--(v17);
		\draw (v17)--(v18);
		
		\draw (v18)--(v19);
		\draw (v18)--(v20);

		\node [scale=2] at (-17.,-1.5) {$\cdots$};
		
		
		\node [scale=2] at (-16.8,-5.45) {$U_{3_{10}}$};
		
		\node [scale=2] at (-17.7,-3.5) {$a$};
		
		\node [draw, shape=circle,fill=black] (v1) at  (-45,-11.5) {};
		\node [draw, shape=circle,fill=black] (v2) at  (-43.5,-13) {};
		\node [draw, shape=circle,fill=black] (v3) at  (-46.5,-13) {};	
		
		\node [draw, shape=circle,fill=black] (v17) at  (-45,-10.5) {};
		\node [draw, shape=circle,fill=black] (v18) at  (-45,-9.5) {};
		
		\node [draw, shape=circle,fill=black] (v19) at  (-46,-8.5) {};
		\node [draw, shape=circle,fill=black] (v20) at  (-44,-8.5) {};	
		
		\node [draw, shape=circle,fill=black] (v21) at  (-46,-9.5) {};
		\node [draw, shape=circle,fill=black] (v22) at  (-44,-9.5) {};

		\draw (v1)--(v2);
		\draw (v1)--(v3);
		\draw(v2)--(v3);
		
		draw (v1)--(v17);
		draw (v17)--(v18);
		
		\draw (v17)--(v21);
		\draw (v17)--(v22);
		
		\draw (v18)--(v19);
		\draw (v18)--(v20);
		
		\node [scale=2.3] at (-44.9,-9.8) {$\cdots$};
		\node [scale=2] at (-45,-8.5) {$\cdots$};
		
		
		\node [scale=2] at (-44.8,-12.45) {$U_{3_{11}}$};
		
		\node [scale=2] at (-45.6,-10.5) {$a$};
		
		\node [draw, shape=circle,fill=black] (v1) at  (-40,-11.5) {};
		\node [draw, shape=circle,fill=black] (v2) at  (-38.5,-13) {};
		\node [draw, shape=circle,fill=black] (v3) at  (-41.5,-13) {};	
		
		\node [draw, shape=circle,fill=black] (v17) at  (-40,-10.5) {};
		\node [draw, shape=circle,fill=black] (v18) at  (-40,-9.5) {};
		
		\node [draw, shape=circle,fill=black] (v19) at  (-41,-8.5) {};
		\node [draw, shape=circle,fill=black] (v20) at  (-39,-8.5) {};	
		
		\node [draw, shape=circle,fill=black] (v21) at  (-41,-9.5) {};
		\node [draw, shape=circle,fill=black] (v22) at  (-39,-9.5) {};
		
		\node [draw, shape=circle,fill=black] (v23) at  (-42,-10.3) {};
		\node [draw, shape=circle,fill=black] (v24) at  (-42,-8.7) {};
		
		\node [draw, shape=circle,fill=black] (v25) at  (-38,-10.3) {};
		\node [draw, shape=circle,fill=black] (v26) at  (-38,-8.7) {};

		\draw (v1)--(v2);
		\draw (v1)--(v3);
		\draw(v2)--(v3);
		
		draw (v1)--(v17);
		draw (v17)--(v18);
		
		\draw (v17)--(v21);
		\draw (v17)--(v22);
		
		\draw (v18)--(v19);
		\draw (v18)--(v20);
		
		\draw (v21)--(v23);
		\draw (v21)--(v24);
		
		\draw (v22)--(v25);
		\draw (v22)--(v26);
		
		\node [scale=2.3] at (-39.9,-9.8) {$\cdots$};
		\node [scale=2] at (-40,-8.5) {$\cdots$};
		
		\node [scale=2] at (-38,-9.2) {$\vdots$};
		\node [scale=2] at (-42,-9.2) {$\vdots$};
		
		\node [scale=2] at (-39.8,-12.45) {$U_{3_{12}}$};
		
		\node [scale=2] at (-40.6,-10.5) {$a$};
		
		\node [draw, shape=circle,fill=black] (v1) at  (-35,-11.5) {};
		\node [draw, shape=circle,fill=black] (v2) at  (-33.5,-13) {};
		\node [draw, shape=circle,fill=black] (v3) at  (-36.5,-13) {};	
		
		\node [draw, shape=circle,fill=black] (v17) at  (-35,-10.5) {};
		\node [draw, shape=circle,fill=black] (v18) at  (-35,-9.5) {};
		
		\node [draw, shape=circle,fill=black] (v19) at  (-36,-8.5) {};
		\node [draw, shape=circle,fill=black] (v20) at  (-34,-8.5) {};	
		
		\node [draw, shape=circle,fill=black] (v21) at  (-36,-9.5) {};
		\node [draw, shape=circle,fill=black] (v22) at  (-34,-9.5) {};
		
		\node [draw, shape=circle,fill=black] (v23) at  (-37,-10.3) {};
		\node [draw, shape=circle,fill=black] (v24) at  (-37,-8.7) {};
		
		\node [draw, shape=circle,fill=black] (v25) at  (-33,-10.3) {};
		\node [draw, shape=circle,fill=black] (v26) at  (-33,-8.7) {};
		
		\node [draw, shape=circle,fill=black] (v27) at  (-36,-11.5) {};
		\node [draw, shape=circle,fill=black] (v28) at  (-34,-11.5) {};

		\draw (v1)--(v2);
		\draw (v1)--(v3);
		\draw(v2)--(v3);
		
		draw (v1)--(v17);
		draw (v17)--(v18);
		
		\draw (v17)--(v21);
		\draw (v17)--(v22);
		
		\draw (v18)--(v19);
		\draw (v18)--(v20);
		
		\draw (v21)--(v23);
		\draw (v21)--(v24);
		
		\draw (v22)--(v25);
		\draw (v22)--(v26);
		
		\draw (v17)--(v27);
		\draw (v17)--(v28);
		
		\node [scale=2.3] at (-34.9,-9.8) {$\cdots$};
		\node [scale=2] at (-35,-8.5) {$\cdots$};
		
		\node [scale=2] at (-35,-11.2) {$\cdots$};
		
		\node [scale=2] at (-33,-9.2) {$\vdots$};
		\node [scale=2] at (-37,-9.2) {$\vdots$};
		
		\node [scale=2] at (-34.8,-12.45) {$U_{3_{13}}$};
		
		\node [scale=2] at (-35.6,-10.5) {$a$};

		\node [draw, shape=circle,fill=black] (v1) at  (-30,-11.5) {};
		\node [draw, shape=circle,fill=black] (v2) at  (-28.5,-13) {};
		\node [draw, shape=circle,fill=black] (v3) at  (-31.5,-13) {};	
		\node [draw, shape=circle,fill=black] (v17) at  (-30,-10.5) {};
		\node [draw, shape=circle,fill=black] (v18) at  (-30,-9.5) {};
		
		\node [draw, shape=circle,fill=black] (v19) at  (-31,-8.5) {};
		\node [draw, shape=circle,fill=black] (v20) at  (-29,-8.5) {};
		
		\node [draw, shape=circle,fill=black] (v29) at  (-31.5,-12.2) {};
		\node [draw, shape=circle,fill=black] (v30) at  (-28.5,-12.2) {};

		\draw (v1)--(v29);
		\draw (v1)--(v30);
		
		\node [scale=2] at (-30.,-12.) {$\cdots$};

		\draw (v1)--(v2);
		\draw (v1)--(v3);
		\draw(v2)--(v3);
		
		\draw (v1)--(v17);
		\draw (v17)--(v18);
		
		\draw (v18)--(v19);
		\draw (v18)--(v20);

		\node [scale=2] at (-30.,-8.5) {$\cdots$};
		
		
		\node [scale=2] at (-29.8,-12.55) {$U_{3_{14}}$};

\node [scale=2] at (-30.6,-10.5) {$a$};		
		
		\node [draw, shape=circle,fill=black] (v1) at  (-25,-11.5) {};
		\node [draw, shape=circle,fill=black] (v2) at  (-23.5,-13) {};
		\node [draw, shape=circle,fill=black] (v3) at  (-26.5,-13) {};	
		
		\node [draw, shape=circle,fill=black] (v17) at  (-25,-10.5) {};
		\node [draw, shape=circle,fill=black] (v18) at  (-25,-9.5) {};
		
		\node [draw, shape=circle,fill=black] (v19) at  (-26,-8.5) {};
		\node [draw, shape=circle,fill=black] (v20) at  (-24,-8.5) {};	
		
		\node [draw, shape=circle,fill=black] (v21) at  (-26,-9.5) {};
		\node [draw, shape=circle,fill=black] (v22) at  (-24,-9.5) {};
		
		\node [draw, shape=circle,fill=black] (v29) at  (-26.5,-12.2) {};
		\node [draw, shape=circle,fill=black] (v30) at  (-23.5,-12.2) {};

		\draw (v1)--(v29);
		\draw (v1)--(v30);
		
		\node [scale=2] at (-25.,-12.) {$\cdots$};	

		\draw (v1)--(v2);
		\draw (v1)--(v3);
		\draw(v2)--(v3);
		
		draw (v1)--(v17);
		draw (v17)--(v18);
		
		\draw (v17)--(v21);
		\draw (v17)--(v22);
		
		\draw (v18)--(v19);
		\draw (v18)--(v20);
		
		\node [scale=2.3] at (-24.9,-9.8) {$\cdots$};
		\node [scale=2] at (-25,-8.5) {$\cdots$};

		\node [scale=2] at (-24.8,-12.55) {$U_{3_{15}}$};
		
		\node [scale=2] at (-25.6,-10.5) {$a$};
		
		\node [draw, shape=circle,fill=black] (v1) at  (-20,-11.5) {};
		\node [draw, shape=circle,fill=black] (v2) at  (-18.5,-13) {};
		\node [draw, shape=circle,fill=black] (v3) at  (-21.5,-13) {};	
		
		\node [draw, shape=circle,fill=black] (v17) at  (-20,-10.5) {};
		\node [draw, shape=circle,fill=black] (v18) at  (-20,-9.5) {};
		
		\node [draw, shape=circle,fill=black] (v19) at  (-21,-8.5) {};
		\node [draw, shape=circle,fill=black] (v20) at  (-19,-8.5) {};	
		
		\node [draw, shape=circle,fill=black] (v21) at  (-21,-9.5) {};
		\node [draw, shape=circle,fill=black] (v22) at  (-19,-9.5) {};
		
		\node [draw, shape=circle,fill=black] (v23) at  (-22,-10.3) {};
		\node [draw, shape=circle,fill=black] (v24) at  (-22,-8.7) {};
		
		\node [draw, shape=circle,fill=black] (v25) at  (-18,-10.3) {};
		\node [draw, shape=circle,fill=black] (v26) at  (-18,-8.7) {};

		\node [draw, shape=circle,fill=black] (v29) at  (-21.5,-12.2) {};
		\node [draw, shape=circle,fill=black] (v30) at  (-18.5,-12.2) {};

		\draw (v1)--(v29);
		\draw (v1)--(v30);
		
		\node [scale=2] at (-20.,-12.) {$\cdots$};	

		\draw (v1)--(v2);
		\draw (v1)--(v3);
		\draw(v2)--(v3);
		
		draw (v1)--(v17);
		draw (v17)--(v18);
		
		\draw (v17)--(v21);
		\draw (v17)--(v22);
		
		\draw (v18)--(v19);
		\draw (v18)--(v20);
		
		\draw (v21)--(v23);
		\draw (v21)--(v24);
		
		\draw (v22)--(v25);
		\draw (v22)--(v26);
		
		\node [scale=2.3] at (-19.9,-9.8) {$\cdots$};
		\node [scale=2] at (-20,-8.5) {$\cdots$};
		
		\node [scale=2] at (-18,-9.2) {$\vdots$};
		\node [scale=2] at (-22,-9.2) {$\vdots$};
		
		\node [scale=2] at (-19.8,-12.55) {$U_{3_{16}}$};
		
		\node [scale=2] at (-20.6,-10.5) {$a$};
		
		\node [draw, shape=circle,fill=black] (v1) at  (-15,-11.5) {};
		\node [draw, shape=circle,fill=black] (v2) at  (-13.5,-13) {};
		\node [draw, shape=circle,fill=black] (v3) at  (-16.5,-13) {};	
		
		\node [draw, shape=circle,fill=black] (v17) at  (-15,-10.5) {};
		\node [draw, shape=circle,fill=black] (v18) at  (-15,-9.5) {};
		
		\node [draw, shape=circle,fill=black] (v19) at  (-16,-8.5) {};
		\node [draw, shape=circle,fill=black] (v20) at  (-14,-8.5) {};	
		
		\node [draw, shape=circle,fill=black] (v21) at  (-16,-9.5) {};
		\node [draw, shape=circle,fill=black] (v22) at  (-14,-9.5) {};
		
		\node [draw, shape=circle,fill=black] (v23) at  (-17,-10.3) {};
		\node [draw, shape=circle,fill=black] (v24) at  (-17,-8.7) {};
		
		\node [draw, shape=circle,fill=black] (v25) at  (-13,-10.3) {};
		\node [draw, shape=circle,fill=black] (v26) at  (-13,-8.7) {};
		
		\node [draw, shape=circle,fill=black] (v27) at  (-16,-11.5) {};
		\node [draw, shape=circle,fill=black] (v28) at  (-14,-11.5) {};

		\node [draw, shape=circle,fill=black] (v29) at  (-16.5,-12.2) {};
		\node [draw, shape=circle,fill=black] (v30) at  (-13.5,-12.2) {};

		\draw (v1)--(v29);
		\draw (v1)--(v30);
		
		\node [scale=2] at (-15.,-12.) {$\cdots$};	

		\draw (v1)--(v2);
		\draw (v1)--(v3);
		\draw(v2)--(v3);
		
		draw (v1)--(v17);
		draw (v17)--(v18);
		
		\draw (v17)--(v21);
		\draw (v17)--(v22);
		
		\draw (v18)--(v19);
		\draw (v18)--(v20);
		
		\draw (v21)--(v23);
		\draw (v21)--(v24);
		
		\draw (v22)--(v25);
		\draw (v22)--(v26);
		
		\draw (v17)--(v27);
		\draw (v17)--(v28);
		
		\node [scale=2.3] at (-14.9,-9.8) {$\cdots$};
		\node [scale=2] at (-15,-8.5) {$\cdots$};
		
		\node [scale=2] at (-15,-11.2) {$\cdots$};
		
		\node [scale=2] at (-13,-9.2) {$\vdots$};
		\node [scale=2] at (-17,-9.2) {$\vdots$};
		
		\node [scale=2] at (-14.8,-12.55) {$U_{3_{17}}$};
		
		\node [scale=2] at (-15.6,-10.5) {$a$};

	\end{tikzpicture}

	\caption{Unicyclic graphs with cycle $C_3$ and radius at most two.}
	\label{u3-rad2}
\end{figure}
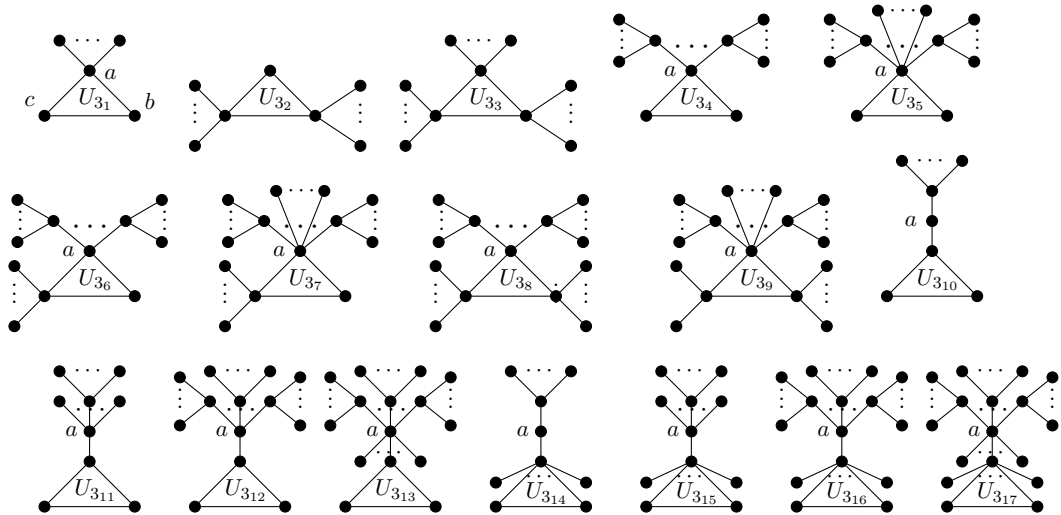

In this subsection, we discuss the class of unicyclic graphs of order at least four with radius at most two, whose cycle is $C_3$. The unicyclic graphs satisfying these conditions are illustrated in Figure \ref{u3-rad2}. Let $C=(a,b,c)$ be the triangle. In the families $U_{3_1}$--$U_{3_9}$ a central vertex $a$ lies on the triangle: in $U_{3_1}$--$U_{3_3}$ the vertex $a$ has no strong children and the sets of support vertices are $\{a\}$, $\{b,c\}$ and $\{a,b,c\}$, respectively, while in $U_{3_4}$--$U_{3_9}$ the vertex $a$ has at least one strong child and the triangle vertices adjacent to leaves are: none, $a$, $c$, $a$ and $c$, $b$ and $c$, and $a$, $b$ and $c$, respectively. In the families $U_{3_{10}}$--$U_{3_{17}}$, a central vertex $a$ lies outside the triangle and is adjacent to the triangle vertex $t$. The vertex $a$ has $k\ge 1$ strong children, no other vertex of the triangle has further neighbours, and leaves may be attached to $a$ and to $t$. Here $k=1$ in $U_{3_{10}}$, $U_{3_{11}}$, $U_{3_{14}}$, $U_{3_{15}}$ and $k\ge 2$ in $U_{3_{12}}$, $U_{3_{13}}$, $U_{3_{16}}$, $U_{3_{17}}$. The vertex $a$ is a support vertex in $U_{3_{11}}$, $U_{3_{13}}$, $U_{3_{15}}$, $U_{3_{17}}$, and $t$ is a support vertex in $U_{3_{14}}$--$U_{3_{17}}$.

As in Subsections \ref{sub-sec-u-5} and \ref{sub-sec-u-4}, we obtain the following observation.

\begin{observation}\label{obs-u3-1}
	A graph $G$ of order at least four is a unicyclic graph containing the cycle $C_3$  and satisfying $\rad(G)\le 2$ if and only if    $G \in U_{3_i}$ for some $1\le i \le 17$.
\end{observation}
\begin{proof}
Let $G$ be unicyclic with triangle $C=(a,b,c)$, of order at least four and radius at most two, and let $v$ be a central vertex. If $v\in V(C)$, say $v=a$, then, as in the proof of Observation~\ref{obs-u5}, the pendant trees at $b$ and $c$ consist of leaves and the pendant tree at $a$ has depth at most two. This gives the families $U_{3_1}$--$U_{3_9}$ (a graph with leaves at $a$ and $b$ only is isomorphic to one in $U_{3_2}$). If $v\notin V(C)$, let $t$ be the vertex of $C$ closest to $v$. Since $\ecc(v)\le 2$, the vertex $v$ is adjacent to $t$, the other two triangle vertices have no neighbours outside $C$, the pendant tree at $t$, other than the branch containing $v$, consists of leaves, and $v$ is adjacent to leaves and strong children only. If $v$ had no strong child, then $t$ would also be a central vertex and $G$ would belong to $U_{3_4}$ or $U_{3_5}$. Otherwise $G$ belongs to one of $U_{3_{10}}$--$U_{3_{17}}$ according to the number of strong children of $v$ and to whether $v$ and $t$ are support vertices. The converse is clear.
\end{proof}

In what follows, we determine exact values of the global coalition number for some classes of graphs introduced above, and obtain bounds for others. To aid our discussion, we consider $G$ to be a unicyclic graph of order $n$ with cycle $C_3 =(a,b,c)$ and radius at most two.
\begin{proposition}\label{prop-u3-1}
	If $G \in U_{3_1} \cup U_{3_2} \cup U_{3_{10}} \cup U_{3_{14}}$, then $GC(G)=n-1$.
\end{proposition}
\begin{proof}
It is easily seen that if $G \in U_{3_1} \cup U_{3_2} \cup U_{3_{10}} \cup U_{3_{14}}$, then the singleton partition of $V(G)$ is not a $gc$-partition of $G$ (in the first two families the vertex $b$, respectively $a$, has no antitwin, while in the last two families the central vertex $a$ has none). Now let $G \in U_{3_1}$, and
let $a$ be the support vertex of $G$. Consider the vertex partition $\pi$ of $G$, where $\lbrace b,c\rbrace \in \pi$,  and the other partition classes are singleton sets. Then $\pi$ is a $gc$-partition of $G$, where the set $\lbrace a\rbrace$ forms a global coalition with all other partition classes. Thus, $GC(G)=n-1$.

Assume next that $G \in U_{3_2}$, and let $\lbrace b,c\rbrace$ be the set of support vertices of $G$. Consider the vertex partition $\pi$ of $G$, where $\lbrace b,c\rbrace \in \pi$,  and the other partition classes are singleton sets. Then $\pi$ is a $gc$-partition of $G$, where each singleton set forms a global coalition with $\lbrace b,c\rbrace$. Thus, $GC(G)=n-1$.

Finally, assume $G \in U_{3_{10}} \cup U_{3_{14}}$. Let $a$ be the central vertex of $G$. Consider the vertex partition $\pi$ of $G$, where $N(a) \in \pi$,  and the other partition classes are singleton sets. Then $\pi$ is a $gc$-partition of $G$, where  each singleton set forms a global coalition with $N(a)$ (Lemma~\ref{lem-uni-lb}). Here $|N(a)|=2$. Thus, $GC(G)=n-1$.
\end{proof}
	\begin{theorem}\label{theo-u3-1}
		Let $G \in U_{3_3}$, and let $a,b,c$ be the support vertices of $G$ such that $n_a \geq n_b \geq n_c$. Then $GC(G)=n_a +n_b +2$.
	\end{theorem}
	\begin{proof}
	Consider the vertex partition $\pi$ of $G$, where $N(c) \in \pi$, and all other partition classes are singleton sets. Then $\pi$ is a $gc$-partition of $G$ of size $n_a + n_b + 2$, where the set $N(c)$ forms a global coalition with all other members. Thus, $GC(G) \geq n_a +n_b +2$. Now we show that $GC(G) \leq n_a + n_b +2$. Let $\pi$ be a $gc$-partition of $G$. Note that each subset $S \subseteq V(G)$ containing $a$, $b$ and $c$ is a global dominating set of $G$, so the vertices $a$, $b$ and $c$ cannot all be in the same partition class. Therefore, we have the following cases.

	\textbf{Case 1.} Each of the vertices $a$, $b$ and $c$ belongs to a distinct partition class.

	Let $V_a$, $V_b$ and $V_c$ be the members of $\pi$ containing $a$, $b$ and $c$, respectively. If $\pi =\lbrace V_a,V_b,V_c\rbrace$, then we are done. Otherwise, we claim that each member in $\pi \setminus \lbrace V_a,V_b,V_c\rbrace$ has a $gc$-partner in $\lbrace V_a,V_b,V_c\rbrace$. Suppose, to the contrary, that this is not the case. Let $A$ be a member of $\pi$ having no $gc$-partner in $\lbrace V_a,V_b,V_c\rbrace$, and let $B$ be a $gc$-partner of $A$. It follows that $A \cup B =L_a \cup L_b \cup L_c$, and hence $V_a=\{a\}$, $V_b=\{b\}$ and $V_c=\{c\}$. Since no two of $\{a\},\{b\},\{c\}$ form a global coalition, each of them has its $gc$-partner in $\{A,B\}$. The partner of $\{a\}$ contains $L_b\cup L_c$, that of $\{b\}$ contains $L_a\cup L_c$ and that of $\{c\}$ contains $L_a\cup L_b$. Two of the three singletons share a partner, which then contains $L_a\cup L_b\cup L_c$, contradicting $A, B\ne\emptyset$. Now let $C$ be an arbitrary member of $\pi $ such that $C \notin \lbrace V_a,V_b,V_c \rbrace$.
	If $C$ is a $gc$-partner of $V_a$, then we must have $(L_b \cup L_c ) \subset (C \cup V_a)$, implying that $\lvert \pi \rvert \leq n_a +4$. In this case, we show that $\lvert \pi \rvert \neq n_a +4$. Suppose to the contrary that $|\pi|=n_a+4$. Then we must have $C \cup V_a = L_b \cup L_c \cup \lbrace a\rbrace$, and each set in $\pi \setminus \lbrace C, V_a \rbrace$ must be a singleton set. The singleton $\{b\}$ has no $gc$-partner among $\{c\}$, the leaves and $C$ (the corresponding unions do not dominate $L_a$ or $L_c$), so $\{b\}$ is a $gc$-partner of $V_a$, which forces $L_c\subseteq V_a$. Similarly $\{c\}$ is a $gc$-partner of $V_a$, which forces $L_b\subseteq V_a$. Thus $V_a=L_b\cup L_c\cup\{a\}$ and $C=\emptyset$, a contradiction. Hence,   $\lvert \pi \rvert \leq n_a +3$. A similar argument shows that if $C$ is a $gc$-partner of $V_b$, then $\lvert \pi \rvert \leq n_b +3$, and if $C$ is a $gc$-partner of $V_c$, then $\lvert \pi \rvert \leq n_c +3$. Hence,  $\lvert \pi \rvert \leq n_a +3$, and so $\lvert \pi \rvert \leq n_a + n_b +2$.

	\textbf{Case 2.} There exists a member $A \in \pi$ such that $\lvert A \cap \lbrace a,b,c\rbrace \rvert =2$.

	Assume that $\lbrace a,b\rbrace \subseteq A$. Let $V_c$ be the member of $\pi$ containing $c$. If there exists a pair $(X,Y)$ of $gc$-partners such that $\lbrace X,Y\rbrace \cap \lbrace A,V_c\rbrace =\emptyset$, then we must have $X \cup Y =L_a \cup L_b \cup L_c$, implying that $\lvert \pi \rvert =4$, as desired. Hence, we may assume that each set in $\pi \setminus \lbrace A, V_c\rbrace$ has a $gc$-partner in $\lbrace A, V_c\rbrace$. Consider an arbitrary member $B \in \pi \setminus \lbrace A, V_c\rbrace$. Assume first that $B$ is a $gc$-partner of $A$. Then we must have $L_c \subset A \cup B$, implying that $\lvert \pi \rvert \leq n_a +n_b + 3$. Now we show that $\lvert \pi \rvert \neq n_a +n_b + 3$. Suppose, to the contrary, that equality holds. Then we must have $A \cup B =L_c \cup \lbrace a,b\rbrace$, and each set in $\pi \setminus \lbrace A, B\rbrace$ must be a singleton set. Let $e \in L_a$. Now it is easy to verify that the set $\lbrace e\rbrace$ has no $gc$-partner: the union of $\{e\}$ with $\{c\}$ or with a leaf misses $b$ or both $b$ and $c$ and therefore does not dominate $L_b\cup L_c$ (note that $|L_b\cup L_c|\ge 2$), the set $\{e\}\cup B\subseteq L_a\cup L_c$ does not dominate $L_b$, and $\{e\}\cup A$ does not dominate $L_c\setminus A\ne\emptyset$, a contradiction. Hence, $\lvert \pi \rvert \leq n_a +n_b + 2$, as desired. Otherwise, we may assume that $B$ is a $gc$-partner of $V_c$. It follows that $(L_a \cup L_b) \subset (B \cup V_c)$, implying that $\lvert \pi \rvert \leq n_c +3$, and so $\lvert \pi \rvert \leq n_a + n_b +2$, as desired. A similar argument shows that if $\lbrace a,c\rbrace \subseteq A$ or $\lbrace b,c\rbrace \subseteq A$, then $\lvert \pi \rvert \leq n_a +n_b +2$. This completes the proof.
\end{proof}
	\begin{proposition}\label{prop-u3-2}
		If $G \in  \left(  \bigcup_{i=4} ^{9} U_{3_i} \right) \cup \left(  \bigcup_{i=11} ^{13} U_{3_i} \right) \cup \left(  \bigcup_{i=15} ^{17} U_{3_i} \right) $, then $GC(G) \geq n-\lvert N(a) \rvert +1$, where $a$ is the central vertex of $G$ (the vertex denoted by $a$ in Figure~\ref{u3-rad2}).
	\end{proposition}
	\begin{proof}
	Consider the vertex partition $\pi$ of $G$, where $N(a) \in \pi$ and all other partition classes are singleton sets. Then $\pi$ is a $gc$-partition of $G$ of size $n- \lvert N(a) \rvert +1$, where the set $N(a)$ forms a global coalition with all other partition classes (Lemma~\ref{lem-uni-lb}). Hence,  $GC(G) \geq n- \lvert N(a) \rvert +1$.
\end{proof}

One could conjecture that the bounds of Propositions~\ref{prop-u5-2}, \ref{prop-u4-3} and~\ref{prop-u3-2} hold with equality. This is not the case in general. Consider the graph $G\in U_{3_4}$ in which $a$ has exactly one strong child $s$, and $s$ has exactly two leaves $x_1,x_2$ (Figure~\ref{fig-counter}, left). Here $n=6$ and $n-|N(a)|+1=4$, but $\{\{a\},\{b\},\{c\},\{s\},\{x_1,x_2\}\}$ is a $gc$-partition: $\{a\}$ and $\{s\}$ are $gc$-partners, and $\{x_1,x_2\}$ is a $gc$-partner of both $\{b\}$ and $\{c\}$. Hence $GC(G)=5$. The difference can be arbitrarily large: if, in addition, $a$ is adjacent to $k$ leaves $y_1,\dots,y_k$ (so $G\in U_{3_5}$), then $n-|N(a)|+1=4$, while $\{a,x_1,x_2\}$ together with the singletons of the remaining $k+4$ vertices is a $gc$-partition, so $GC(G)\ge k+5$. Note that in this example $s$ is also a central vertex and $\{a,x_1,x_2\}=N(s)$, so the lower bound is again of the form given by Lemma~\ref{lem-uni-lb}, but for a different central vertex. A similar situation occurs in $U_{4_5}$ when every strong child of $a$ has exactly one leaf (Figure~\ref{fig-counter}, right): the set consisting of $c$ and all leaves at distance two from $a$ is a $gc$-partner of every other singleton, which gives $GC(G)\ge n-k$ for $k$ strong children, whereas $n-|N(a)|+1=n-k-1$. In Section~\ref{sec-concl}, we formulate a corrected conjecture.

\begin{figure}[!htbp]
	\centering
	\begin{tikzpicture}[scale=0.9, every node/.style={font=\small}]
		\tikzstyle{v}=[draw, circle, inner sep=0pt, minimum size=5pt]
		\node[v, fill=red]           (a)  at (0,0) {};
		\node[v, fill=blue]          (b)  at (-1.1,-1.5) {};
		\node[v, fill=green!60!black](c)  at (1.1,-1.5) {};
		\node[v, fill=orange]        (s)  at (0,1.5) {};
		\node[v, fill=violet]        (x1) at (-0.9,2.7) {};
		\node[v, fill=violet]        (x2) at (0.9,2.7) {};
		\draw (a)--(b)--(c)--(a); \draw (a)--(s); \draw (s)--(x1); \draw (s)--(x2);
		\node at (0.35,-0.1) {$a$}; \node at (-1.45,-1.5) {$b$}; \node at (1.45,-1.5) {$c$};
		\node at (0.35,1.5) {$s$}; \node at (-1.25,2.75) {$x_1$}; \node at (1.25,2.75) {$x_2$};
		\begin{scope}[xshift=6cm]
		\node[v, fill=red]           (a)  at (0,0) {};
		\node[v, fill=blue]          (b)  at (-1.3,-1.2) {};
		\node[v, fill=violet]        (c)  at (0,-2.4) {};
		\node[v, fill=orange]        (d)  at (1.3,-1.2) {};
		\node[v, fill=green!60!black](s1) at (-1,1.3) {};
		\node[v, fill=cyan!70!black] (s2) at (1,1.3) {};
		\node[v, fill=violet]        (x1) at (-1,2.6) {};
		\node[v, fill=violet]        (x2) at (1,2.6) {};
		\draw (a)--(b)--(c)--(d)--(a); \draw (a)--(s1); \draw (a)--(s2); \draw (s1)--(x1); \draw (s2)--(x2);
		\node at (0.35,0.05) {$a$}; \node at (-1.65,-1.2) {$b$}; \node at (0.35,-2.45) {$c$}; \node at (1.65,-1.2) {$d$};
		\node at (-1.4,1.3) {$s_1$}; \node at (1.4,1.3) {$s_2$}; \node at (-1.4,2.65) {$x_1$}; \node at (1.4,2.65) {$x_2$};
		\end{scope}
	\end{tikzpicture}

	\caption{Two unicyclic graphs for which $GC(G)>n-|N(a)|+1$. The classes of a $gc$-partition are indicated by colors: left, a graph in $U_{3_4}$ with $GC(G)=5>4$, right, a graph in $U_{4_5}$ with $GC(G)=6>5$.}
	\label{fig-counter}
\end{figure}

\section{A connection with fair division of indivisible goods}\label{sec-fair}
In this section, we show that $gc$-partitions admit a natural interpretation in the language of fair division of indivisible goods. For some background on fair division, see e.g. the survey \cite{amanatidis}.

In the basic model, a set $M$ of indivisible goods is to be allocated to $k$ agents. An allocation is a partition $(A_1,\dots,A_k)$ of $M$ into bundles, and agent $i$ evaluates a bundle $S$ by a monotone valuation $v_i(S)$. One of the most studied fairness criteria is \emph{envy-freeness up to any good} (EFX), which requires $v_i(A_i)\ge v_i(A_j\setminus\{g\})$ for all agents $i,j$ and every good $g\in A_j$, being a relaxation of envy-freeness (EF), where one needs $v_i(A_i)\ge v_i(A_j)$. The other one \emph{maximin share guarantee} (MMS), which requires $v_i(A_i)\ge \mathrm{MMS}^k_i$, where $\mathrm{MMS}^k_i$ is the maximum, over all partitions $(P_1,\dots,P_k)$ of $M$, of $\min_j v_i(P_j)$. Valuations taking only the values $0$ and $1$ (\emph{binary} or \emph{dichotomous} valuations) form an important special case \cite{babaioff,halpern}, and allocations of the vertex set of a graph subject to graph-theoretic constraints have been studied, e.g., in \cite{bouveret}. On the other hand, the variants with goods assigned to the edges with some additional constraints imposed, like acyclicity or connectedness, were recently studied, e.g. by Anholcer et al. \cite{anholcer-epistemic,anholcer-rainbow}.

Let $G$ be a graph of order $n\ge 2$, let the goods be the vertices of $G$, and consider $k$ agents sharing the binary valuation
\[
v_G(S)=\begin{cases}1, & \text{if $S$ is a global dominating set of $G$},\\ 0, & \text{otherwise.}\end{cases}
\]
This setting may describe a situation in which two complementary networks connect a set of objects. Every agent wishes to control all the objects, where control means being able to reach each of them in both networks, $G$ and $\overline{G}$. The problem with a single network corresponds to ordinary coalitions \cite{haynes3}, while the setting with several networks on the same vertex set can be described by factor domination \cite{brigham}. The valuation $v_G$ is monotone and binary, so the following holds.

\begin{enumerate}[(i)]
\item
$\mathrm{MMS}^k(v_G)=1$ if and only if $V(G)$ can be partitioned into $k$ global dominating sets, i.e., if and only if $k\le d_g(G)$. Thus, the global domatic number of Zelinka and Sampathkumar is precisely the largest number of agents each of whom can be guaranteed their maximin share with value one. For $k\le d_g(G)$, any global domatic partition with $k$ classes (obtained from a maximum one by combining the classes) is an allocation which is simultaneously MMS and EF, thus also EFX.

\item
For $k>d_g(G)$ we have $\mathrm{MMS}^k(v_G)=0$, so every allocation satisfies the MMS criterion. Moreover, every allocation in which no agent receives a global dominating set is EF (and thus EFX), because all agents value all bundles and all sub-bundles at $0$. In other words, if the number of agents exceeds $d_g(G)$, the criteria under consideration do not distinguish between allocations.
\end{enumerate}

Global coalition partitions change the situation in case (ii). Let us say that an allocation $(A_1,\dots,A_k)$ with $v_G(A_i)=0$ for all $i$ is \emph{completable by an ally} if for every agent $i$ there is an agent $j\ne i$ with $v_G(A_i\cup A_j)=1$. In other words, none of the agents can be satisfied individually, but everyone can find another agent such that their combined bundles will be sufficient. This means that the $gc$-partitions of $G$ are exactly the allocations completable by an ally, and $GC(G)$ is the largest number of agents for which the vertices of $G$ can be divided so that the division is not fair for any agent individually, but everyone can find an ally and combine the bundles with them to become satisfied.

This means that the results of the previous sections can be interpreted as follows. Theorem~\ref{gc-bound} provides an upper bound for the number of agents that could be possibly allies of a given agent. Corollary~\ref{cor-2dg} says that if there is a partition that is fair for $k$ agents, then there exists also an ally-completable partition for $2k$. Proposition~\ref{prop-mu} says that to find an allocation completable by an ally, it is enough to assign a maximal unsatisfying bundle to one agent and a single good to every other agent. Theorem~\ref{rad-c-gc} shows that for networks of radius at least three, it does not matter whether one considers only one network or two networks simultaneously. Theorem~\ref{thm-antitwin} says that every good can be assigned to a different agent if and only if one of the vertices corresponding to the goods has an antitwin (in both networks).

The quantity $\min_i\max_{j\ne i}v(A_i\cup A_j)$ is a relaxation of the $\ell$-out-of-$d$ maximin shares of \cite{BNT}, where an agent is guaranteed the union of several bundles. The difference is that, here, the second bundle is chosen based on the agent's preferences rather than competitively. We do not investigate this direction further, but we would like to mention two questions. First, what can be said when the agents have different binary valuations, for instance, when agent $i$ is satisfied by the dominating sets of its own network $G_i$, and the networks of all the agents share the vertex set? In particular, the case $G_1=G$, $G_2=\overline{G}$ corresponds to global domination, as considered in the previous sections of the paper. Another possible research direction is as follows. The $gc$-partitions constructed in Proposition~\ref{prop-mu} and Lemma~\ref{lem-uni-lb} are very unbalanced: one bundle is large, and every other agent receives precisely one good. It would be interesting to study more balanced $gc$-partitions, where the classes have (almost) equal sizes, a constraint sometimes imposed in fair division problems.

	\section{Concluding remarks and open problems}\label{sec-concl}
	We introduced global coalitions and global coalition partitions, proved that every graph of order at least two admits such a partition, and determined the global coalition number for several classes of graphs. The main structural result of the paper is Theorem~\ref{rad-c-gc} together with Corollary~\ref{cor-join}: the new parameter differs from the coalition number of $G$ and of $\overline{G}$ only when both $G$ and $\overline{G}$ are connected and have radius at most two. We conclude the paper with the following problems for future research. 
	\begin{problem}
		Characterise the graphs $G$ for which there exists a $gc$-partition $\pi$ of $G$ and a class $A\in\pi$ whose number of $gc$-partners satisfies the bound of Theorem~\ref{gc-bound} with equality.
	\end{problem}

	Proposition~\ref{prop-multi} and Theorem~\ref{thm-tripartite} suggest the following. 
	
	\begin{conjecture}\label{conj-multi}
	For every complete multipartite graph $G=K_{n_1,\dots,n_m}$ with $n_1\ge\dots\ge n_m$ and $m\ge 2$, the global coalition number $GC(G)$ equals the lower bound of Proposition~\ref{prop-multi}.
	\end{conjecture} 
	
	Equivalently, by Corollary~\ref{cor-join}, Conjecture~\ref{conj-multi} determines the coalition number of every disjoint union of complete graphs. 
	
     \begin{problem}
     \label{prob-uni}In Propositions \ref{prop-u5-2}, \ref{prop-u4-3} and \ref{prop-u3-2} we showed that if
	$$ G \in \left(\bigcup_{i=5} ^{10} U_{5_i}\right) \cup \left(\bigcup_{i=4} ^{10} U_{4_i}\right) \cup \left(  \bigcup_{i=4} ^{9} U_{3_i} \right) \cup \left(  \bigcup_{i=11} ^{13} U_{3_i} \right) \cup \left(  \bigcup_{i=15} ^{17} U_{3_i} \right),$$
then $GC(G) \geq n-\lvert N(a) \rvert +1,$  where $a$ denotes  the central vertex of $G$.

One might suppose that equality always holds, which is false, as shown at the end of Section~\ref{sec-gl-uni}: in the families $U_{3_4}$, $U_{3_5}$, $U_{3_6}$, $U_{4_5}$ and $U_{4_6}$ a different vertex of eccentricity two, or a maximal non-global-dominating set that is not an open neighborhood, may give a better bound. Determine $GC(G)$ for the graphs in these five families.
\end{problem}

\noindent{}Let us formulate another conjecture.

     \begin{conjecture}\label{conj-uni}
     If $G$ belongs to one of the families $U_{5_5}$--$U_{5_{10}}$, $U_{4_4}$, $U_{4_7}$--$U_{4_{10}}$, $U_{3_7}$--$U_{3_9}$, $U_{3_{11}}$--$U_{3_{13}}$ or $U_{3_{15}}$--$U_{3_{17}}$, then $GC(G)=n-|N(a)|+1$, where $a$ is the central vertex of $G$ indicated in Figures~\ref{u5-rad2}--\ref{u3-rad2}.
     \end{conjecture}
     We have no structural explanation of the different behaviour of the families in Problem~\ref{prob-uni} and in Conjecture~\ref{conj-uni}. Finding one would probably be the first step towards a proof.
 
	\begin{problem}
		Investigate the behaviour of the global coalition number under graph operations, such as Cartesian and lexicographic products. 
	\end{problem}
	
By Corollary~\ref{cor-join}, for the join, this amounts to determining the coalition number of disjoint unions of graphs.

\begin{problem}
For any graph $G$ and any $gc$-partition $\pi$ of $G$, we can define the global coalition graph $GCG(G,\pi)$, whose vertices correspond
		bijectively to the elements of $\pi$. Two vertices are adjacent in $GCG(G,\pi)$ if and only if their corresponding subsets in $\pi$ together form a global coalition in $G$. This leads to the following question: Is every graph realisable as a global coalition graph?
	\end{problem}
	
For ordinary coalitions, this question was answered affirmatively in \cite{haynes1}, and self-coalition graphs were studied in \cite{haynes4}. By Theorem~\ref{rad-c-gc}, every coalition graph of a graph of radius at least three is also a global coalition graph.

Let us conclude this section with a computational problem.

	\begin{problem}
Determine the computational complexity of deciding whether $GC(G)\ge k$ for a given graph $G$ and integer $k$. By Theorem~\ref{rad-c-gc}, for graphs of radius at least three, this is equivalent to the corresponding question for the coalition number.
	\end{problem}
	
	\begin{flushleft}
		\textbf{{\large Conflicts of interest}}\vspace{-3.5mm}
	\end{flushleft}
	The authors declare no conflicts of interest.
	\begin{flushleft}
		\textbf{{\large Data availability}}\vspace{-3.5mm}
	\end{flushleft}
	Not applicable -- no data were produced during the research.
	\begin{flushleft}
		\textbf{{\large Disclosure}}\vspace{-3.5mm}
	\end{flushleft}
	The research was initiated, and the initial version of this paper was prepared by Nazli Besharati, Doost Ali Mojdeh and Mohammad Reza Samadzadeh (\url{https://arxiv.org/abs/2509.15386v1}). After Marcin Anholcer joined the team, several proofs were reformulated, and a few results were added. A completely new part is Section \ref{sec-fair}. Various AI tools were used to aid the bibliography search, check the consistency of notions and notation, and correct language mistakes in the final version of this paper. The authors take full responsibility for the contents of the paper.


\baselineskip12pt
	\begin{thebibliography}{99}

\bibitem{alikhani} S. Alikhani, D. Bakhshesh, H.R. Golmohammadi and E.V. Konstantinova, Connected coalitions in graphs, Discuss. Math. Graph Theory 44 (2024), no. 4, 1551--1566.

\bibitem{ABG-total} S. Alikhani, D. Bakhshesh and H.R. Golmohammadi, Total coalitions in graphs, Quaest. Math. 47 (2024), no. 11, 2283--2294.

\bibitem{ABGK} S. Alikhani, D. Bakhshesh, H.R. Golmohammadi and S. Klav\v{z}ar, On independent coalition in graphs and independent coalition graphs, Discuss. Math. Graph Theory 45 (2025), no. 2, 533--544.

\bibitem{amanatidis} G. Amanatidis, H. Aziz, G. Birmpas, A. Filos-Ratsikas, B. Li, H. Moulin, A.A. Voudouris and X. Wu, Fair division of indivisible goods: Recent progress and open questions, Artificial Intelligence 322 (2023), 103965.

\bibitem{anholcer-epistemic}
M. Anholcer, M. Bartkowiak, B. Bosek and J. Grytczuk, Epistemic fair division of independence structures,
\href{https://arxiv.org/abs/2606.11494}{arXiv:2606.11494 [math.CO]} (2026).

\bibitem{anholcer-rainbow}
M. Anholcer, M. Bartkowiak, B. Bosek, J. Grytczuk, Z. Lonc and P. Rz\k{a}\.zewski, A rainbow partition theorem for trees and connected maximin share allocations of chores,
\href{https://arxiv.org/abs/2609.38628}{arXiv:2609.38628 [math.CO]} (2026).


\bibitem{babaioff} M. Babaioff, T. Ezra and U. Feige, Fair and truthful mechanisms for dichotomous valuations, in: Proceedings of the 35th AAAI Conference on Artificial Intelligence (AAAI 2021), 5119--5126.

\bibitem{BNT} M. Babaioff, N. Nisan and I. Talgam-Cohen, Competitive equilibrium with indivisible goods and generic budgets, Math. Oper. Res. 46 (2021), no. 1, 382--403.

\bibitem{BHP-trees} D. Bakhshesh, M.A. Henning and D. Pradhan, On the coalition number of trees, Bull. Malays. Math. Sci. Soc. 46 (2023), no. 3, Art. 95.

\bibitem{BHP-chains} D. Bakhshesh, M.A. Henning and D. Pradhan, Singleton coalition graph chains, Comput. Appl. Math. 43 (2024), no. 2, Art. 85.

\bibitem{barat} J. Bar\'at and Z.L. Bl\'azsik, General sharp upper bounds on the total coalition number, Discuss. Math. Graph Theory 44 (2024), no. 4, 1567--1584.

\bibitem{Bondy}
J.A. Bondy and  U.S.R. Murty, Graph Theory, Springer, 2008.

\bibitem{bouveret} S. Bouveret, K. Cechl\'arov\'a, E. Elkind, A. Igarashi and D. Peters, Fair division of a graph, in: Proceedings of the 26th International Joint Conference on Artificial Intelligence (IJCAI 2017), 135--141.

\bibitem{BCD-survey} R.C. Brigham, J.R. Carrington and R.D. Dutton, Global domination, in: T.W. Haynes, S.T. Hedetniemi and M.A. Henning (eds.), Topics in Domination in Graphs, Developments in Mathematics 64, Springer, Cham, 2020, 497--519.

\bibitem{brigham}
R.C. Brigham and R.D. Dutton, Factor domination in graphs, Discrete Math. 86 (1990), no. 1-3, 127--136.

\bibitem{chartrand}
G. Chartrand, H. Jordon, V. Vatter and P. Zhang, Graphs \& Digraphs, 7th ed., Chapman and Hall/CRC, Boca Raton, 2024.

\bibitem{desormeaux} W.J. Desormeaux, P.E. Gibson and T.W. Haynes, Bounds on the global domination number, Quaest. Math. 38 (2015), no. 4, 563--572.

\bibitem{halpern} D. Halpern, A.D. Procaccia, A. Psomas and N. Shah, Fair division with binary valuations: One rule to rule them all, in: Web and Internet Economics (WINE 2020), Lecture Notes in Computer Science 12495, Springer, Cham, 2020, 370--383.

\bibitem{haynes1}
T.W. Haynes, J.T. Hedetniemi, S.T. Hedetniemi, A.A. McRae and R. Mohan, Coalition graphs, Commun. Comb. Optim. 8 (2023), no. 2, 423--430.

\bibitem{haynes2}
T.W. Haynes, J.T. Hedetniemi, S.T. Hedetniemi, A.A. McRae and R. Mohan, Coalition graphs of paths, cycles, and trees, Discuss. Math. Graph Theory 43 (2023), no. 4, 931--946.

\bibitem{haynes3}
T.W. Haynes, J.T. Hedetniemi, S.T. Hedetniemi, A.A. McRae and R. Mohan, Introduction to coalitions in graphs, AKCE Int. J. Graphs Comb. 17 (2020), no. 2, 653--659.

\bibitem{haynes4}
T.W. Haynes, J.T. Hedetniemi, S.T. Hedetniemi, A.A. McRae and R. Mohan, Self-coalition graphs, Opuscula Math. 43 (2023), no. 2, 173--183.

\bibitem{haynes-ub} T.W. Haynes, J.T. Hedetniemi, S.T. Hedetniemi, A.A. McRae and R. Mohan, Upper bounds on the coalition number, Australas. J. Combin. 80 (2021), no. 3, 442--453.

\bibitem{haynes-hs}
T.W. Haynes, S.T. Hedetniemi and P.J. Slater, Fundamentals of Domination in Graphs, Marcel Dekker, New York, 1998.

\bibitem{HHH-core} T.W. Haynes, S.T. Hedetniemi and M.A. Henning, Domination in Graphs: Core Concepts, Springer Monographs in Mathematics, Springer, Cham, 2023.

\bibitem{HeMo-24} M.A. Henning and D.A. Mojdeh, Double coalitions in graphs, Bull. Malays. Math. Sci. Soc.  48 (2025), no. 2, Art. 51.

\bibitem{HeMo-25} M.A. Henning and D.A. Mojdeh, Double coalitions in regular graphs, Graphs Combin.  41 (2025), no. 4, Art. 74.

\bibitem{jafari} A. Jafari, S. Alikhani and D. Bakhshesh, $k$-coalitions in graphs, Australas. J. Combin. 92 (2025), no. 2, 194--209.

\bibitem{mojdeh-ac}
D.A. Mojdeh, M. Alishahi and M. Chellali, Trees with the same global domination number as their square, Australas. J. Combin. 66 (2016), no. 2, 288--309.

\bibitem{mojdeh-s}
D.A. Mojdeh and M.R. Samadzadeh, Perfect coalition in graphs, Electron. J. Graph Theory Appl. 13 (2025), no. 2, 343--360.

\bibitem{olariu} S. Olariu, No antitwins in minimal imperfect graphs, J. Combin. Theory Ser. B 45 (1988), no. 2, 255--257.

\bibitem{samadzade-m}
M.R. Samadzadeh and D.A. Mojdeh, Independent coalition in graphs: existence and characterization,  Ars Math. Contemp. 24 (2024), no. 3, \#P3.09.


\bibitem{samadzade-mn}
M.R. Samadzadeh, D.A. Mojdeh and R. Nadimi, Paired coalition in graphs, AKCE Int. J. Graphs Comb. 22 (2025), no. 1, 43--54.


\bibitem{sampathkumar}
E. Sampathkumar, The global domination number of a graph, J. Math. Phys. Sci. 23 (1989), no. 5, 377--385.

\bibitem{zelinka} B. Zelinka, Complementarily domatic number of a graph, Math. Slovaca 38 (1988), no. 1, 27--32.

\end{thebibliography}
\end{document}